\documentclass[a4paper,final,reqno]{report}
\usepackage{euscript,amsfonts,amsmath,amsthm}
\usepackage{hhline,array}
\usepackage{verbatim}
\newtheorem{theorem}{Theorem}[section]
\newtheorem{corollary}[theorem]{Corollary}
\newtheorem{conjecture}[theorem]{Conjecture}
\newtheorem{lemma}[theorem]{Lemma}

\theoremstyle{definition}
\newtheorem{definition}[theorem]{Definition}

\theoremstyle{remark}
\newtheorem{remark}[theorem]{Remark}

\theoremstyle{remark}

\theoremstyle{remark}

\theoremstyle{plain}

\newtheoremstyle{note}
  {3pt}
  {3pt}
  {}
  {}
  {\itshape}
  {:}
  {.5em}
  {}

\theoremstyle{note}

\newtheoremstyle{citing}
  {3pt}
  {3pt}
  {\itshape}
  {}
  {\bfseries}
  {.}
  {.5em}
  {\thmnote{#3}}

\theoremstyle{citing}

\newtheoremstyle{break}
  {9pt}
  {9pt}
  {\itshape}
  {}
  {\bfseries}
  {.}
  {\newline}
  {}

\theoremstyle{break}

\swapnumbers
\theoremstyle{plain}

\newcommand{\raz}{//}                  

\let\lvert=|\let\rvert=|

\chardef\bslash=`\\ 

\newlength{\halfspatsiya}
\newcommand{\oequal}%
{\mbox{\textcircled{\mbox{\hspace{0.8\halfspatsiya}= \;}}}}

\newcommand{\eval}[2][\right]{\relax
  \ifx#1\right\relax \left.\fi#2#1\rvert}

\newcommand{\Io}{I.Jo\'o{}}

\newcommand{\bysame}{\mbox{\rule{3em}{.4pt}}\,}

\newcommand{\no}{$no^0$}
\newcommand{\bvp}{b.v.p.\ }
\newcommand{\ef}{e.f.\ }
\newcommand{\eaf}{e.a.f.\ }
\newcommand{\ev}{e.v.\ }
\newcommand{\chv}{ch.v.\ }
\newcommand{\fss}{f.s.s.\ }
\newcommand{\spf}{sp.f.\ }
\newcommand{\gspf}{g.sp.f.\ }
\newcommand{\sa}{s.a.\ }

\makeindex

\makeatletter
\newcommand{\INDEX}[2]{\index{#1@#2}}
\@addtoreset{equation}{section}
\renewcommand{\p@equation}{}

\newcommand{\l@abcd}[2]{\hbox to\textwidth{#1\dotfill #2}}
\makeatother

\begin{document}

\begin{titlepage}

\thispagestyle{empty}


\begin{center}
\renewcommand{\thefootnote}{\fnsymbol{footnote}}
{\large A.M.Minkin
\footnote{
Partially supported by the Russian fund of fundamental researches,
project N 97-01-00566,
and by the International Soros Science Education Program,
grant N 149d "Soros Associate Professor"}
}\\
\setcounter{footnote}{0}

\vspace*{6mm}
{\large\textbf{
Equiconvergence theorems for differential operators}}\\
\vspace*{6mm}
\end{center}
\end{titlepage}

\tableofcontents

\date{}

\begin{abstract}
The paper is a survey dedicated to topic in the title.
In chapter 1 we expose the most advanced equiconvergence results for Birkhoff- or Stone-regular differential operators.
Considerable part of them was obtained by the Saratov mathematical school 
but is published in the literature that is hard to come by.

We present also an author's (commutator) approach to equiconvergence
and derive a stonger form of the Riemann localization principle
as well as a first equiconvergence result (not equisummability)
for multidimensional Schr\"odinger operator.

Chapter 2 contains a full proof of the equiconvergence on the whole interval, which constitutes a true generalization of the Tamarkin-Stone theorem. 

Given Birkhoff-regular ordinary differential operator $L$ in $L^2(0,1)$ and continuous function $f$, which belongs to closure of $D_L$ in $C[0,1]$,
we establish necessary and sufficient
conditions for uniform equiconvergence on $[0,1]$ of the 
eigenfunction expansion of $L$ and of trigonometric Fourier integral of the modified function 
\[
\tilde{f}(x)=\left\{
\begin{array}{ll}
		    f(0), & x<0;\\
                    f(x), & 0\le x \le1;\\
		    f(1), & x>1.
\end{array}
\right.
\]
These conditions consist of uniform converge to $0$ on $[0,\delta]$ for some (any) $0 < \delta <1$
of certain singular integrals acting upon specified linear combinations of functions $f(x)$ and $f(1-x)$.

Chapters 3 and 4 apply our approach to equiconvergence 
to singular self-adjoint differential operators, generalizing
well-known results of A.G.Kostuchenko, and to general serii in eigenfunctions.
\end{abstract}

\centerline{\textbf\large{Preface}}

Many central problems of the spectral theory of linear operators
are concentrated
around the problem of eigenfunction expansions. From one hand it
accumulates questions of eigenvalues and eigenfunctions asymptotics,
from the other
it connects mathematics with many physical problems of string
and membrane vibrations, of quantum mechanics and so on.
There are two most elaborated parts of this theory,
firstly, the
spectral theory of self-adjoint operators (ordinary, singular and
in partial derivatives)
and secondly,
that of boundary value
problems in a finite interval.

Herewith already at the beginning of the 20th century
G.D.Birkhoff
\INDEX{Birkhoff}{G.D.Birkhoff}
discovered an important class of \textit{regular}
higher order boundary value problems.
 Just immediately J.Tamarkin
\INDEX{Tam}{J.Tamarkin}
observed that for two such
problems the difference of eigenfunction expansions converges to
zero in any interior point of the main interval. This phenomenon
was called \textit{equiconvergence} and it makes possible
to reduce numerous questions of point and uniform convergence to
those of some model, usually, trigonometric system.

This remarquable result has several predecessors in the case of
second order operators. We shall mention
\INDEX{Ste}{V.A.Steklov}V.A.Steklov,
\INDEX{Hob}{E.W.Hobson}E.W.Hobson
and
\INDEX{Haa}{A.Haar}A.Haar.
Later on it was generalized by
\INDEX{Tam}{J.Tamarkin}J.Tamarkin
and by
\INDEX{Sto}{M.Stone}M.Stone
and yielded a large field of investigations
with a lot of off-shoots and generalizations. In the present review we set
ourselves a task of exposing the current state-of-arts in this domain
with a strong
emphasis to describe the most advanced achievements,
at least to the best of our knowledge.
During the past two decades the author himself developped new
approaches to these questions and obtained solutions of several long
standing problems. However, they are published in a literature which is hard
to come by. Partially they remain yet unpublished or appeared only in
conferences proceedings. Therefore we place here the most important of them
with complete demonstrations. In particular, a solution is given to the
problem of \emph{equiconvergence on the whole interval} obtained in 1992.
We also present a complete investigation of higher order
singular self-adjoint quasidifferential operators under minimal possible
restrictions.

History of the question is exposed but with no claim of
completeness.
Hence, the bibliography is rather large but not
exhaustive.

During the exposition we formulate several conjectures reflecting
our own
understanding of the subject in order to stimulate further
investigations whether these conjectures will happen to be
true or not.

Due to the lack of time as well as author's insufficient
knowledge we omit some important topics,
for instance, opearator bundles,
differential expressions with multiple roots of the
characteristic equation and all the more with varying
multiplicity roots (equations with a turning point), operators in
partial derivatives. We only touch the latter once.
These questions deserve separate review or reviews.

During the preparation we have fruitfull discussions
with our colleagues. They also provided us many materials used in
the paper. In particular,
Professor A.P.Khromov
\INDEX{Khro}{A.P.Khromov}
has kindly given us permission to use his
unpublished notes on
the problem of equiconvergence.

 Therefore we take an opportunity to thank
A.P.Khromov, \INDEX{Khro}{A.P.Khromov}
B.\'E.Ku\-nyav\-skii, \INDEX{Kun}{B.\'E.Kunyavskii}
S.N.Kup\-t\-sov, \INDEX{KupS}{S.N.Kuptsov}
G.V.Radzievskii, \INDEX{Radz}{G.V.Radzievskii}
V.S.Rykhlov \INDEX{Rykh}{V.S.Rykhlov}
and I.Yu.Trushin \INDEX{Tru}{I.Yu.Trushin}
though the list may have been
considerably increased.

Throughout the paper the reader will often meet citations of
N.P.Kuptsov's \INDEX{KupN}{N.P.Kuptsov}
[1925-1995] results.
However, equiconvergence reflects only one of the numerous fields
of interest and activity of this universal mathematician.
He published very selectively but his impact
on the developement of the Saratov school on spectral theory is difficult
to overestimate.
Therefore we dared to dedicate him this review as a small tribute
to his memory.

In the body of the text we use some notations which has become standard.
\begin{itemize}
\item b.v.p.  --- boundary value problems,
\item e.f.    --- eigenfunctions,
\item e.a.f.  --- eigenfunctions and associated eigenfunctions,
\item e.v.    --- eigenvalues,
\item ch.v.   --- characteristic values,
\item f.s.s.  --- fundamental system of solutions,
\item sp.f.   --- spectral function,
\item g.sp.f. --- generalized spectral function,
\item s.a.    --- self-adjoint,
\item span    --- minimal closed subspace, containing a given set of elements,
\item $[a]:=a+O(1/\rho)$ --- the Birkhoff's symbol. \INDEX{Birkhoff}{G.D.Birkhoff}
\item $Entier(h)$  ---  the largest integer $\leq h$.
\end{itemize}


\chapter{Introduction}\label{Chap.intro}

\textit{\hfill To the memory of N.P.Kuptsov}
\INDEX{KupN}{N.P.Kuptsov}

\section{Birkhoff-regular problems}
\INDEX{Birkhoff}{G.D.Birkhoff}

\subsection{Early results}

Let us consider a differential operator $L$ in $L^2(0,1)$ defined by a
two-point \bvp ($(D=-id/dx)$):%
\hfill
\begin{equation}
l(y)\equiv D^ny+\sum_{k=0}^{n-2}p_k(x)D^ky=\lambda y,\ \ 0\le x\le 1,\ \
p_k\in L(0,1)  \label{eqB.1}
\end{equation}
and $n$ linearly independent normalized boundary conditions \cite[p.65--66]
{Nai:ldo}:
\begin{gather}
U_j(y)\equiv V_j(y)+\ldots =0,\;\;j=0,\ldots ,n-1,  \label{eqB.2} \\
V_j(y)\equiv b_j^0D^jy(0)+b_j^1D^jy^{(}1).  \notag
\end{gather}
\noindent
Here the ellipsis takes place of lower order terms at $0$ and at $1$.
Further $b_j^0,b_j^1$ are column vectors of length $r_j$, where
\begin{equation*}
0\le r_j\le 2,\;\sum_{k=0}^{n-1}r_j=n,\ \ rank[b_j^0b_j^1]=r_j.
\end{equation*}
This form of normalized boundary conditions was first introduced by
\INDEX{Salaff}{S.Salaff}S.Salaff \cite[p.356--357]{sal68}.
It is evident that $r_j=0$ implies the
absence of order $j$ conditions. In the case $r_j=2$ we merely put
\begin{equation*}
\left[ b_j^0b_j^1\right] =\left[
\begin{array}{cc}
1 & 0 \\
0 & 1
\end{array}
\right] .
\end{equation*}

At the beginning of the 20th century
\INDEX{Birkhoff}{G.D.Birkhoff}G.D.Birkhoff discovered a famous broad class
of \bvp  with remarkable spectral properties \cite{bir08a,bir08b}. Recall
its definition.
\begin{definition}
Let $q=Entier(n/2)$,
$\varepsilon _j:=\exp (2\pi ij/n)$, $k=0,..,n-1;$
\begin{equation}
b^i=(b_j^i)_{j=0}^{n-1}=\left[
\begin{array}{l}
{b_0^i} \\
\vdots \\
b_{n-1}^i
\end{array}
\right] ,\;\mathbf{B}_k^i=[b_j^i\cdot \varepsilon _k^j]_{j=0}^{n-1},\;i=0,1;
\label{eqB.3}
\end{equation}
\medskip
\begin{equation}
\theta (b^0,b^1,L)=\det [\mathbf{B}_k^0,k=0,\ldots ,q-1|\mathbf{B}%
_k^1,k=q,\ldots ,n-1]  \label{eqB.4}
\end{equation}
The vertical line $|$ separates columns with superscripts $0$ and $1$. We
shall call boundary conditions (\ref{eqB.2}) and the corresponding operator $%
L$
\INDEX{Birkhoff}{G.D.Birkhoff}\textbf{Birkhoff-regular}
and write $L\in (R)$ if\hfill
\begin{equation}
\left\{
\begin{array}{llll}
\theta (b^0,b^1,L) & \neq & 0, & n=2q, \\
\theta (b^0,b^1,L) & \neq & 0\;\mbox{and}\;\theta (b^1,b^0,L)\neq 0, &
n=2q+1.
\end{array}
\right.  \label{eqB.5}
\end{equation}
\end{definition}
\noindent
This form of Birkhoff-regularity was invented by
\INDEX{Birkhoff}{G.D.Birkhoff}
\INDEX{Salaff}{S.Salaff} S.Salaff\cite[p.361]{sal68}
who has done a first serious
investigation of the nature of the regularity determinants.

It is worth noting here that recently it was shown \cite{min95a} that
\begin{theorem}[A.M.Minkin]\label{Thsquare}
\INDEX{Min}{A.M.Minkin}
$L\in (R)\Longleftrightarrow L^2\in (R)$.
\end{theorem}
This property serves to reduce odd order problems to the even ones without
any separate treatment of the former.

Let $\{\lambda _j\}_1^\infty $ be the set of all \ev of $L$, $G(x,\xi
,\varrho )$ be its \INDEX{Green}{G.Green}Green function, i.e. the resolvent's $R_\lambda :=(L-\lambda
I)^{-1}$ kernel,
\begin{equation*}
\varrho =\lambda ^{1/n},\qquad |\varrho |=|\lambda |^{1/n},
\end{equation*}
where\hfill

\begin{equation}
\arg \varrho =\arg \lambda /n,\qquad 0\leq \arg \lambda \leq 2\pi ,
\label{eqB.6}
\end{equation}
if $n$ is even and

\begin{equation}
\left\{
\begin{tabular}{lclllcll}
$\arg \varrho$ &$=$&$\arg \lambda /n,\qquad$
&
$-\pi /2$ & $\leq$ &$\arg \lambda$ & $\leq$ & $ \pi /2$; \\
$\arg \varrho$ &$=$&$\pi -\frac{\pi -\arg \lambda }n,\qquad$
&
$\pi /2$  & $\leq $& $|\arg\lambda|$ & $\leq$ & $ \pi$
\notag%
\end{tabular}
\right.  \label{eqB.7}
\end{equation}
if it is odd.
Then

\begin{equation}
\varrho \in S_0=S_1\cup S_2,\quad
S_k=\left\{ (k-1)\frac{\pi}{n}
            \leq \arg \varrho <
	    \frac{k\pi}{n}
 \right\} ,\qquad n\text{ even}
\label{eqB.8}
\end{equation}

\begin{equation}
\varrho \in S_0=S_1\cup S_2=\left\{ |\arg \varrho |\leq \pi /2n\right\} \cup
\left\{ |\pi -\arg \varrho |\leq \pi /2n\right\} ,\qquad n\text{ odd.}
\label{eqB.9}
\end{equation}

Ch.v. $\varrho _j=\lambda _j^{1/n}$ of the operator $L\in (R)$ tend
asymptotically to one arithmetic progression if $n$ is odd and to
two ones if it is even:

\begin{eqnarray}
\varrho _j &=&2\pi j+c+o(1),\quad n=2q+1,  \label{eqB.10} \\
\varrho _j^{\prime } &=&2\pi j+c^{\prime }+o(1),\quad \varrho _j^{\prime \prime
}=2\pi j+c^{\prime \prime }+o(1),\quad n=2q,  \notag
\end{eqnarray}
where constants $c, c^{\prime },c^{\prime \prime }$ are defined
by the leading
coefficients $b_{\hat \imath }^0,b_{\hat \imath }^1$ in boundary conditions
and $\quad j=\pm N,\pm (N+1),\ldots $. The \INDEX{Green}{G.Green}Green function admits the
following remarkable estimate from above off some small $\delta $%
-neighborhoods of the \chv $\varrho _j:$

\begin{equation}
G(x,\xi ,\varrho )=O(\varrho ^{-(n-1)}).  \label{eqB.11}
\end{equation}

Let $S_r(f)$ be the $r$-th partial sum of \eaf expansion including
all summands with $|\varrho _j|\leq r$:

\begin{equation}
S_r(f)=S_r(f,L):=-\frac 1{2\pi i}\int_{\Gamma _r}\int_0^1G(x,\xi ,\varrho )f(\xi )\cdot
n\varrho ^{n-1}\,d\varrho ,  \label{eqB.12}
\end{equation}

\begin{equation}
\Gamma _r=\{\varrho \in S_0,\ \ |\varrho |=r\}.  \label{eqB.13}
\end{equation}

It is suitable to take the integral over contour $\Gamma _r$ in the
principal value sense if it intersects some poles of the
\INDEX{Green}{G.Green}Green
function, i.e. when some \chv $\varrho _j\in \Gamma _r$. In this case the
corresponding summands are taken with the factor $1/2$.

\INDEX{Birkhoff}{G.D.Birkhoff}G.D.Birkhoff has deeply investigated
convergence of such expansions for sufficiently smooth functions
( $f$ is of bounded variation, $f\in V(0,1)$ ).
It happens so that these expansions behave
themselves like the trigonometric ones in all interior points of the main
interval $(0,1)$. For instance, the sum (\ref{eqB.12}) converges to $\frac
12(f(x+0)+f(x-0))$. However, at the end points it converges to some linear
combinations of the limiting values of $f$ at these points whose
coefficients are defined by the boundary forms (\ref{eqB.2}).

Further,
\INDEX{Tam}{J.Tamarkin}J.Tamarkin $(p_k\in C(0,1))$ and
\INDEX{Sto}{M.Stone}M.Stone $(p_k\in L(0,1))$ \cite{tam12,Tam:expan,sto26}
have established the following fundamental result.

\begin{theorem}[J.Tamarkin-M.Stone]
\label{ThA} Let $L\in (R)$ and
\begin{equation}
\sigma _r(f)=\frac 1\pi \int_\mathbb{R}
\frac{\sin r(x-\xi )}{x-\xi }f(\xi )\,d\xi
\label{eqB.14}
\end{equation}
denotes the $r$-th partial sum of the trigonometric Fourier
\INDEX{Fou}{J.Fourier}
integral of
a given function $f\in L(0,1)$. Then
\begin{equation}
\lim\limits_{k\rightarrow \infty }||S_{r_k}(f)-\sigma _{r_k}(f)||_{C(K)}=0
\label{eqB.15}
\end{equation}
for any compact $K$ in $(0,1)$.
\end{theorem}

Here the radii $r_k$ are taken in such a way that
\begin{equation}
dist(\Gamma _{r_k},\left\{ \varrho _j\right\} )\geq \varepsilon >0.
\label{eqB.16}
\end{equation}

$_{}$ For instance, we can take $r_k=2\pi k+\alpha $ with some appropriate $%
\alpha $.

\INDEX{Tam}{J.Tamarkin}J.Tamarkin also established another useful and
important results.

\begin{theorem}[J.Tamarkin]
\label{ThA0} Let $L\in (R),\quad f^{(j)}(x)$ be absolutely continuous for
$j=0,\ldots ,n-1;\;f^{(n)}\in L(0,1)$ and $f$ satisfies boundary conditions
(~\ref{eqB.2}~)~. Then
\begin{equation}
\lim\nolimits_{k\rightarrow \infty }\left\| S_{r_k}(f)-f\right\| _{C[0,1]}=0.
\label{eqB.17}
\end{equation}
\end{theorem}

\begin{theorem}[J.Tamarkin]
\label{ThB}Let $L\in (R)$ and denote $L^0$ an operator defined by the
simplest differential expression $D^n$ and
the leading boundary conditions (\ref{eqB.2})
\begin{equation}
V_j(y)=0,\quad j=0,\ldots n-1.  \label{eqB.19}
\end{equation}
Then
\begin{equation}
\lim\nolimits_{k\rightarrow \infty }\left\|
S_{r_k}(f,L)-S_{r_k}(f,L^0)\right\|
_{C[0,1]}=0.  \label{eqTwh}
\end{equation}
for any function $f\in L(0,1)$.
\end{theorem}

\INDEX{Tam}{J.Tamarkin}J.Tamarkin's investigations were summarized in the
book \cite{Tam:expan}, published on the eve of the 1917 year revolution in
Russia. Therefore its results became known and accessible only after their
reprinting in abbreviated and abbriged form in the fundamental article \cite
{tam28}.

Of course, there were predecessors for this result, namely such theorem has
been earlier established for second order operators in the pioneering works
of
\INDEX{Ste}{V.A.Steklov}V.A.Steklov,
\INDEX{Hob}{E.W.Hobson}E.W.Hobson and
\INDEX{Haa}{A.Haar}A.Haar\cite{ste07,ste10,hob08,haa10}.

\subsection{Main problems}

Theorem \ref{ThA} is remarkable because it completely reduces the question
of Birkhoff's series convergence on any internal compact to an analogous one
\INDEX{Birkhoff}{G.D.Birkhoff}
for the model, namely trigonometric system. The latter is a classic problem
and it is elaborated during the last two hundred years (if not more). In the
meantime theorem \ref{ThA} generated question of the ways of further
investigations. In fact, there are four main directions:

\begin{enumerate}
\item  to consider more general differential expressions, preserving
Birkhoff-re\-gu\-la\-ri\-ty;  \INDEX{Birkhoff}{G.D.Birkhoff}

\item  to consider more general boundary conditions, including, for
instance, those with a Stieltjes integral \INDEX{Sti}{T.Stieltjes}
\begin{equation}
V_j(y)+\int\limits_0^1y^{(j)}(x)d\sigma _j(x)=0,  \label{eqB.18}
\end{equation}
requiring regularity of the
leading boundary conditions (\ref{eqB.19}),
where $d\sigma _j$ denotes a vector-column of height $r_j$
of finite measures
which are continuous at the end points $0,1$;

\item  irregular boundary value problems;

\item  equiconvergence on the whole interval.
\end{enumerate}

\subsection{Polynomial pencils}

Observe that the first two problems were investigated at once by
\INDEX{Tam}{J.Tamarkin}J.Tamarkin himself \cite{tam28}. In particular, he
considered a polynomial pencil
\begin{equation}
y^{(n)}+p_1(x,\varrho )y^{(n-1)}+\ldots +p_n(x,\varrho )y=0,
 \label{eqB.20}
\end{equation}
\begin{equation}
L_i(y)\equiv \sum_{s=0}^n\varrho ^sL_i^{(s)}(y)=0,
\ \ i=1,\ldots,n,
\label{eqB.21}
\end{equation}
where he put
\begin{equation}
L_i^{(s)}(y)\equiv \sum_{l=1}^n\left[
a_{il}^{(s)}y^{(l-1)}(0)+b_{il}^{(s)}y^{(l-1)}(1) +\int_0^1\alpha
_{il}(x)y^{(n-1)}(x)\,dx\right]=0.
\label{eqB.22}
\end{equation}

Herewith it is assumed that the so called characteristic equation
\begin{equation*}
\varphi ^n+p_{10}(x)\varphi ^{n-1}+\ldots +p_{n-1,0}(x)\varphi +p_{n0}(x)=0
\end{equation*}
admits $n$ continuous nonintersecting roots $\varphi _1(x),\ldots ,\varphi _n(x)$
together with a lot of other awkward restrictions upon the coefficients in
the expression (\ref{eqB.20}) and boundary conditions (\ref{eqB.21}).

However, observe that under such generality
(an integral, containing $y^{(n-1)(x)}$ in (\ref{eqB.22}) )
many important instances are lost
or hidden in contrast with more simple situations when they become
transparent. Concretely, the hypothesis that the
boundary forms (\ref{eqB.21})
are polynomials in $\varrho $ yield a very small (from the first glance)
restriction upon the factors $\alpha _{il}$: \emph{their derivative must be
continuous and of bounded variation} \cite[p.30]{tam28}.
Then the main part of the characteristic determinant
(~see (\ref{eqS.2})~)
of
the \bvp  (\ref{eqB.20})-(\ref{eqB.21}) occurs to be a quasipolynomial
where the coefficients by the leading exponential terms are nontrivial
determinants. Demanding them not to vanish
\INDEX{Tam}{J.Tamarkin}J.Tamarkin just transfers the notion of regularity to
this general situation.

Let us, however, assume that these forms don't depend on $\varrho $, i.e.
\begin{equation}
L_i(y)\equiv L_i^{(0)}(y),\qquad L_i^{(s)}(y)=0,\quad s>0.  \label{eqB.23}
\end{equation}

Then it is possible to extract the main part of the characteristic
determinant only for the normalized boundary conditions (\ref{eqB.19}) but
not for the general ones (\ref{eqB.21}). In other words, in the case (\ref
{eqB.23})
\INDEX{Tam}{J.Tamarkin}J.Tamarkin's regularity conditions implicitly demands
that these boundary conditions are equivalent to the
\textit{standard normalized ones}
(\ref{eqB.19}) plus, possibly, some lower order terms.
Hence such a generality in the case of standard \bvp
(not pencils) happens to be apparent.

\section{Stone-regularity}

\subsection{Historical remarks}

Of course, researchers of the beginning of the century understood very
well importance of the problems 1--4 above. However all attempts to
investigate irregular two-point \bvp  yielded only the class of \emph{%
decomposing} boundary conditions when $m$ conditions are taken in one end
point and $n-m$ in another, $m\neq n-m,$ since otherwise such boundary
conditions are regular. In the latter case they are called Sturmian
\INDEX{SturmK}{J.C.Sturm}
conditions and necessarily $n$ is even. It happens so that the
\INDEX{Green}{G.Green}Green
function of the \emph{decomposing}
boundary conditions
has an exponential growth in $\varrho
$ and the associated \ef expansions behaves themselves like Taylor
series or exponential series in the complex plane. Let us note contribution
to the field due to
\INDEX{Khro}{A.P.Khromov}A.P.Khromov,
\INDEX{Eber}{W.Eberhard}W.Eberhard,
\INDEX{Fre}{G.Freiling}G.Freiling,
\INDEX{Ben}{H.Benzinger} H.Benzinger,
\INDEX{SchuB}{B.Schultze}B.Schultze
and
\INDEX{Wol}{M.Wolter}M.Wolter \cite{ebe64,hro66,hro76,fre81,wol83}. Earlier
articles and extensive bibliography may be found in the book \cite{Nai:ldo}
or in the articles just cited.

The problem of finding \emph{good} boundary conditions consists mainly of
difficulties with \INDEX{Green}{G.Green}Green's function estimate from below. Namely, the
advantage of the resolvent's approach used by
\INDEX{Birkhoff}{G.D.Birkhoff}G.D.Birkhoff and his successors, leans heavily upon
the explicit formula
\begin{equation}
G(x,\xi ,\varrho )=\frac{\Delta (x,\xi ,\varrho )}{\Delta (\varrho )}
\label{eqS.1}
\end{equation}
which stems from the method of variation of constants. At the moment we
shall need and recall only a formula for the denominator (which is usually
referred to as \emph{the characteristic determinant}):
\begin{equation}
\Delta (\varrho )=\left| U_j(y_k)\right| _{j,k=0}^{n-1}.
\label{eqS.2}
\end{equation}
Here $\left\{ y_j\right\} _{j=0}^{n-1}$ stands for some \fss of the
equation (\ref{eqB.1}).

When Birkhoff-regularity conditions are violated we are unable to
\INDEX{Birkhoff}{G.D.Birkhoff}
estimate the characteristic determinant from below for \emph{arbitrary
summable coefficients} $p_k$ of the expression (\ref{eqB.1}).

In order to get around this difficulty
\INDEX{Khro}{A.P.Khromov}A.P.Khromov in 1962 and
\INDEX{Ben}{H.Benzinger} H.Benzinger in 1970
introduced a class of S-regular or Stone-regular \bvp  \cite{hro62,ben70}.
Roughly speaking they started to consider not \bvp  but rather operators
because in this approach it is assumed that characteristic determinant
admits an asymptotic expansion. Since the only known situation when it is
possible is when the coeffients in (\ref{eqB.1}) are smooth, $p_k\in
C^\infty (0,1),$ we shall take this hypothesis throughout if otherwise is not
explicitly assumed.
Of course, only a
finite but enough large smoothness is needed but we shall omit here details.
Then there exists a \fss with an exponential asymptotics such that
exponentials are factored by asymptotic power series in $\varrho ^{-1}$. Then
the characteristic determinant happens to be a finite sum of such
exponentials. Now it is possible to indicate its main part
\begin{equation}
\Delta (\varrho )=\Delta _0(\varrho )+\ldots ,\qquad \varrho \in S_k  \label{eqS.3}
\end{equation}
where
\begin{equation}
\Delta _0(\varrho )=\sum c_i(\varrho )\exp (\varrho \sigma _i),\qquad \varrho \in S_k
\label{eqS.4}
\end{equation}
and the exponents $\sigma _i$ have the largest real part.
The sum in (\ref{eqS.4}) consists of two ($n=2q+1,$ $%
i=1,2$ ) or three ($n=2q,i=1,2,3$) summands.
In the even case $\Re \sigma
_i>\Re \sigma _{i+1}$.
In the odd one $\sigma _1$ contains all $\varepsilon
_j$ such that $\Re (i\varrho \varepsilon _j)\geq 0$ throughout the sector $S_k$
under consideration while $\sigma _2$ differs from it by a summand $%
\varepsilon _q$. The latter is the unique value such that $\Re (i\varrho
\varepsilon _j)$ changes sign in the corresponding
$\varrho $-sector
$S_k \ (0\le
j\le n-1$).

Call a function $c(\varrho )$ an asymptotic function of order $\alpha $ ($%
\alpha $ real) in the sector $S_k$ if
\begin{equation}
\exists d=\lim c(\varrho )/\varrho ^\alpha ,\qquad \varrho \rightarrow \infty ,\qquad
\varrho \in S_k  \label{eqS.5}
\end{equation}
and
\begin{equation}
d\neq 0.  \label{eqS.6}
\end{equation}

\begin{definition}[A.P.Khromov, H.Benzinger]
\INDEX{Khro}{A.P.Khromov}
\INDEX{Ben}{H.Benzinger}
\label{Sreg} A \bvp  is called Stone-regular\INDEX{Sto}{M.Stone}\linebreak
(shortly S-regular)
in the
sector $S_k$ if the coefficients $c_i(\varrho )$ in (\ref{eqS.4}) are
asymptotic power functions of orders $\alpha _i$,
 respectively (hence, the limits $d_i\neq 0$).
Then the corresponding operator $L$ is called of type ($%
\alpha _1,\alpha _3$) if $n$ is even or of type ($\alpha _1,\alpha _2$) if
$n$ is
odd.
\end{definition}

The Birkhoff's regularity corresponds to the case when\INDEX{Birkhoff}{G.D.Birkhoff}
\begin{equation*}
\alpha _1=\alpha _2=\chi ,\quad n\quad \text{odd;\qquad }\alpha _1=\alpha
_3=\chi ,\quad n\quad \text{even; \qquad },
\end{equation*}
where the quantity
\begin{equation}
 \label{eq:total}
\chi :=\sum_{j=0}^{n-1}jr_j
\end{equation}
is called
a \textit{total} order of the \bvp \cite[p.194]{shk83}.
At the first glance this definition depends on the choice of the sector
$S_k,k=1,2$.
However,
\INDEX{Eber}{W.Eberhard}W.Eberhard and
\INDEX{Fre}{G.Freiling}G.Freiling proved independence of the orders $\alpha
_i$ of the sector's choice \cite{ebfr78}. Later
\INDEX{SchuB}{B.Schultze}B.Schultze made a complement to the theory adding
to the definition of S-regularity the case when $d_2=0$ but $\alpha _2\leq
(\alpha _1+\alpha _3)/2$ \cite{bschu79a,bschu79b}.

\INDEX{Green}{G.Green}Green function of S-regular problem obeys a polynomial estimate

\begin{equation}
G(x,\xi ,\varrho )=O(|\lambda |^a).  \label{eqS.7}
\end{equation}

Generally speaking, it is \emph{worse} than (\ref{eqB.11}). Presently we are
able to prove rigorously that $a>-(n-1)/n$ for irregular \bvp
 but this fact falls out of the review's goals
and we plan to expose it elsewhere. In view of (\ref{eqS.7}) the next
theorem looks quite natural.

\begin{theorem}
\label{ThC}\cite[Theorem 5]{hro62} Let $L$ be a S-regular operator of the
type $(\alpha _1,\alpha _2)$. If $f\in D_{L^m}$ with
$m=\frac 1nEntier\left[
\sum\nolimits_{j=0}^{n-1}r_j-\min (\alpha _1,\alpha _2)-(n-1)\right] +2$,
then
\begin{equation}
\lim\nolimits_{k\longrightarrow \infty }\left\| S_{r_k}(f)-f(x)\right\|
_{C(0,1)}=0  \label{eqS.9}
\end{equation}
The exponent $m$ here is exact.
\end{theorem}

Advanced results in that theory may be found in
in A.A.Shkalikov's \INDEX{Shka}{A.A.Shkalikov}
articles \cite{shk83,shktr94}, see also \cite{ben70}.
However, no
statements like theorems \ref{ThA},\ref{ThB} were proved. In subsection
below we give some of
\INDEX{Khro}{A.P.Khromov}A.P.Khromov's results omitting theorems concerning
Riesz summability of such expansions.

\subsection{Finite functions}

At first let us consider the case of equal orders.

\begin{theorem}
\label{ThHr5}\cite[Theorem 5]{hro62}
Let $L$ and $L^{\prime }$ be S-regular
differential operators of types $(\alpha ,\alpha )$ and $(\alpha ^{\prime
},\alpha ^{\prime })$, respectively. Given a a number $\delta ,\ \ 0<\delta
\leq 1/2$ and a summable function $f$ which vanishes off the interval $%
K=[\delta ,1-\delta ]$ the following relation holds

\begin{equation}
\lim\nolimits_{k\longrightarrow \infty }\left\|
S_{r_k}(f,L)-S_{r_k}(f,L^{\prime })\right\| _{C(K)}=0.  \label{eqS.10}
\end{equation}
\end{theorem}

In the second theorem it is assumed that these operators have the same type $%
(\alpha _1,\alpha _2)$ but these numbers possibly differ.

\begin{theorem}
\label{ThHr6}\cite[Theorem 6]{hro62}
Assume that $c_i(\varrho )=c_i^{\prime
}(\varrho )[1],\ \ i=1,2$. Then for any summable function $f$ vanishing off
some interval $[\delta _1,1-\delta _2]\subset (0,1)$ the relation (\ref
{eqS.10}) remains valid. The interval $K$ of equiconvergence is as follows.

If $|\alpha _2-\alpha _1|\leq 1$ then $K=[\delta _3,1-\delta _3].$

If $|\alpha _2-\alpha _1|>1$ then


\begin{tabular}{||>{$}l<{$}|>{$}l<{$}|>{$}l<{$}|>{$}l<{$}||}
\hline
n & K & \varepsilon  & \alpha _2-\alpha _1 \\
\hline
\hline
\mu +1 & \left[ \varepsilon ,1-\delta _3\right]  &
1 - \vert \alpha _2-\alpha_1 |^{-1}
-\delta _2 & \alpha _2-\alpha _1>+1 \\
\hline
4\mu +3 & \left[ \delta _3,\varepsilon \right]  & \vert\alpha _2-\alpha
_1|^{-1}+\delta _1 & \alpha _2-\alpha _1<-1\\
\hline
\end{tabular}

and it is subject to the evident restrictions that
$K\subset (0,1)$ and its left
end is less than the right one.
\end{theorem}

There are also examples demonstrating sharpness of the theorem's conditions.
The third theorem deals with a \bvp  where some of the coefficients $p_k$
vanish:

\begin{theorem}
\label{ThHr7}\cite[Theorem 7]{hro62}
Assume that $L$ is defined by expression
\[
l(y)=y^{(n)}+p_{n-k}y^{(k)}+\dots +p_n(x)y
\]
 and similarly for $L^{\prime }$
with $l^{\prime }(y)$ of the same type.
Suppose also that
\begin{equation}
c_i=c_i^{\prime }\left[ 1+O(\varrho ^{k+1-n})\right] ,{}\qquad 0\leq k\leq
n-2-|\alpha _2-\alpha _1|  \label{eqS.11}
\end{equation}
Implicitly (\ref{eqS.11}) implies that both differential operators in
question are of the same type $(\alpha _1,\alpha _2)$. Then again (\ref
{eqS.10}) holds for the function $f$ and the interval $K$ as in the theorem
\ref{ThHr6}.
\end{theorem}

\begin{remark}
The proof of these results is difficult but it is easy to guess the interval
of equiconvergence applying our theorem \ref{ThMgenser1} (see below).
\end{remark}

Of course, it is possible to consider \bvp
with a polynomial growth of the resolvent from an abstract point of view,
i.e. a priori assuming the inequality (\ref{eqS.7}) to be fulfilled
without regard of how this could be obtained.
However, we dare to conjecture that
\begin{conjecture}
\label{conSreg}Given a \bvp  subject to estimate (\ref{eqS.7}) with smooth
coefficients $p_k$, then it is Stone-regular. Moreover, we believe that this
assertion is also valid for general boundary conditions (\ref{eqB.18}).
\end{conjecture}

\subsection{Strongly irregular boundary conditions}

Since \bvp  with decomposing boundary conditions have an exponential growth
of the resolvent it is hardly possible to expect any kind of equiconvergence
with a trigonometric series expansion. However,
\INDEX{SchuB}{B.Schultze}B.Schultze succeeded in finding such a phenomenon
for a class of \bvp  which are \emph{partially decomposing}. Further we
describe this unexpectable result (see \cite{bschu79b}).

Let
$0<m<n,\quad \alpha $
be an
$m\times n,\ \beta $
and
$\gamma $ be $%
(n-m)\times n$ matrices, respectively.

\begin{definition}
\cite{bschu79b}
A two-point boundary conditions are called
\emph{strongly irregular}
if they are equivalent to boundary conditions of the form:
\begin{equation}
My^{\vee }(0)+Ny^{\vee }(1)=0,\qquad y^{\vee }(x):=\left( D^jy(x)\right)
_{j=0}^{n-1}  \label{eqBS.1}
\end{equation}
with
\begin{equation}
M=\left( \frac \alpha \gamma \right) ,\qquad N=\left( \frac 0\beta \right)
,\qquad \text{ when \quad }m>n-m  \label{eqBS.2}
\end{equation}
or with
\begin{equation*}
M=\left( \frac \alpha 0\right) ,\qquad N=\left( \frac \gamma \beta \right)
,\qquad \text{ when \quad }n-m>m.
\end{equation*}
Here $0$ stands for a zero-matrix of appropriate size which is clear from
the context.
\end{definition}

To be definite we shall confine ourselves in what follows to the first case
and consider differential operators $L$ and $L^{\prime }$ defined by the
equation (\ref{eqB.1}) and boundary conditions (\ref{eqBS.2}) (operator $L$)
and
\begin{equation}
M^{\prime }y^{\vee }(0)+Ny^{\vee }(1)=0,\qquad M^{\prime }=\left( \frac
\alpha 0\right)  \label{eqBS.3}
\end{equation}
(operator $L^{\prime }$ ).

\begin{definition}
If $A=\left[ a_1,\ldots a_m\right] $ is a
$k\times m$-matrix of rank$\ k,\quad m\geq k$
with columns
$a_i(i=1,\ldots m)$,
we define the
\textbf{weight}
of the matrix\textbf{\ }$A$ as in \cite[definition 3]{bschu79a}:
\begin{equation*}
weight(A):=\max \left\{ i_1+\cdots +i_k\quad |\quad \det
\left[
a_{i_1},\ldots ,a_{i_k}\right] \neq 0
\right\} .
\end{equation*}
\end{definition}

\begin{theorem}[B.Schultze]
\INDEX{SchuB}{B.Schultze}
\label{ThBSchu}\cite[Remark to Theorem 6]{bschu79b}Assume that the weight of
the matrix $\beta $ does not increase, if an arbitrary column of $\beta $ is
replaced by an arbitrary column of the matrix $\gamma $. Then the
equiconvergence relation (\ref{eqS.10}) is valid for any $f\in L(0,1)$ and
any compact $K\subset (0,1)$.
\end{theorem}

Evidently this theorem has theorem \ref{ThB} as a counterpart and
indicates that despite an exponential growth the
\INDEX{Green}{G.Green}Green functions of both
operators $L$ and $L^{\prime }$ are very close to one another.

\section{General differential expressions}

Investigation of the convergence of \ef expansions at the end
points or in their neirborhoods met serious difficulties and has been
examined only for sufficiently smooth functions, say for functions of
bounded variation, $f\in V[0,1]$ \cite{tam28,bir08b}.Among the last papers
let us note \cite{kalu94}. Therefore writers on the topic concentrated their
efforts upon various generalizations of the differential expression $l(y)$
replacing it by
\begin{equation}
l_1(y)=D^ny+Fy  \label{eqG.1}
\end{equation}
where $F$ denotes a linear operator dominated in a certain sense by the
first summand. Here we want to distinguish the following results.

\subsection{Nonsmooth coefficient by the $(n-1)$th derivative}

First, let us take
\begin{equation}
l_2(y)=D^ny+Fy,\quad Fy=\sum_{k=0}^{n-1}p_k(x)D^ky.  \label{eqNS.1}
\end{equation}
If
\begin{equation}
p_{n-1}(x)\in C^{n-1}[0,1],  \label{eqNS.2}
\end{equation}
then upon substitution
\begin{equation}
y=Vz,\quad V(x)=\exp \left( -\frac in\int_0^xp_1(s)ds\right)  \label{eqNS.3}
\end{equation}
we pass
to the standard form (\ref{eqB.1}%
) of the differential expression. However, if (\ref{eqNS.2}) is broken,
namely
\begin{equation}
p_{n-1}(x)\in L[0,1],  \label{eqNS.4}
\end{equation}
the situation changes abruptly.
\INDEX{Rykh}{V.S.Rykhlov}V.S.Rykhlov established in this case an existence
of a \fss with an exponential asymptotics as long as in 1977 \cite{ryh77}.
This was a breakthrough in the theory and soon after that he built a
nontrivial analogue of the \INDEX{Birkhoff}{G.D.Birkhoff}Birhoff's theory.
Of course, here also appeared
the
\INDEX{Birkhoff}{G.D.Birkhoff}Birhkhoff-regularity conditions
in a slightly modified form: one has
only to replace $b_j^1$ by $b_j^1\cdot V(1)$ in
(\ref{eqB.3})-(\ref{eqB.4}). The following
main result belongs to this author \cite{ryh83,ryh84}.

\begin{theorem}[V.S.Rykhlov]
\INDEX{Rykh}{V.S.Rykhlov}
\label{ThRy}Given a differential operator $L\in (R)$ defined by (\ref{eqNS.1}%
) with summable coefficents
$p_j(x),\ \ j=0,\dots,n-1$.
Then the equiconvergence (\ref{eqB.15})
remains valid provided one of the following relations is fulfilled:

\begin{enumerate}
\item  $p_{n-1}(x)\in L^q[0,1],\quad f(x)\in L^p[0,1],\qquad \frac 1p+\frac
1q<1;$

\item  $p_{n-1}(x)\in H_\infty ^\alpha [0,1],\quad f(x)\in H_1^\beta
[0,1],\qquad \alpha +\beta >1;$

\item  $p_{n-1}(x)\in H_1^\alpha [0,1],\quad f(x)\in H_\infty ^1[0,1],\qquad
\alpha +\beta >1.$
\end{enumerate}

Moreover, the following estimate with a modified
\INDEX{Dir}{L.Dirichlet}Dirichlet kernel holds
\begin{equation}
\left\| S_{r_k}(f,L)-\left( V\sigma _{r_k}V^{-1}\right) (f)\right\|
_{C(K)}\leq \left( \frac{\log r}{\log ^{\alpha +\beta }r}+\frac 1{\log ^\alpha
r}+\frac 1{\log ^\beta r}\right)   \label{eqNS.5}
\end{equation}
for any compact $K\subset (0,1)$.
\end{theorem}

Here

$H_1^\alpha [0,1]=\left\{ g(x)\in L[0,1]\quad |\quad \varpi _1(g,\delta)
=
O\left( \log ^{-\alpha }\frac 1\delta \right) ,\ \alpha >0\right\} ,$\\[5pt]

$H_1^0[0,1]=L[0,1];$\\[5pt]

$H_\infty ^\alpha [0,1]=\tilde H_\infty ^\alpha [0,1]
\qquad (0\leq \alpha \leq 1);$\\[5pt]

$H_\infty ^\alpha [0,1]=\tilde H_\infty ^\alpha [0,1]\cup V[0,1]
\qquad (\alpha >1);$\\[5pt]

$\tilde H_\infty ^\alpha [0,1]=\left\{ g(x)\in C[0,1]\quad |\ \varpi
_1(g,\delta )=O\left( \log ^{-\alpha }\frac 1\delta \right) ,\ \alpha
>0\right\} $\\[5pt]

$\tilde H_\infty ^0[0,1]=L_\infty [0,1].$\\[5pt]

Concrete examples clearly demonstrate sharpness of the inequalities $\frac
1p+\frac 1q<1,\quad \alpha +\beta >1$ in the theorem's statement. Obviously,
(\ref{eqNS.5}) signifies \emph{an} \emph{equiconvergence with a rate}.
Presently there appeared several new results in this direction, see, for
instance \cite{ryh96a}, but we shouldn't go into further details.

\subsection{Integral operators with a Green-type kernel}
\INDEX{Green}{G.Green}

Next, consider the equiconvergence problem for integral operators of the form
\begin{equation}
Af(x)=\int_0^1A(x,t)f(t)dt.  \label{eqI.1}
\end{equation}
Obviously, any minimal function system in $L^2(0,1)$ occurs to be an
\ef family of some integral operator. Hence, here we deal with a
most general situation.

Set
\begin{equation}
A_{s,j}(x,\xi )=\overline{D_\xi ^j\left( \overline{D_x^sA(x,\xi )}\right) }%
,\qquad x\neq \xi  \label{eqI.2}
\end{equation}
\begin{equation}
\Delta A_{s,j}(x):=A_{s,j}(x,\xi )|_{\xi =x+0}^{\xi =x-0}  \label{eqI.3}
\end{equation}
and suppose the following conditions to be fulfilled:

\begin{enumerate}
\item[i)]  the derivatives $A_{s,j}(x,\xi )$ are continuous whenever $t\leq x
$ or $x\leq t$ ($s,j=0,\ldots ,n)$;

\item[ii)]  the jumps $\Delta A_{s,j}(x)\in C^{n-1-j}[0,1],$ ($s,j=0,\ldots
,n-1)$;

\item[iii)]  $A_{s,0}(x,\xi ),\ s=0,\dots,n-1$ are continuous and the
             $(n-1)$th derivative $A_{n-1,0}(x,\xi )$ is discontinuous
             at the line $x=\xi $%
\begin{equation}
\Delta A_{s,0}(x)=i\cdot \delta _{s,n-1};\quad s=0,\ldots ,n-1;
\label{eqI.4}
\end{equation}

\item[iv)]  operator $A$ admits no zero \ev
\end{enumerate}

Then the inverse $A^{-1}$ occurs to be an integro-differential operator of
the form
\begin{equation}
l_3(y)=(E+N)\left( D^ny+\alpha y\right)  \label{eqI.5}
\end{equation}
with some boundary conditions
\begin{equation}
U_j(y)=\int_0^1y(x)\varphi _j(x)dx,\quad j=0,\ldots ,n-1  \label{eqI.6}
\end{equation}
where $E$ denotes an identity operator in $L^2(0,1)$ and $N$ stands for an
integral operator like (\ref{eqI.1}) with a kernel $N(x,t)$. The latter is
separately continuous in both triangles $t\leq x$ or $x\leq t$; $U_j(y)$ are
$n$ linear independent forms of variables
\[
D^jy(0),D^jy(1),\quad j=0,\ldots,n-1;
\]
$\alpha $
is some complex number. Note that the kernel
$N(x,t)$ and
the functions $\varphi _j(x)$ may be calculated efficiently through the initial
kernel $A(x,t)$.

Hence, any kernel $A(x,t)$ subject to the aforementioned restrictions
i)--iv) constitutes a
\INDEX{Green}{G.Green}Green function of the integro-differential operator (%
\ref{eqI.5})--(\ref{eqI.6}).

Evidently, the boundary conditions (\ref{eqI.6})
may be called \emph{natural}.

\begin{theorem}[A.P.Khromov]
\label{ThHrI2}\cite{hro81}%
\INDEX{Khro}{A.P.Khromov}Assume that

\begin{enumerate}
\item[1)]  the integral operator $A$ satisfies restictions
\textrm{i)-iv)} and in
addition (\ref{eqI.4}) is also valid for $s=n$;

\item[2)]  the leading forms $U_j(y)$ in \emph{natural} boundary conditions
(\ref{eqI.6}) are Birkhoff-regular; \INDEX{Birkhoff}{G.D.Birkhoff}

\item[3)]  uniformly in $\xi \in [0,1]$%
\begin{equation}
\underset{0}{\overset{1}{Var}}_xA_{n0}(x,\xi )\text{ is bounded.}
\label{eqI.7}
\end{equation}
\end{enumerate}

Then (\ref{eqB.15}) is true where now $S_{r_k}(f)$ stands for the partial
sum of the \ef expansion associated with the operator $A$. The sum
includes all summands whose \ev are greater than $r_k^{-n}$ in
absolute value.
Recall that if $\mu $ is an \ev of the integro-differential operator
(\ref{eqI.5})--(\ref{eqI.6})
then $\lambda =\mu ^{-1}$ is an \ev of $A$.

If in addition the kernel $A(x,\xi )$ is symmetric,
$A(x,\xi )=\overline{A(\xi ,x)},$
then condition 2) may be omitted \cite{min77b}.
It merely follows from the
regularity of \sa boundary conditions \cite[$n$ even]{sal68,fie72},
and \cite[$n$ odd]{min77a}.
\end{theorem}

Conditions of this theorem are exact, no one of conditions 1)--3) may be
removed. Add some words concerning apriori restrictions iii)--iv). Condition
iv) is needed because otherwise equiconvergence for a function $f_0$ from
the $A$ kernel forces it to be rather specific, namely to have a uniformly
convergent trigonometric series. Now about iii). Given an integral operator $%
A$ with sufficiently smooth kernel, it is possible to indicate an integral
operator $B$ retaining the same \eaf system and such that iii) is
satisfied with some even $n$.

\emph{Hence, condition iii) singles out a \textbf{canonical} operator among
all integral operators with the same \eaf system admitting equiconvergence
with a trigonometric series expansions. }

It is worth noting here that perhaps the first time such an
integral operator appeared in
\INDEX{Lan}{R.Langer}R.Langer's article \cite{lan26} where the case $n=1$
was treated.

A partial case of integral operators constitute finite convolution
operators:
\begin{equation}
Af(x)=\int_0^1A(x-t)f(t)dt,\ \ 0\leq x \leq 1.
\label{eqI.8}
\end{equation}
From the theorem \ref{ThHrI2} it follows directly

\begin{theorem}[A.P.Khromov]
\label{ThHrI4}\cite{hro81}%
\INDEX{Khro}{A.P.Khromov} Assume that

\begin{enumerate}
\item  the function $A(x)\in C^{2n}$ for $x\geq 0$ and $x\leq 0$;

\item  $D^jA(+0)-D^jA(-0)=i\delta _{j,n-1},\qquad j=0,\ldots ,n;$

\item  there exists an inverse $A^{-1}.$
\end{enumerate}

Then relation (\ref{eqB.15}) holds true.
\end{theorem}

Note that
\INDEX{Pal}{B.V.Pal'tsev}B.V.Pal'tsev has earlier investigated
a convolution operator
whose kernel coincides with a restriction of the Fourier transform
\INDEX{Fou}{J.Fourier}
of a rational function \cite{pal72}:
\begin{equation}
A(x)=\int_{\mathbb{R}}\frac{P(t)}{Q(t)}e^{-ixt}dt,  \label{eqI.9}
\end{equation}
where $P(t)=t^p+at^{p-1}+\ldots $, $Q(t)=t^q+bt^{q-1}+\ldots $ are
polynomials in $t$ of orders $p$ and $q$, respectively.

\begin{theorem}[B.V.Pal'tsev]
\label{ThPal}%
\INDEX{Pal}{B.V.Pal'tsev}Let $q^{+} (q^{-})$ denote the number of roots of
$Q(z)$
in the upper (lower) half-plane. Assume that
\[
n:=q-p\geq \min \left\{q^{+},q^{-}\right\}
\]
and the poly\-no\-mi\-al
\[
P(t)\overline{Q(t)},\quad t\in \mathbb{R}
\]
has only real coefficents.
Define the function $A(x)$ through equality
(\ref{eqI.9}).
Then (\ref{eqB.15}) is valid for an integral convolution operator
$A_1(f)$
with the kernel $A_1(x-t)$ where
\begin{equation*}
A_1(x)=\frac{(-1)^n}{2\pi }A(x)\exp \left( \theta x\right) ,\qquad \theta
=\frac in(a-b).
\end{equation*}
\end{theorem}

\subsection{Functional--differential perturbation}

Let us pass now to a general approach when perturbation $F$ is an abstract
linear operator. Let $F$ be a bounded linear operator from the
\INDEX{H\"older}{O.H\"older}H\"older space $%
C^\gamma [0,1]$ into the space $L(0,1)$ and impose a restriction
\begin{equation}
0\leq \gamma < n-1.  \label{eqR.1}
\end{equation}

Consider a functional--differential operator $L$ defined by
(\ref{eqG.1})
and boundary conditions (\ref{eqB.18}).
For a summable function $f$ let
$\tilde f(x)$ be its extension to the whole axis which
vanishes off $[0,1]$.
Denote
$
\varpi ^q\left( \delta ,f\right)
$
the $q$th modulus of continuity of $f$ in
$L(0,1)$:
\begin{equation*}
\varpi ^q\left( \delta ,f\right) :=\sup_{0\leq h\leq \delta }\int_{\mathbb{R}%
}\left| \Delta _h^q\left( \tilde f,x\right) \right| dx,
\end{equation*}
\begin{equation*}
\Delta _h^q\left( \tilde f,x\right) :=\sum_{s=0}^q(-1)^qC_q^s\tilde f(x+hs).
\end{equation*}

\INDEX{Radz}{G.V.Radzievskii}G.V.Radzievskii and
\INDEX{Gom}{A.M.Gomilko}A.M.Gomilko established
\begin{theorem}[A.M.Gomilko, G.V.Radzievskii]
\label{ThRG1}\cite{gora91}Let the boundary forms $V_j(y)$ be
Birkhoff-regular. Then\INDEX{Birkhoff}{G.D.Birkhoff}
\begin{equation}
\left\| S_{r_k}(f)-\sigma _{r_k}^\pi (f)\right\| _{C\left[ \delta ,1-\delta
\right] }\leq d_q\cdot \frac{r_k}{1+r_k\delta }\varpi ^q\left( \frac
1{r_k},f\right).   \label{eqR.2}
\end{equation}

Here $q\equiv 1$ if there are measures in (\ref{eqB.18}) and $q=1,2,\ldots $
otherwise.
\begin{equation}
\sigma _r^\pi (f):=\sum_{\left| j\right| \leq r/2\pi
}(f,e_j)_{L^2(0,1)}e_j,\quad e_j:=\exp (2\pi ijx),\ \ 0\leq x\leq 1
\label{eqR.3}
\end{equation}
--- stands for a partial sum of a trigonometric Fourier series.
\INDEX{Fou}{J.Fourier}
\end{theorem}

Further,
\INDEX{Radz}{G.V.Radzievskii}G.V.Radzievskii investigated the case of
abstract perturbations $F$ acting from $C^\gamma [0,1]$ into $L^p(0,1),\quad
1\leq p<\infty $ or into the Sobolev
\INDEX{Sob}{S.L.Sobolev}
space $W_1^1$ \cite{rad94a}. This
generalization overlaps such important situation as
\begin{equation}
Fy=Ny^{(n-1)}+\sum_{k=0}^{n-2}p_k(x)D^ky,  \label{eqR.4}
\end{equation}
$N$ being an integral operator with a smooth kernel,
\begin{equation*}
D^{n-1}N(x,t)\in L_1\left( \left[ 0,1\right] \times \left[ 0,1\right]
\right) .
\end{equation*}
Herewith he also gave estimates of the rate of convergence of \eaf
expansions in terms of modulus of continuity of the function in question or
in terms of certain special $K$-functionals.

However, it is necessary to note that
\INDEX{Rykh}{V.S.Rykhlov}V.S.Rykhlov's operator hasn't been covered yet by
this approach. It corresponds to the \emph{critical} case
\begin{equation}
\gamma =n-1,  \label{eqR.5}
\end{equation}
when the summand $Fy$ can not be viewed of as a \textit{weak}
perturbation of the
main term $D^ny$. Hence, it is highly desirable to extend
their abstract
approach to the case (\ref{eqR.5}).

\subsection{Unconditional equiconvergence}

Recently perturbations similar to (\ref{eqR.4}) were also investigated by
\INDEX{Bas}{A.G.Baskakov}A.G.Baskakov and
\INDEX{Katza}{T.K.Katzaran}T.K.Katzaran in \cite{baka88} from another point
of view. Namely, they set
\begin{equation*}
Fy=Ky^{(n-1)}+K_0\left( y^{(n-1)}+\alpha y\right) ,\qquad \alpha \neq 0,\
\alpha \in \mathbb{C}
\end{equation*}
where $K$ denotes a finite-dimensional operator acting from
$C(0,1)$ into
$L^2(0,1)$
and $K_0$
is a Hilbert-Schmidt
\INDEX{Hil}{D.Hilbert}
\INDEX{Schmidt}{E.Schmidt}
operator in $L^2(0,1)$.
They
consider a functional-differential operator $L$ defined by (\ref{eqG.1}),(%
\ref{eqI.6}) with $\varphi _j\in L^2(0,1)$
and a priori assume that these
boundary conditions are strongly regular \cite[p.71]{Nai:ldo}.

Next, this operator is compared with a model one $\tilde L$,
generated by the
differential expression $y^{(n)}$ with the same boundary conditions.
They proved that $L$ is spectral in
\INDEX{Dun}{N.Dunford}N.Dunford\INDEX{SchJ}{J.Schwartz}-J.Schwartz'
sense and all its \ev are simple except, perhaps, for a
finite number.

Let $P_\sigma (\tilde P_\sigma )$ denote the spectral projector on the
portion of the spectrum of operator $L(\tilde L)$ falling into the set $%
\sigma .$ Let us state one of their main results.

\begin{theorem}[A.G.Baskakov, T.K.Katzaran]
\label{ThBK}There exists a simultaneous enumeration of the spectra of both
operators $L$ and $\tilde L$ such that
\begin{equation*}
\left\| P_\sigma -\tilde P_{\tilde \sigma }\right\| _{L^2(0,1)\rightarrow
L^2(0,1)}\longrightarrow 0
\end{equation*}
whenever
$\min \left\{ \left| \lambda _j\right| \quad |\ \lambda _j\in
\sigma \right\} \longrightarrow \infty $.
\end{theorem}

Note that such phenomenon hasn't been observed earlier even in the case of
ordinary differential operators. Let us clarify that spectrality of
both operators yields only an estimate from above
\begin{equation*}
\left\| \vphantom{\tilde P_\sigma}
               P_\sigma \right\| _{L^2(0,1)\rightarrow L^2(0,1)},
\left\| \tilde P_\sigma \right\| _{L^2(0,1)\rightarrow L^2(0,1)}
\leq C
\end{equation*}
uniformly with respect to all choices of the subset
$\sigma \in \mathbb{R}$,
whence
\begin{equation*}
\left\| P_\sigma -\tilde P_\sigma \right\| _{L^2(0,1)\rightarrow
L^2(0,1)}\leq 2C.
\end{equation*}
Of course, this is weaker than the theorem's \ref{ThBK} assertion.

\subsection{General boundary conditions}

B.v.p. in a finite interval seem to be useful models of general
nonself-adjoint operators. It is especially true for problems
with general functionals in boundary conditions. Below we give an account of
some known results in this direction related to the equiconvergence problem.

\subsection{First order \protect\bvp }

\INDEX{Sed}{A.M.Sedletskii}A.M.Sedletskii \cite{sed75,sed82} investigated a
differentiation operator
\begin{equation*}
l(y)=y^{\prime },\qquad -1\leq x\leq 1
\end{equation*}
with a difficult to study \emph{smeared }condition:
\begin{equation}
U(y)=\int_{-1}^1\frac{k(t)}{\left( 1-|t|\right) ^\alpha }y(t)dt=0\qquad
(0<\alpha <1).  \label{eqC.1}
\end{equation}
He obtained an equiconvergence on the whole interval with a \emph{shearing}
weight.

\begin{theorem}[A.M.Sedletskii]
\INDEX{Sed}{A.M.Sedletskii}
If $\underset{-1}{\overset{1}{Var}}k<\infty ,\quad k(1-0)\cdot k(-1+0)\neq 0,
$ then for any $f\in L[-1,1]$%
\begin{equation}
\lim_{r\longrightarrow \infty }
\left\|
\left( 1-|x|\right) \cdot
 \left[
    S_r(f)-\sigma _r(f)
 \right]
\right\| _{C[-1,1]}=0.  \label{eqC.2}
\end{equation}
Of course, here $r\longrightarrow \infty $ remaining at least at a fixed
positive distance of \ev $\lambda_k$.
Recall that the latter lies in a strip
$
 |\Im\lambda_k|\leq const.
$
\end{theorem}
For systems of exponentials with a half-bounded spectrum,
$
 \inf\Im\lambda_k>-\infty,
$
such kind of results as (\ref{eqC.2})
was also obtained in \cite[Theorem 4.1]{HNP81}.

Further
\INDEX{Kab}{S.N.Kabanov}S.N.Kabanov \cite{kab90a}
studied a differentiation operator with
a boundary condition of general form:
\begin{equation}
\alpha y(-1)+\int_{-1}^1y^{\prime }(t)h(t)dt=0  \label{eqC.3}
\end{equation}

\begin{theorem}[S.N.Kabanov]
\INDEX{Kab}{S.N.Kabanov}
\label{ThKab1}Assume that $h(t)\in L^q[-1,1]$,\\
$M_1h(t),\  M_1\tilde h(t)\in V[-1,1]$ where
\begin{equation*}
M_1h(t)=\int_{-1}^t\frac \partial {\partial \xi }\frac{(\tau -t)^{\xi
+\alpha -1}}{\Gamma (\xi +\alpha )}\mid _{\xi =0}h(\tau )d\tau ,
\end{equation*}
$\tilde h(t):=h(-t)$ and $M_1h(1)\cdot M_1\tilde h(1)\neq 0$.
Then for any
$f\in L^p[-1,1],\quad \frac 1p+\frac 1q=1,$
and for any
$\delta \in (0,1/2)$
relation (\ref{eqB.15}) is satisfied.
\end{theorem}

\INDEX{Amv}{O.I.Amvrosova}O.I.Amvrosova \cite{amv83} considered a higher
order operator
\begin{equation}
l(y)=y^{(n)}  \label{eqC.4}
\end{equation}
with \emph{smeared} boundary conditions containing power singularities
\begin{equation}
U_j(y)=\int_{-1}^1\varphi _j(t)y(t)dt+\int_{-1}^1\frac{k_j(t)}{\left(
1-|t|\right) ^{\alpha _j}}y^{(p_j)}(t)dt=0,\quad j=1,\ldots ,n;
\label{eqC.5}
\end{equation}

\begin{eqnarray*}
n-1 \geq p_1\geq \ldots \geq p_n\geq 0,\quad 0<\alpha _j<1, \\
\underset{-1}{\overset{1}{Var}}k_j <\infty ,\quad \underset{-1}{\overset{1%
}{Var}}\varphi _j<\infty ,
\end{eqnarray*}
singled out a class of regular boundary conditions and obtained for it the
following result.

\begin{theorem}[O.I.Amvrosova]
\INDEX{Amv}{O.I.Amvrosova}
Let the boundary conditions (\ref{eqC.5}) be regular. Then (\ref{eqC.2}) is
valid for any function $f\in L[-1,1]$.
If besides
$
f\in L^p[-1,1],\ p>1,
\quad \alpha _\nu >\frac 1q,\quad\frac 1p+\frac 1q=1
$
then
\begin{equation}
\lim_{r\longrightarrow \infty }\left\| \left( 1-|x|\right) ^\gamma \left[
S_r(f)-\sigma _r(f)\right] \right\| _{C[-1,1]}=0  \label{eqC.6}
\end{equation}
for any $\gamma >\frac 1p$.
\end{theorem}

Later
\INDEX{Kab}{S.N.Kabanov}S.N.Kabanov carried over his first order theorem
\ref{ThKab1}
to the operator (\ref{eqC.4}) \cite{kab90a}
 with general boundary conditions:
\begin{equation}
U_j(y)=\sum_{k=0}^{n-1}a_{j,k}y^{(k)}(-1)+
\int_{-1}^1y^{(n)}(t)h_j(t)dt=0,\quad
j=1,\ldots ,n.  \label{eqC.7}
\end{equation}
 He studied a rather difficult case
when the numbers below don't vanish,
\begin{equation}
\beta _j=D_1^{-\alpha_j } h_j(t)_{|t=1}\neq 0,\quad
\gamma _j=D_1^{-\alpha_j} \tilde h_j(t)_{|t=1}\neq 0,
\quad j=1,\dots,n,
\label{eqC.8}
\end{equation}
where
\begin{equation*}
D_1^{-\alpha }h=\frac 1{\Gamma (\alpha )}\int_t^1(\tau -t)^{\alpha -1}h(\tau
)d\tau .
\end{equation*}
In this case he singled out a class of regular boundary conditions and
arrived at the following result.

\begin{theorem}[S.N.Kabanov]
\INDEX{Kab}{S.N.Kabanov}
Let the leading terms in the boundary conditions (\ref{eqC.7}) be regular,
$f\in L^p[-1,1],\quad
h_j\in L^q[-1,1],\quad \frac 1p+\frac 1q=1$.
Then for any $\delta \in (0,1/2)$ relation (\ref{eqB.15}) is satisfied.
\end{theorem}

For a general first order integro-differential operator
\begin{equation}
(E+N)(y^{\prime }+\tilde \alpha y)  \label{eqC.9}
\end{equation}
with general boundary condition (\ref{eqC.3}) it is also possible to
establish an equiconvergence.

\begin{theorem}[S.N.Kabanov]
\INDEX{Kab}{S.N.Kabanov}
\cite{kab90b}Assume that
\begin{enumerate}

\item[a)] $N(x,t)=\sum_1^3N_i(x,t),$ $N_i(x,t)$
are continuous for $t\leq x$ and $t\geq x,$%
\begin{equation*}
N_1(x,x-0)-N_1(x,x+0)=\tilde \alpha ,
\end{equation*}
\begin{equation*}
N_{2,t}^{\prime }(x,x-0)=\varphi (x)s(t),\quad \varphi (x)\in C[-1,1],\quad
s(t)\in L^q[-1,1],
\end{equation*}
\begin{equation*}
N_3(x,t)=v(x)h(t),\quad v(x)\in C[-1,1];
\end{equation*}

\item[b)]  for some $\alpha \in (0,1]\qquad \underset{-1}{\overset{1}{Var}}%
D_1^{-\alpha }h(t)<\infty ,\quad \underset{-1}{\overset{1}{Var}}D_1^{-\alpha
}\tilde h(t)<\infty $

\item[c)]  $D_1^{-\alpha }h(1)\cdot D_1^{-\alpha }\tilde h(1)\neq 0.$
\end{enumerate}
Then, if $h\in L^q{[-1,1]}$,
for any function
$f\in L^p[-1,1],\quad \frac 1p+\frac 1q=1,$
and any $\delta \in (0,1/2)$ relation (%
\ref{eqB.15}) is satisfied.
\end{theorem}

Necessity of conditions a) here is also justified. More precisely, it is
shown that the most general first order integro-differential operator with
boundary condition (\ref{eqC.3}) can be brought to the form (\ref{eqC.9})
where a) is fulfilled.

Next
\INDEX{Amv}{O.I.Amvrosova}O.I.Amvrosova \cite{amv84} studied a fractional
differentiation operator
\begin{equation}
l(y)=D^\alpha y=\frac d{dx}\int_0^x\frac{(x-)^{\alpha -1}}{\Gamma (\alpha )}%
y(t)dt,\qquad 0<\alpha <1,\quad x\in [-1,1]  \label{eqC.10}
\end{equation}
with boundary condition
\begin{equation}
U(y)=\int_{-1}^1\frac{k(t)y(t)}{(1-|t|)^{\beta +1}}dt=0,  \label{eqC.11}
\end{equation}
when
\begin{equation*}
\underset{-1}{\overset{1}{Var}}k(t)<\infty ,\quad 0<\beta +1\leq \alpha <1,
\end{equation*}
\begin{equation*}
k(0+0)\neq k(0-0),\quad k(-1+0)\cdot k(1-0)\neq 0.
\end{equation*}
Assuming the restrictions above to be fulfilled she obtained

\begin{theorem}[O.I.Amvrosova]
\INDEX{Amv}{O.I.Amvrosova}
Let $f\in L[-1,1],\quad D^\beta f(x)$ be absolutely continuous in $[-1,1]$.
Then
\begin{equation}
\lim_{r\longrightarrow \infty }\left\| \left( 1-|x|\right) |x|^{1+\gamma
}\left[ S_r(f)(x)-\sigma _r(f)(x)\right] \right\| _{C[-1,1]}=0
\label{eqC.12}
\end{equation}
where $\gamma $ stands for any positive number.
\end{theorem}

At last let us note an interesting
\INDEX{Sed}{A.M.Sedletskii}A.M.Sedletskii's article \cite{sed91} where he
investigated a uniform convergence of \ef expansions for a
differentiation operator with a
\INDEX{Sti}{T.Stieltjes}Stieltjes integral in boundary condition
\begin{equation*}
U(y)=\int_{-1}^1y(t)d\sigma (t)=0,
\end{equation*}
where
\begin{equation*}
d\sigma (t)=\frac{b(1-|t|)}{(1-|t|)^\alpha }k(t)dt
\end{equation*}
with a weakly oscillating function $b(t)$. This case is much more difficult
then the power singularity case, i.e. when
$b(t)\equiv 1$.
Hence here opens
an opportunity for a movement towards the most general form of boundary
conditions involving any kind of singularities at the end points.

Note that this subsection reproduces
\INDEX{Khro}{A.P.Khromov}
A.P.Khromov's review on equiconvergence
presented at the 7th Saratov winter school in 1994 \cite{hro95b} and is put
here under his kind permission.

\subsection{Asymptotic formulas for partial sums}

In the theorem \ref{ThA} we have no estimates of the rate of
equiconvergence. However, in 1967
\INDEX{KupN}{N.P.Kuptsov}N.P.Kuptsov \cite{nkup67a} has already indicated a
possibility of obtaining asymptotic formulas for the remainder term in that
theorem. Later he accomplished calculations in the second order case though
they have never been published (see remark in \cite[p.41]{osk73}). In this
case the asymptotic formulas are very complicated. Therefore in 1973
\INDEX{Os}{G.P.Os'kina}G.P.Os'kina established such formulas in every
subinterval $[\delta ,\pi -\delta ]\subset (0,\pi )$ in the article
\cite{osk73} accomplished under
\INDEX{KupN}{N.P.Kuptsov}N.P.Kuptsov's supervisorship. She considered
differential expression
\begin{equation*}
y^{(n)}+p_{n-2}y^{(n-2)}+\ldots +p_0y,\quad 0\leq x\leq \pi
\end{equation*}
and the simplest one, $y^{(n)}$,
with one and the same set of boundary
conditions (\ref{eqB.2}) but at the points $0$ and $\pi $.
Let $S_{r_k}(f)$
and $S_{r_k}^0(f)$ be the corresponding partial sums of \eaf
expansions and set
\begin{equation*}
Qy=l(y)-y^{(n)}.
\end{equation*}
Her main result reads as follows.

\begin{theorem}[G.P.Os'kina]
\label{ThOs}
Fix $\delta \in (0,1/2)$. Let $n$ be even and
$\delta \leq x\leq \pi -\delta$.
Then
\begin{equation}
S_{r_k}(f)-S_{r_k}^0(f)=\frac 1{\pi n}\int_0^\pi \int_0^\pi L_{r_k}(x,\xi
,t)p_{n-2}(\xi )f(t)dtd\xi +O\left( \frac 1{r_k}\right)   \label{eqOs.1}
\end{equation}
with an explicit though complicated expression for the kernel $L_{r_k}(x,\xi
,t)$:
\begin{multline}
L_r(x,\xi ,t)=\\
\int_r^\infty \frac 1\eta \left\{ \sin \eta \left[ \left|
x-\xi \right| +\left| \xi -t\right| \right] +\cos \eta \left| x-\xi \right|
\cdot I_1-\cos \eta \left| \xi -t\right| \cdot I_2\right\} d\eta
\label{eqOs.2}
\end{multline}
where we set for brevity
\begin{equation}
I_1=\sum_{j=k+1}^{3k-1}\overline{\varepsilon _j}\exp \left( \eta \varepsilon
_j\left| \xi -t\right| \right) ,\qquad n=4k,  \label{eqOs.3}
\end{equation}
\begin{equation}
I_2=\sum_{j=k+1}^{3k-1}\varepsilon _j\exp \left( \eta \varepsilon _j\left|
x-\xi \right| \right) ,\qquad n=4k.  \label{eqOs.4}
\end{equation}
In the case $n=4k+2$ one must take $3k$ as the upper bound in the sums (\ref
{eqOs.3})-(\ref{eqOs.4}) and replace there $\varepsilon _j$ by $\varepsilon
_{j+1/2}$.
\end{theorem}

\section{Equiconvergence and uniform minimality}

Eigenfunction systems may also be viewed of as an interesting and important
example of families of functions. During the last 20 years such
function-theoretic approach has been elaborated by
\INDEX{Il}{V.A.Il'in} V.A.Il'in and his school in numerous works,
see, for instance,
\cite{Il:spectral}. Below we shall briefly describe some of their results
about equiconvergence.

First, consider a function family
\begin{equation}
U=\left\{ u_k\right\} _{k=1}^\infty   \label{eqIl.1}
\end{equation}
and assume that they are \ef of the maximal operator generated in
$L^2(0,1)$ by the differential expression (\ref{eqNS.1})
with summable coefficients:
\begin{equation}
lu_k+\lambda _ku_k=\theta _ku_{k-1}, \qquad 0\leq x\leq1,  \label{eqIl.2}
\end{equation}
where the number $\theta _k$ takes two values, $0$ --- then $u_k$ is an
\ef, or $1$---then we require in addition that $\lambda _k=\lambda _{k-1}$
and it is an associated function$.$ Set $\theta _1=0$. In the case $n=1$ we
come to an exponential system
\begin{equation}
\left\{ \exp (i\lambda _kx)\right\} ,\quad 0\leq x\leq 1  \label{eqIl.3}
\end{equation}
--- a classical object of the function theory.

\subsection{A priori restrictions}

To facilitate an exposition let for simplicity $n$ be even and set
\begin{equation}
\mu _k:=\left( (-1)^{(n+2)/2}\cdot \lambda _k\right) ^{1/n}=\left( \varrho
e^{i\varphi }\right) ^{1/n}=\varrho ^{1/n}e^{i\varphi /n},\quad -\pi <\varphi \leq \pi .
\label{eqIl.4}
\end{equation}
Fix $p\geq 1$ and impose three a priori restrictions:

\begin{enumerate}
\item[A1)]  the system (\ref{eqIl.1}) is closed and minimal in $L^p(0,1)$,

\item[A2)]  $|\Im \mu _k|\leq C_1$

\item[A3)]  $\sum_{r\leq |\mu _k|\leq r+1}1\leq C_2,\qquad \forall r\geq 0.$
\end{enumerate}

Enumerate all \chv $\mu _k$ in the ascending modulus order.
Observe that
A1) yields existence of a unique biorthogonal system
\begin{equation}
\left\{ u_k^\prime\right\} _{k=1}^\infty ,\qquad u_k^\prime\in L^q(0,1),
\qquad \frac1p+\frac 1q=1.  \label{eqIl.5}
\end{equation}
Set
\begin{equation}
S_r(f)=\sum_{|\mu _k|\leq r}(f,u_k^\prime)u_k  \label{eqIl.6}
\end{equation}
Conditions A1)--A3) are assumed to be valid throughout the subsection
without further mentioning.

\begin{theorem}[V.A.Il'in]
\label{ThIl1}In order for the equiconvergence
\begin{equation}
\lim_{r\longrightarrow \infty }\left\| S_r(f)-
\left( V\sigma _r^\pi V^{-1}\right) (f)\right\| _{C(K)}=0
\label{eqIl.7}
\end{equation}
to be valid for any $f\in L^p(0,1)$ and any compact $K\subset (0,1)$ it is
necessary and sufficient that for any compact $K_0\subset(0,1)$%
\begin{equation}
\left\|
u_k\right\| _{L^p(K_0)}\cdot
\left\| u_k^\prime\right\| _{L^q(0,1)}
\leq C(K_0).  \label{eqIl.8}
\end{equation}
All the ingredients in the subtrahend in (\ref{eqIl.7}) are defined above in
(\ref{eqNS.3}), (\ref{eqR.3}).
\end{theorem}

In particular, equiconvergence holds provided the system (\ref{eqIl.1}) is
uniformly minimal (shortly, $U\in (UM)$ ), i.e. $K_0=[0,1]$
in (\ref{eqIl.8}).
This result is generalized to the matrix case \cite{dya96}
but we omit the statement due to its awkwardness.

\INDEX{Lom}{I.S.Lomov}I.S.Lomov \cite{lom95a} investigated equiconvergence
on the whole interval for two \ef subsystems of the form (\ref
{eqIl.1})--(\ref{eqIl.2}). To avoid a long introduction let us explain that
in the case of two self-adjoint operators $L_1,L_2$ in $L^2(0,1)$ he arrives
at an estimate
\begin{equation}
\left\| S_r(f,L_1)-S_r(f,L_2)\right\| _{C[0,1]}\leq C\left\| f\right\|
_{V[0,1]}  \label{eqLom.1}
\end{equation}
for any $f\in V[0,1]$. This estimate stems also from the convergence of both
series in question to some limiting values, see discussion in the beginning
of the chapter. However, methods developped in the aforementioned paper are
powerful and seem to serve well provided boundary conditions would be taken
into account.

\subsection{General series in eigenfunctions}

Given a family (\ref{eqIl.1}),(\ref{eqIl.2}) of \ef,
assume that conditions A2)--A3) are fulfilled and omit the first one: A1).
Consider a general series
\begin{equation}
\sum_{k=1}^\infty c_ku_k(x)  \label{eqIl.9}
\end{equation}
and assume that it converges to some summable function $f$ on subinterval $%
J\subset (0,1)$ in the weak sense:
\begin{equation}
\lim_{r\rightarrow \infty }\int_JS_r(x)\overline{\varphi (x)}dx=\int_Jf(x)%
\overline{\varphi (x)}dx  \label{eqIl.10}
\end{equation}
for any $\varphi (x)$ such that $\varphi ,l^{*}(\varphi )\in L^2(J).$ Here
\begin{equation*}
S_r(x):=\sum_{\left| \mu _k\right| \leq r}c_ku_k(x).
\end{equation*}
Assume also that
\begin{equation}
c_k\cdot \left\| u_k\right\| _{L^2(J)}\longrightarrow 0\text{ \quad }as\text{%
\quad k}\rightarrow \infty .  \label{eqIl.11}
\end{equation}

\begin{theorem}[A.M.Minkin]
\label{ThMgenser1}\cite{min80c}%
\INDEX{Min}{A.M.Minkin}Under the conditions (\ref{eqIl.10})-(\ref{eqIl.11})
the following relation holds
\begin{equation}
\lim_{r\longrightarrow \infty }\left\| S_r(x)-\sigma _r(f)\right\| _{C(K)}=0
\label{eqIl.12}
\end{equation}
for any compact $K\subset J$. Here $f$ is extended as $0$ to the whole
axis off $J$.
\end{theorem}

Compare two preceding theorems in the case $p_{n-1}(x)\equiv 0.$ In the
theorem \ref{ThMgenser1} completeness is omitted. Convergence in some weak
sense is used instead. Minimality isn't required at all. An analogue of
uniform minimality --- inequality (\ref{eqIl.8})--- is replaced by an
obviously necessary and very weak condition (\ref{eqIl.11}). The advantages
are evident: the domain of equiconvergence of a given \ef
expansion is easily to determine through the condition (\ref{eqIl.11}). Of
course, (\ref{eqIl.8}) with
$K_0=J$
yields (\ref{eqIl.11}) but not the converse.

Note also that
\INDEX{Moi}{E.I.Moiseev}E.I.Moiseev \cite{moi84,moi87}
has accomplished a deep investigation of concrete sine/cosine
and exponential  systems in $L^p(0,\pi), 1<p<\infty$:
\begin{enumerate}
\item  the system (\ref{eqIl.3}) with
$\lambda _k = k-\frac \beta2sign\ k,\ \ k\in \mathbb{Z}$;

\item  general systems
$G_{\beta,\gamma}:=\left\{
\sin \left((k+\beta /2)x+\gamma /2\right)
\right\} _{k=1}^\infty $

\item  sine systems
$S_\beta:=\left\{
\sin \left( (k+\beta /2)x \right)
\right\} _{k=1}^\infty $

\item  cosine systems
$C_\beta=1 \cup \left\{
 \cos \left( (k+\beta /2)x \right)
\right\} _{k=1}^\infty $

\end{enumerate}

He established delicate estimates of their biorthogonal systems using
difficult calculations. From his results it follows that
$S_\beta \in (UM)$ if $>\frac1p-2$ but is complete only
if $\beta\leq\frac1p$.

Analogously,
$C_\beta \in (UM)$ if $>\frac1p-1$ but is complete only
if $\beta\leq\frac1p+1$.
We omit formulation for the general system $G_{\beta,\gamma}$ to shorten
an exposition, see details in \cite{moi87}.

Note also the paper
\cite{ilmo92} of V.A.Il'in
\INDEX{Il}{V.A.Il'in}
and
E.I.Moiseev
\INDEX{Moi}{E.I.Moiseev}
where an important partial case ($ \beta=0$)
of the system $G_{\beta,\gamma}$ was considered
where the corresponding function family
happens to be a mixture of two sets of \ef of distinct \bvp .

Hence, V.A.Il'in's
\INDEX{Il}{V.A.Il'in}
theorem \ref{ThIl1}
holds for
any $f\in L^p(0,\pi)$ if the parameter $\beta$ is such
that the corresponding
system is complete and uniformly minimal.

However, for $p=2$ even a stronger result is valid.
Namely, equiconvergence with a trigonometric series expansion
holds for
any $f$ in the span of the corresponding system in $L^p(0,\pi)$
if the parameter $\beta$ is such that the corresponding
system is only uniformly minimal, see our theorem \ref{ThMgenser2}
below.

Note that actually theorem \ref{ThMgenser2} is true for any
$p,\ \ 1<p<\infty$
but this generalization needs \textit{almost orthogonality}
of the Birkhoff's \fss in $L^p(0,1)$ which is valid\INDEX{Birkhoff}{G.D.Birkhoff}
but we haven't yet published this result.

Theorem \ref{ThIl1} poses a natural question:

\emph{\ when and under what assumptions do conditions A1)--A3) hold? }

\noindent This problem seems to be very difficult.
Observe that these requirements are
obviously fulfilled for strongly regular two-point \bvp
as well as for
unconditional bases from exponentials (\ref{eqIl.3})
with the A3) condition being fulfilled.
Complete description of the latter systems  with
\textbf{arbitrary} exponents $\lambda_k$
is given in \cite{min91}.

For the beginning it would be important to
write down explicitly classes of irregular two-point \bvp
satisfying A1)--A3), at least separately.

A deep problem of uniform minimality of \eaf families
also deserves  a
separate treatment.
It seems that E.I.Moiseev's
\INDEX{Moi}{E.I.Moiseev}
results cited above
give a certain foundation for our conjecture about connection between
basicity and uniform minimality.

At first, let us introduce a distance between two normalized
in $L^p(0,1), \ 1<p<\infty,$ \ef families
$U:=\{ u_k \},\ \tilde{U}:=\{ \tilde{u}_k \}$,
satisfying (\ref{eqIl.2}):
\[
  d(U,\tilde{U}):=
\sup_k \left[
             \sum_{j=0}^{n-1}
		   \left(
                       |u_k^{[j]}(0)-\tilde{u}_k^{[j]}(0)|
		    \right)
		        +|\mu_k-\tilde\mu_k|
       \right].
\]

\begin{conjecture}
Call the system $U$ \textbf{stably} complete (incomplete)
if there exists a small $\varepsilon>0$ such that
completeness (incompleteness) preserves
for any other \ef
system $\tilde{U}$ with distance $\leq \varepsilon$.
Then
\begin{itemize}
\item An \ef family (\ref{eqIl.1}) is a basis in $L^p(0,1)$
if and only if
it is uniformly minimal and stably complete,
\item it is a basis in its span
if and only if
it is uniformly minimal and stably incomplete,
\item in the case $p=2$ basicity should be replaced by unconditional
basicity.
\end{itemize}
\end{conjecture}

Let us also indicate how to get around the condition (\ref{eqIl.10}).
\begin{theorem}[A.M.Minkin]\label{ThMgenser2}
\INDEX{Min}{A.M.Minkin}
Assume
that the family $U$ of \eaf (\ref{eqIl.1}),(\ref{eqIl.2})
is uniformly minimal
in $L^2(0,1)$, let
$U^\prime=\{ u_k^\prime \}$
be some its biorthogonal system which obeys an inequality
\[
  \| u_k \|_{L^2(0,1)} \cdot \| u_k^\prime \|_{{L^2(0,1)}}
\leq C,\ \ \forall k.
\]
 It is not important, whether $U$ is complete or not.
Denote $E=spanU \subset {L^2(0,1)}$, let
\[
  S_r(f):=\sum_{|\mu_k|\leq r} (f,u_k^\prime) u_k
\]
be a partial sum of the corresponding \ef expansion.
Then (\ref{eqIl.7}) remains true for any $f\in E$ and for any compact
$K\subset (0,1)$.
\end{theorem}
\begin{proof}
Fix a compact $K\subset (0,1)$
 and at first assume that $f$ is a linear combination of \eaf.
Then considerations in the theorem \ref{ThMgenser1}  yield
an upper estimate:
\[
 |S_r(f)(x)-\left( V\sigma _r^\pi V^{-1}\right) (f)(x)|
\leq
C_1\cdot \sup_k |(f,u_k^\prime)|\cdot \| u_k \|_{L^2(0,1)},\ \ x\in K.
\]
Here we need no assumption (\ref{eqIl.10}) because the series in question
reduces to a finite sum.
Next, using uniform minimality
we immediately derive that this difference constitutes a
uniformly bounded family of linear operators acting from
$L^2(0,1)$ to $C(K)$.
But on a dense (in $E$ ) subset of linear combinations of \ef
\[
 \left( V\sigma _r^\pi V^{-1}\right) (f) \longrightarrow f
\]
in $C(K)$ because \eaf are enough smooth. It remains to apply
the
\INDEX{Ban}{S.Banach}\INDEX{Steinh}{G.Steinhaus}Banach-Steinhaus theorem.
\end{proof}

\section{Singular self-adjoint operators}

\subsection{Self-adjoint expressions}

First recall that given a differential expression (\ref{eqNS.1}) in some
interval $G=(a,b)$ of the real axis its formally adjoint or
Lagrange-adjoint
\INDEX{Lag}{J.L.Lagrange}
is defined as follows
\begin{equation}
l^{*}(z)=z^{\{n\}}=Dz^{\{n-1\}}+\overline{p_0(x)}z^{\{0\}},  \label{eqSA.2}
\end{equation}
where
\begin{equation}
z^{\{0\}}=z,\quad z^{\{j\}}=Dz^{\{j-1\}}+\overline{p_j(x)}z^{\{0\}},\qquad
j=1,\ldots ,n-1.  \label{eqSA.3}
\end{equation}
It remains a \emph{differential }expression provided coefficients are enough
smooth
\begin{equation}
p_j^{(j)}\in L_{loc}^1(G).  \label{eqSA.4}
\end{equation}
Hence in the theory of \sa differential operators it is natural to
write the operation $l(y)$ in the form ($q=Entier(n/2)$ )
\begin{equation}
l(y)=D^ny+\sum_{\mu =0}^{q-1}\left( l_{2\mu }(y)+l_{2\mu +1}(y)\right) ;
\label{eqSA.5}
\end{equation}
\begin{equation*}
l_{2\mu }(y)=D^\mu \left( f_\mu (x)D^\mu y\right) ,\qquad l_{2\mu
+1}(y)=iD^\mu \left\{ Dg_{n-\mu }+\overline{g_{n-\mu }}D\right\} D^\mu y.
\end{equation*}
The coefficients $f_\mu (x),g_{n-\mu }$ are locally summable. Therefore (\ref
{eqSA.5}) is understood as a quasidifferential expression (q.d.) (see,
\cite[Ch.V]{Nai:ldo}). The most general form of
\sa expressions is
given below in chapter 3 (formula (3.0.1) there).
Further we assume that $l(y)$ is \sa,
$l^{*}=l$. Denote $L_{\min }$ the minimal symmetric operator
in $L^2(G)$ obtained by restricting $l(y)$ to the set of smooth enough
functions with compact support and then taking its closure.

\subsection{Classification}

Recall that in the theory of \sa operators a differential expression $l(y)$
is called $\emph{regular}$ (don't mix with the Birkhoff-regularity)
\INDEX{Birkhoff}{G.D.Birkhoff}
if the
following two conditions are fulfilled:

\begin{enumerate}
\item[i)]  the interval $G$ is finite;

\item[ii)]  all the coefficients of $l(y)$ are summable on $G$.
\end{enumerate}
Otherwise it is called \emph{singular}.

In the regular case one needs to add $n$ \sa boundary conditions at the end
points in order to define a \sa differential operator in $L^2(G)$. However,
it is now known that such \bvp  are Birkhoff-regular (see \cite{sal68,fie72}
\INDEX{Birkhoff}{G.D.Birkhoff}
if $n$ even and \cite{min77a} if $n$ is odd). A short proof which is valid
for all $n$ at once is given in \cite{min93a}. Hence, in the \sa case
equiconvergence obviously holds.
Now let us pass to singular expressions. First of all recall some basic
facts from the abstract spectral theory.

\subsection{Spectral function}

In the singular case the number $n_\lambda $ of solutions from $L^2(G)$ of
the equation
\begin{equation}
l(y)=\lambda y  \label{eqSA.6}
\end{equation}
is called a defect number. It is stable in the upper/lower half-plane
$
\mathbb{C}_{\pm },\ n_\lambda
\equiv n_{\pm },
\ \lambda \in \mathbb{C}_{\pm }
$.
The pair $(n_{+},n_{-})$ is called \emph{a defect index}
of the differential expression $l(y)$.

If they coincide, $n_{+}=n_{-}=m,$ then there exists a \sa extension $L$ in
$L^2(G)$ of the symmetric operator $L_{\min }$ defined by $m$ \sa boundary
conditions. However, the latter should be understood in some special sense
(see details in \cite[Ch.V]{Nai:ldo}).

Let $E_t,-\infty <t<\infty,\ E_t\leq E_s,\ t<s$ be the
corresponding resolution of identity. Normalize it as folows
\begin{equation*}
E_t=\frac 12\left( E_{t+0}+E_{t-0}\right) ,\qquad \lim_{t\rightarrow -\infty
}E_t=0,\quad \lim_{t\rightarrow \infty }E_t=I.
\end{equation*}
In the sequel these conditions are assumed to be fulfilled without further
mentioning.

The ortoprojector $E_t$ turns out to be an integral operator with a kernel
\begin{equation*}
\Theta (x,s,t),\quad x,s\in G,\qquad -\infty <t<\infty
\end{equation*}
see \cite[Chapter 13]{DunSchw:II}. There
exists an explicit formula for the spectral function $\Theta $ but it
shouldn't be needed in theorems' statements and therefore omitted. If
\begin{equation}
n_{+}\neq n_{-}  \label{eqSA.7}
\end{equation}
then there also exists a \sa extension $L^{+}$ of $L_{\min }$ but in some
ambient space $H^{+}$, containing $H=L^2(G)$ as a proper subspace. Let $%
E_t^{+}$ be the resolution of identity associated with $L^{+}$ and $P^{+}$
be the ortoprojector onto $H$ in $H^{+}$. Setting
\begin{equation*}
F_t:=P^{+}E_t^{+}|H
\end{equation*}
we obtain a nonorthogonal resolution of identity. $F_t$ is also
an integral operator
and its kernel $\Theta (x,s,t)$ is called a generalized spectral
function (\gspf), the same title is also applied to the family $F_t$
itself but usually this fact causes no ambiguity. The only difference
between $E_t$ and $F_t$ is that instead of the
\INDEX{Par}{A.M.Parceval} Parceval equality $F_t$ obeys the
\INDEX{Bessel}{F.W.Bessel}Bessel inequality
\begin{equation*}
\int_{\mathbb{R}}d_t(F_tf,f)\leq (f,f),\qquad \forall f\in H.
\end{equation*}
Inequality (\ref{eqSA.7}) may realize not only for odd $n$ but also
in the
even case provided the coefficents of $l(y)$ are complex-valued functions.

\subsection{Schr\protect\"odinger operator}
\INDEX{Schro}{E.Schr\protect\"odinger}

To begin with let $G=\mathbb{R},\quad n=2,\ \ q(x)$
be a real valued measurable function,
$q\in L_{loc}^1(G)$
and let
\[
l_4(y)=-y^{\prime \prime}+q(x)y.
\]
Let $\Theta _0(x,s,t)$ be a \spf, corresponding to the
zero potential, $q(x)\equiv 0$,\quad $G=\mathbb{R},$ namely
\begin{equation}
\Theta _0(x,s,t)=\left\{
\begin{tabular}[t]{ll}
$\frac 1\pi \frac{\sin r(x-s)}{x-s},$ & $r=\sqrt{t},\quad t>0,$ \\
$0$ & $t\leq 0.$%
\end{tabular}
\right.  \label{eqShL.1}
\end{equation}

\begin{theorem}[B.M.Levitan]
\label{ThLM1}%
\INDEX{Lev}{B.M.Levitan}%
Fix a compact $K\subset \mathbb{R}$. Then
\begin{equation}
\lim_{t\longrightarrow +\infty }\left[ \Theta (x,s,t)-\Theta
_0(x,s,t)\right] =0  \label{eqShL.2}
\end{equation}
uniformly with respect to $x,s\in K$.
\end{theorem}

Next, consider a half-bounded interval $G=[0,\infty )$
and assume that a real
valued function $q\in L[0,a]$ for any $a>0$.
Since the potential is real the
defect index equals $(1,1)$ or $(2,2)$. In the first case add to
equation
\begin{equation}
l_4(y)=\lambda y  \label{eqShL.3}
\end{equation}
a boundary condition at the origin
\begin{equation}
y^{\prime }(0)=hy(0),\quad h\neq \infty ,\quad h\text{ real.}
\label{eqShL.4}
\end{equation}
In the second case add in addition a boundary condition at infinity (we omit
details).

Denote $\Theta ^h(x,s,t)$ the \spf of the problem (\ref{eqShL.3})--(\ref
{eqShL.4}) and let $\Theta _1^h(x,s,t)$ be a \spf of the same problem with
a zero potential.

\begin{theorem}[B.M.Levitan, V.A.Marchenko]
\label{ThLM2}%
\INDEX{Lev}{B.M.Levitan}%
\INDEX{Mar}{V.A.Marchenko} Let $b>0$. Then uniformly with respect to $x,s\in
[0,b]$%
\begin{equation}
\lim_{t\longrightarrow +\infty }\left[ \Theta ^h(x,s,t)-\Theta
_1^h(x,s,t)\right] =0,  \label{eqShL.5}
\end{equation}
\begin{equation}
\lim_{t\longrightarrow +\infty }\left[ \Theta _1^h(x,s,t)-\Theta
_0(x,s,t)\right] =I_1+I_2+I_3,  \label{eqShL.6}
\end{equation}
where we set
\begin{equation*}
I_1=\frac h\pi \int_0^\infty \frac{\sin \nu (x+s)}{\nu ^2+h^2}\nu \ d\nu ,
\end{equation*}
\begin{equation*}
I_2=-\frac{2h^2}\pi \int_0^\infty \frac{\cos \nu (x-s)}{\nu ^2+h^2}\ d\nu ,
\end{equation*}
\begin{equation*}
I_3=\left\{
\begin{tabular}{ll}
$0,$ & $h\geq 0,$ \\
$h^2e^{h(x+s)},$ & $h<0.$%
\end{tabular}
\right.
\end{equation*}
\end{theorem}
When $h=0$ the boundary condition (\ref{eqShL.4}) should be understood as
\begin{equation}
y(0)=0.  \label{eqShL.7}
\end{equation}
Together with $q(x)\equiv 0$ it corresponds to the \spf
\begin{equation}
\Theta _1^\infty (x,s,t)=\frac 2\pi \int_0^r\sin \nu x\cdot \sin \nu s\ d\nu
,\qquad t>0  \label{eqShL.8}
\end{equation}
and $\Theta _1^\infty \equiv 0,\quad t\leq 0.$

\begin{theorem}[B.M.Levitan, V.A.Marchenko]
\label{ThLM3}%
\INDEX{Lev}{B.M.Levitan}%
\INDEX{Mar}{V.A.Marchenko}Let $b>0$. Then uniformly with respect to $x,s\in
[0,b]$%
\begin{equation}
\lim_{t\longrightarrow +\infty }\left[ \Theta ^\infty (x,s,t)-\Theta
_1^\infty (x,s,t)\right] =0  \label{eqShL.9}
\end{equation}
\end{theorem}

Conditions of the theorems \ref{ThLM1}--\ref{ThLM3} means exactly
equiconvergence for $\delta _s$--- the
\textbf{delta-function}
at the point $s$.
Moreover, uniformity of their statements with respect to $s\in [0,b]$ means
that equiconvergence also holds for measures with finite support
\begin{equation}
\lim_{r\longrightarrow \infty }
\left\| S_r(d\mu )- S_r^0 (d\mu)\right\| _{C[0,b]}=0,
\qquad \forall b>0.
\label{eqShL.10}
\end{equation}
Here
$S_r^0=\int_{-r^n}^{r^n}dE_t^0$
where the resolution of identity  $E_t^0$
corresponds to the case of zero potential. In particular, it is possible to
take $f\in L[0,b],\ \ f\equiv 0,\quad x>b$ in (\ref{eqShL.10})  instead of
$d\mu $.

Note that the same assertion is valid for $f\in L^2(G)$. Concretely,

\begin{theorem}[B.M.Levitan, V.A.Marchenko]
\label{ThLM4}%
\INDEX{Lev}{B.M.Levitan}%
\INDEX{Mar}{V.A.Marchenko}Let either $G=\mathbb{R}$ or $G=[0,\infty )$ and in
the latter case the boundary conditions (\ref{eqShL.4}) or (\ref{eqShL.7})
are added to the equation (\ref{eqShL.3}).
Then for any compact $K\subset G$,
the regular end point $0$ included (if present)
\begin{equation}
\lim_{r\longrightarrow \infty }\left\| S_r(f)-S_r^0(f)\right\| _{C(K)}=0.
\label{eqShl.11}
\end{equation}
\end{theorem}

\begin{remark}
Observe that in the theorems \ref{ThLM2}--\ref{ThLM4} the difference of
\spf or \ef expansions vanishes uniformly up to one
(regular) end point. This is much more subtle fact than the question of
equiconvergence on internal compacts which may also be established by
the methods of the articles \cite{lev53,mar55} as well.
\end{remark}

As far as we learned from our elder colleagues
an equiconvergence theorem for singular \sa operators was also
obtained by
\INDEX{KupN}{N.P.Kuptsov}N.P.Kuptsov independently of
\INDEX{Lev}{B.M.Levitan}%
\INDEX{Mar}{V.A.Marchenko}B.M.Levitan and V.A.Marchenko.
However, his proof has never been published and seems to be
irrevocably lost.

\subsection{Higher order}

Constructions of
\INDEX{Lev}{B.M.Levitan}%
\INDEX{Mar}{V.A.Marchenko}B.M.Levitan and V.A.Marchenko extensively used the
theory of hyperbolic equations as well as the theory of the transformation
operators. Therefore it was impossible to transform them directly
to the high-order
case. Only ten years later
\INDEX{Kost}{A.G.Kostuchenko }A.G.Kostuchenko made a breakthrough and
succeeded to generalize their results to the even order \sa operators \cite
{kos:dr,kos68}. Unlike \cite{lev53,mar55} he bypassed a use of the
transformation operators theory which are \emph{unbounded} for $n>2$.
Instead, he applied the theory of parabolic equations
\begin{equation*}
\frac{\partial u}{\partial t}=l(u),
\end{equation*}
where $l$ is a \sa differential expression of the form (\ref{eqB.1}).

\begin{theorem}[A.G.Kostuchenko]
\label{ThKos1}%
\INDEX{Kost}{A.G.Kostuchenko }Let $G=\mathbb{R}$ and $L_{\min }$ be defined by
a \sa expression (\ref{eqB.1}) of order $n=2q,\ q>1$ with real
coefficients. In addition, assume that $L_{\min }$ is half-bounded and
\begin{equation}
p_{n-2}(x)\text{ is a piece-wise smooth function,}\ \
p_j\in L_{loc}^1(G),\ \ 0\leq j<n-2.  \label{eqShL.12}
\end{equation}
\end{theorem}

Let $L$ be its \sa extension in $L^2(G)$ with a \spf $\Theta (x,s,t)$.
Set $r:=t^{1/n},\ t>0$.
Then the theorem's \ref{ThLM1} assertion is valid.
\begin{theorem}[A.G.Kostuchenko]
\label{ThKos2}%
\INDEX{Kost}{A.G.Kostuchenko }
Let $G=[0,\infty )$ and all other theorem's
\ref{ThKos1} assumptions be satisfied.
Assume also that $0$ is a regular
end point, i.e.
$p_j\in L[0,a),\ \forall a>0$
and the defect index
equals $(q,q)$. Then any \sa extension $L$ of $L_{\min }$ is generated by
$q$
\sa boundary conditions at $0$ \cite[pp.212--214]{Nai:ldo}:
\begin{equation}
B_j(y)\equiv \sum_{i=1}^nb_{ij}y^{(j-1)}(0)=0,\qquad j=1,\ldots ,q;
\label{eqShL.13}
\end{equation}
\begin{equation}
\sum_{j=1}^qb_{ij}\overline{b_{k,n+1-j}}-\sum_{j=1}^qb_{i,n+1-j}\overline{%
b_{kj}}=0,\quad i,k=1,\ldots q.  \label{eqShL.sa}
\end{equation}
\INDEX{Kost}{A.G.Kostuchenko }A.G.Kostuchenko also assumes these forms to
obey two additonal complicated restrictions which we shall omit for
simplicity.
Denote $\Theta _1(x,s,t)$ a \spf of a model \sa operator
generated by the expression $D^ny$ and \sa boundary conditions (\ref
{eqShL.13}). Then uniformly with respect to $x,s\in [0,b]$%
\begin{equation}
\lim_{t\longrightarrow +\infty }\left[ \Theta (x,s,t)-\Theta
_1(x,s,t)\right] =0  \label{eqShL.14}
\end{equation}
for any $b>0$.
\end{theorem}
He also established equiconvergence for square summable functions. The
statement repeats that of the theorem \ref{ThLM4} and therefore is omitted
here.

It is difficult to underestimate the importance of his contribution.
However, observe that self-adjointess of $l(y)$ yields some smoothness of
the coefficients in addition to (\ref{eqShL.12}) as is readily seen from
(\ref{eqSA.2})-(\ref{eqSA.5}).
Moreover,
\INDEX{Kost}{A.G.Kostuchenko }A.G.Kostuchenko treats the general case of
nonhalf-bounded $L_{\min }$ by its squaring. This operation requires
existence of additional $n$ derivatives of each of the coefficients $p_j$.

At the beginning of eighties we obtained equiconvergence theorems for
singular \sa high-order equations \emph{without any unnecessary a priori
restrictions}.

\begin{theorem}[A.M.Minkin]
\cite{min:phd}\label{ThMsing1}%
\INDEX{Min}{A.M.Minkin}Let $G$ be a finite or infinite interval
of $\mathbb{R}$%
, $l(y)=y^{[n]}$ be a general quasi-differential \sa
expression of the form
\begin{eqnarray}
y^{[0]}=y,\quad y^{[1]}=Dy^{[0]},  \label{eqShL.15} \\
y^{[j]}=Dy^{[j-1]}+\sum_{k=0}^{j-2}p_{j-1,k}(x)y^{[k]},\quad j=2,\ldots ,n
\notag
\end{eqnarray}
with complex-valued coefficients such that
\begin{equation}
p_{i,k}\in L_{loc}^1(G),\qquad p_{i,k}(x)=\overline{p_{n-1-k,n-1-i}(x)}.
\label{eqShL.16}
\end{equation}
Let $L_{\min }$ be the corresponding minimal symmetric operator in $L^2(G)$
and $F_t$ be some its \gspf. For instance, $F_t$ may be a restriction of a
\spf of some \sa extension of $L_{\min }$ in a larger interval $G_1\supset
G$. Set
\begin{equation*}
\overset{0}{L^1}(G):=\left\{ f\in L^1(G)\quad |\quad f\equiv 0\quad \text{%
near the boundary}\right\}
\end{equation*}
Let $f\in \overset{0}{L^1}(G)\cup L^2(G),\quad g\in L^2(G)\cup L^1(\mathbb{R})$
and assume that $f(x)\equiv g(x)$ for almost all $x\in \Omega =(a,b)\subset G
$. Then
\begin{equation}
\lim_{r\rightarrow \infty }\left\| \int_{-r^n}^{r^n}dF_tf-\sigma
_r(g)\right\| _{C(K)}=0  \label{eqShL.17}
\end{equation}
for any compact $K\subset \Omega $.
\end{theorem}

Notice that (\ref{eqShL.17}) means simultaneously equiconvergence and
localization. The case $f\in L^2(G)$ appeared for the first time in \cite
{min80b}. Sketch of the proof is published in \cite{min80a} provided that
the operator measure $dF_t$ is discrete. Proof of the general case closely
follows the lines of the discrete one, see \cite{min:phd}.

For $n=2q$ and $G=[0,\infty )$ this theorem's statement may be improved
in order to include the end point $0$. Namely, require that
\begin{equation*}
p_{i,k}\in L[0,a),\qquad \forall a>0
\end{equation*}
and assume that the \sa expression $l(y)=y^{[n]}$ has
a defect index $(m,m)$
in $L^2(G),\quad q\leq m\leq n.$ Since the coefficients are complex this is
really a requirement
(see, \cite[p.175]{AkhGl:II}). A deep investigation of
defect indices for general symmetric systems of singular differential
equations has been accomplished in an article
of
\INDEX{Kog}{V.I.Kogan} V.I.Kogan and
\INDEX{Rofe}{F.S.Rofe-Beketov}F.S.Rofe-Beketov \cite{kogro75}.
We refer the
interested reader to it for more information and details.

\begin{theorem}[A.M.Minkin]
\label{ThMsing2}%
\INDEX{Min}{A.M.Minkin}Let $L$ be a \sa operator in $L^2[0,\infty ),$
defined by (\ref{eqShL.15}) and $m$ boundary conditions at least $q$ from
which are given at the origin like (\ref{eqShL.13}). Fix $b>0$ and let $f\in
L[0,b],\quad f(x)\equiv 0,\ x>b$. Then for any $\varepsilon \in (0,b)$%
\begin{equation}
\lim_{r\rightarrow \infty }\left\| S_r(f,L)-S_r(f,L_b)\right\|
_{C[0,b-\varepsilon ]}=0  \label{eqShL.18}
\end{equation}
where $L_b$ stands for an ordinary \sa differential operator in $L^2[0,b]$
generated by $l(y)$ and decomposing \sa boundary conditions $q$ of which
coincide with (\ref{eqShL.13}) and the other $q$ are taken at the
end point $b$.
\end{theorem}

Let us stress the fact that in the theorem \ref{ThMsing1} the
coefficients' requirements are the least possible (see, (\ref{eqShL.16}) ).
Moreover, it covers at once q.d. equations of arbitrary order as well as
\sa extensions going outside the space $L^2(G)$. Earlier these questions
haven't been considered at all.
In addition theorem \ref{ThMsing2}
removes unnecessary restrictions on the coefficients of the boundary forms
$B_j(y)$ imposed in \cite{kos68}.

Later
\INDEX{Imam}{V.I.Imamberdiev } V.I.Imamberdiev established a \spf
asymptotics of an odd order \sa operator developing the parabolic
equations method, see \cite{ima93}.

\subsection{Kato condition}

Looking closely at the theorems \ref{ThLM2}--\ref{ThLM4} as well as
at theorems \ref{ThKos2}, \ref{ThMsing2}
one sees that they doesn't constitute a \emph{full }%
generalization of
\INDEX{Tam}{J.Tamarkin}Tamarkin's theorem \ref{ThB}. Indeed, only \emph{one
end point} is included! Therefore it is quite natural to pose a question
of whether these results are valid throughout the \emph{whole infinite
interval} $G$.

Even in the second-order case this problem remains open in the case of
arbitrary potential $q(x)$. However, recently this important question has
been answered in affirmative in a series of articles due to
\INDEX{Il}{V.A.Il'in}V.A.Il'in,
\INDEX{An}{I.Antoniu}I.Antoniu and
\INDEX{Kri}{L.V.Kritskov}L.V.Kritskov \cite{ilan95}--\cite{ilkri95}.

They considered a class of potentials in $\mathbb{R}$ satisfying
\INDEX{Kato}{T.Kato}\textbf{Kato} condition:
\begin{equation}
\sup_{-\infty <x<\infty }\int_x^{x+1}|q(s)|ds\leq C.  \label{eqKato}
\end{equation}
Let us state their results.

\begin{theorem}[V.A.Il'in]
\cite{il95a}\label{ThKato1}%
\INDEX{Il}{V.A.Il'in}Consider a \sa
\INDEX{Schro}{E.Schr\protect\"odinger} operator $L$ in $L^2(\mathbb{R})$
with potential $q(x)$ satisfying
\INDEX{Kato}{T.Kato} Kato condition.
Let $\Theta (x,s,t)$ be its \spf and define $%
\Theta _0(x,s,t)$ as in (\ref{eqShL.1}). Then there exists $T>0$ such that
for some finite constant $C_T$%
\begin{equation}
\sup_{t\geq T}\sup_{x,s\in \mathbb{R}}\left| \Theta (x,s,t)-\Theta
_0(x,s,t)\right| =C_T.  \label{eqKato.1}
\end{equation}
In addition if $1\leq p\leq 2,\quad f\in L^p(\mathbb{R})$ then
\begin{equation}
\lim_{r\longrightarrow \infty }\left\| S_r(f)-\sigma _r(f)\right\| _{C(\mathbb{R%
})}=0.  \label{eqKato.2}
\end{equation}
\end{theorem}

This theorem was preceded by result due to
\INDEX{An}{I.Antoniu}I.Antoniu and
\INDEX{Il}{V.A.Il'in}V.A.Il'in \cite{ilan95} where the
\INDEX{Hill}{G.Hill}Hill operator ($q(x)$
is continuous periodic function on $\mathbb{R})$ was investigated.
Afterwards the
theorem \ref{ThKato1} was carried over to the
\INDEX{Schro}{E.Schr\protect\"odinger}Schr\"odinger operator with a matrix potential $%
Q(x)$ satisfying
\INDEX{Kato}{T.Kato}Kato condition \cite{kur96}.
\INDEX{Kur}{A.V.Kurkina} A.V.Kurkina proved an inequality
\begin{eqnarray*}
\sup_{t\geq T}\sup_{x,s\in \mathbb{R}}\sum_{k=1,k\neq j}^m\left\{ \left| \Theta
_{jk}(x,s,t)\right| +\left| \Theta _{jj}(x,s,t)-\Theta _0(x,s,t)\right|
\right\} &\leq &C_T<\infty , \\
j &=&1,\ldots ,m
\end{eqnarray*}
for the components $\Theta _{jk}(x,s,t)$ of the \spf $\Theta (x,s,t)$.

Further
\INDEX{Il}{V.A.Il'in}V.A.Il'in and
\INDEX{An}{I.Antoniu}I.Antoniu investigated the so called liouvillian,
generated by a \sa
\INDEX{Schro}{E.Schr\protect\"odinger}Schr\protect\"odinger operator
satisfying
\INDEX{Kato}{T.Kato}Kato condition \cite{ilan96}. It is important for physical
applications but we have to omit the statement because it requires
introducing a lot of preliminary notions.

Of course, it would be important to generalize these results
to higher orders
as well as to clarify necessity of the
\INDEX{Kato}{T.Kato}Kato condition in the question
of equiconvergence on the whole interval $G$.

\section{Multidimensional Schr\protect\"odinger-type operator}
\INDEX{Schro}{E.Schr\protect\"odinger}

During the past 20 years the author developped a rather general approach to
equiconvergence problems. It will be discussed in other chapters. Here we
bring to the reader's attention its application to operators in partial
derivatives obtained in a joint article with
\INDEX{Shu}{L.A.Shuster}L.A.Shuster \cite{mishu91}.

Let $D=[-1,1]^m,\ m>1,\ q(x)$ be a real valued summable function in $D$.
Set
\begin{equation*}
Ly=(-\Delta )^ny+q(x)y,\qquad \Delta =\frac \partial {\partial x_1^2}+\ldots
+\frac \partial {\partial x_m^2}.
\end{equation*}
At first $L$ is defined on trigonometric polynomials
\begin{equation*}
e_s(x):=\exp (i\pi \left\langle s,x\right\rangle ),\quad \left\langle
s,x\right\rangle =\sum_{j=1}^ms_jx_j,\quad s\in \mathbb{Z}^m.
\end{equation*}
Under restriction $2n>m$ there exists its
\INDEX{Fried}{K.O.Friedrichs}Friedrichs extension
which will
also be written as $L$
and happens to be a half-bounded \sa operator in $L^2(D)$,
satisfying
periodic boundary conditions
$y(x+2s)=y(x),\ \ s\in \mathbb{Z}^m$.
Denote its
spectrum $\sigma (L)$ and let $\varphi (k)$ be the number of
integer solutions
from $\mathbb{Z}^m$ of the equation
\begin{equation}
|s|^2=k,\text{ }k\text{ natural,}  \label{eqMiShu2}
\end{equation}
where $|s|^2:=\left\langle s,s\right\rangle $. We shall need the following

\begin{definition}
For any continuous $2$-periodic function on $D$ set
\begin{equation*}
\hat{f}(s)
:=\int_Df(x)\cdot \overline{e_s(x)}dx,
\end{equation*}
\begin{equation*}
\left\| f\right\| _A:=\sum_{s\in \mathbb{Z}^m}\left|
\hat{f}%
(s)\right| .
\end{equation*}
Introduce a space $A$ of all absolutely convergent series on
$\mathbb{R}^m/(2\mathbb{Z})^m$
as the subspace of continuous
$2$-periodic function $f$ on $D$
such that the $A$-norm is finite, see \cite{Kah:abs}.
\end{definition}

It is known that in absence of the potential $q(x)$ the
\INDEX{Laplace}{P.Laplace}Laplace operator $%
-\Delta $ on the torus $\mathbb{R}^m/(2\mathbb{Z})^m$ has \ev of high
multiplicity. In the presence of potential, generally speaking, such a
multiple \ev splits into a group of neighboring \ev. The
theorem below is due to
\INDEX{Shu}{L.A.Shuster}L.A.Shuster and gives a rigorous description
of this phenomenon.

\begin{theorem}[L.A.Shuster]
Take $a\in (0,\frac 12\pi ^2)$ and let $2n>m+3.$ Then

\begin{enumerate}
\item  there exists an integer $k(a)$ such that
\begin{equation}
Card
\left\{
\lambda \in \sigma (L)\ |\ \left| \lambda ^{1/n}-k\pi
^2\right| \leq a
\right\}
=\varphi(k),  \forall k\geq k(a).\label{eqMiShu3}
\end{equation}

\item  Denote $H(k)$ the spectral subspace of $L$ spanned by all
\ef with \ev in the cluster (\ref{eqMiShu3}). If in
addition $n>m+1$ than there exists a basis $\left\{ h_j(x)\right\}
_{j=1}^{\varphi (k)}$ in $H(k)$ such that
\begin{equation}
h_j(x)=\exp (i\pi \left\langle s,x\right\rangle )+O\left(
k^{-(n-m-1)}\right) ,  \label{eqMiShu4}
\end{equation}
where $s\in \mathbb{Z}^m$ runs over all solutions of the equation (\ref
{eqMiShu2}). The symbol $O()$ is understood here in the $A$-norm sense and
the constants are absolute.
\end{enumerate}
\end{theorem}

Under this theorem's assumptions let $P_k$ be the ortoprojector onto $H(k)$
in $L^2(D)$ and set
\begin{equation*}
\tau _r(f):=\sum_{k(a)\leq k\leq r^2}P_kf.
\end{equation*}
Denote $P_0$ an ortoprojector onto the set of all \ef of $L$
corresponding to a finite number of first \ev $\lambda $ such that $%
\lambda <\left( k(a)\pi ^2-a\right) ^n$ and set $S_r(f):=P_0(f)+\tau _r(f).$
In addition denote
\begin{equation*}
\sigma _r^\pi (f)=\sum_{|s|\leq r}(f,e_s)e_s
\end{equation*}
the $r$-th partial sum of a multiple trigonometric Fourier series.
\INDEX{Fou}{J.Fourier}

\begin{theorem}[A.M.Minkin]
\INDEX{Min}{A.M.Minkin}Assume that $2n>m+3.$ Then $\forall f\in L^2(D)$%
\begin{equation}
\lim_{r\rightarrow \infty }r^{2n-2m-3/2}\left\| S_r(f)-\sigma _r^\pi
(f)\right\| _A=0.  \label{eqMiShu5}
\end{equation}
\end{theorem}

Obviously, the theorem's assertion claims equiconvergence with rate provided
that $2n>m+3$ and divergence with rate, otherwise. Moreover, it is
established in the $A$-norm which is principally stronger than the $C$-one.
Note also that equiconvergence holds for \emph{nonsmooth, namely, square
summable} functions. Recall that known results requires considerable
order of the Riesz
\INDEX{Riesz}{F.Riesz}
typical means of the function in question,
(see \cite[p.70--76]{alilnik76} ).

It seems to us natural to join together \ef corresponding to the
same cluster. However, it is likely that this idea hasn't been employed
earlier in the spectral theory of operators in partial derivatives. We think
that strong smoothness requirements on the function $f$ usual in that theory
stem from the fact that one tries to obtain convergence of the series itself
and not of an appropriate series \emph{with brackets}, instead.

\section{General equiconvergence principles}

In this section we briefly outline several general approaches for the
equiconvergence problem.

\subsection{Iteration of the resolvent's equation}

Given a
\INDEX{Ban}{S.Banach}Banach space $B$ with a norm
$\left\| {\ \  }\right\| _B$,
a dense
lineal $D\subset B$ endowed with a second norm
$\left\| {\ \  }\right\| _D$,
consider linear operators $A$ and $Q$ mapping $D$ into $B$.
Note that $D$ isn't
necessarily closed in the norm $\left\| {\ \  }\right\| _D$.
Consider a family of
concentric circonferences $C_n$ centered at the origin with
radii
$r_n\longrightarrow \infty ,\quad
r_1\leq r_2\leq \ldots$.

Set
\[
\alpha _n:=\max_{\lambda \in C_n}\left\| QR_\lambda (A)\right\|
_{B\rightarrow B}
\]
and
\[
\beta _n:=
\max_{\lambda \in C_n}\sup_{f\in B,\ f\neq 0}
\left\|
R_\lambda (A)f
\right\| _D /
\| f \|_{B}.
\]
In \cite{nkup67a}
\INDEX{KupN}{N.P.Kuptsov}N.P.Kuptsov established the following general
theorem.

\begin{theorem}[N.P.Kuptsov]
\label{ThKupgen}Assume that $\alpha _n\longrightarrow 0,\ \
n\longrightarrow \infty$.
Then there exists a natural $N$ such that for $%
n\geq N\quad \exists R_\lambda (A+Q)$.

If in addition $\alpha
_n\longrightarrow 0$ and $\beta _n\alpha _n=O\left( 1/r_n\right)$,
then for
any $f\in B$%
\begin{equation*}
\left\| \frac 1{2\pi i}\int_{C_n}\left( R_\lambda (A+Q)f-R_\lambda
(A)f\right) d\lambda \right\| _D=o(1),\qquad n\longrightarrow \infty .
\end{equation*}
\end{theorem}

Take for instance $A=D^n$ in $[0,1]$ with regular two-point boundary
conditions and let
$B=L[0,1],\quad D=D_A,\quad
\left\| y\right\|
_D:=\left\| y\right\| _{C[0,1]}$.

Set
\begin{equation*}
Qy=\int_0^1D^{n-2}y(t)d_t\sigma (x,t)
\end{equation*}
where $\underset{0}{\overset{1}{Var}}_x\sigma (x,t)=:q(x)\in L[0,1]$ and the
theorem \ref{ThKupgen} applies.

The proof rests on the formula
\begin{equation*}
-\frac 1{2\pi i}\int_{C_n}
\bigl(
R_\lambda (A+Q)f-R_\lambda (A)f
\bigr)
d\lambda
=\frac 1{2\pi i}\int_{C_n}R_\lambda (A)QR_\lambda (A)fd\lambda +V_nf
\end{equation*}
with the remainder's estimate
\begin{equation*}
\left\| V_nf\right\| _D=O\left( \alpha _n\cdot \left\| f\right\| _B\right),
\end{equation*}
which arises after iterating one time the identity connecting resolvents of
the main and perturbed operators.

\subsection{Commutator approach}

In \cite{min80c} and in subsequent papers we developped a new machinery in
order to handle the equiconvergence problems. Below we shall illustrate it
briefly in the simplest situation.

Let $D^{\prime }(\mathbb{T})$ be the space of generalized functions on the
one-dimensional torus
\[
\mathbb{T}:\mathbb{=R}\backslash \mathbb{Z},\qquad
D^{\prime}(\mathbb{T})=\left(C^\infty (\mathbb{T})\right)^\prime.
\]
For any $F\in D^{\prime }(\mathbb{T})$ set $%
\hat{F}(l):=F(e_{-l}),\ \
e_{l}=\frac 1{\sqrt{2\pi }}\exp (ilx)$,
\begin{equation*}
S_r(F):=\sum_{\left| l\right| \leq r/2\pi }\hat{F}(l)\cdot
e_{l}.
\end{equation*}
Introduce a space $PF(\mathbb{T})$ of pseudofunctions on $\mathbb{T}$ as a
subspace of $D^{\prime }(\mathbb{T})$ with vanishing Fourier
\INDEX{Fou}{J.Fourier}
coefficients as $%
\left| l\right| \longrightarrow \infty ,$%
\begin{equation*}
\left\| F\right\| _{PF}=\sup_l\left| \hat{F}(l)\right| .
\end{equation*}
Let $K$ be a compact in $\mathbb{T}$. Introduce a seminorm in $C(K)$:
\begin{equation*}
\left\| f\right\| _{A(K)}:=\inf_{g|_K\equiv f}\left\| g\right\| _A.
\end{equation*}
Denote $A(K)$ the lineal in $C(K)$ with finite seminorm
$\| {\ \ } \|_{A(K)}$.

\begin{theorem}[A.M.Minkin, localization principle]
\label{ThMinloc} Given $F\in D^{\prime }(\mathbb{T}%
),\ F|\Omega =0$
for some open set $\Omega \subset \mathbb{T}$.
Then
\begin{equation*}
\lim_{r\rightarrow \infty }\left\| S_r(F)\right\| _{A(K)}=0
\end{equation*}
for any compact $K\subset \Omega $.
\end{theorem}

Its proof relies on a proposition which generalizes one theorem
which goes back to
\INDEX{Raih}{A.Rajchman}A.Rajchman \cite[p.194--196]{Bar:ts}.

\begin{lemma}
Let $\gamma \in C^1\left( \mathbb{
T}\right) $ be such that $\gamma ^{\prime
}\in A$. Then for any $F\in PF$
\begin{equation*}
\lim_{r\rightarrow \infty }\left\| \left[ S_r,\gamma \right] (F)\right\|
_A=0,
\end{equation*}
where $[\ ,]$ stands for the operators' commutator.
\end{lemma}

\begin{proof} Obviously,
\begin{equation*}
\widehat{\left( \gamma \cdot F\right)}(l)=
\sum_{j+k=l}\hat{F}(k)\cdot\hat{\gamma}(l).
\end{equation*}
Hence,
\begin{equation*}
S_r\left( \gamma \cdot F\right) =
\sum_{l,k\in P_1}\hat{F}(k)%
\hat{\gamma }(l-k)\cdot e_l(x),
\end{equation*}
where
$P_1=\left\{ (l,k)\ \ |\ \ |l|\leq r/2\pi ,\ \ -\infty <k<\infty
\right\} $.
Conversely, expanding $\gamma (x)$ into an absolutely convergent
trigonometric Fourier
\INDEX{Fou}{J.Fourier}
series we arrive at identity
\begin{equation*}
\gamma \cdot S_r\left( F\right) =
\sum_{l,k\in P_2}\hat{F}(k)%
\hat{\gamma }(l-k)\cdot e_l(x),
\end{equation*}
where $P_2=\left\{ (l,k)\ \ |\ \ |k|\leq r/2\pi ,\ \ -\infty <l<\infty
\right\} $. Then
\begin{equation}
S_r\left( \gamma \cdot F\right) -\gamma \cdot S_r\left( F\right)
=\sum
_{(l_1,l_2)\in P_1\setminus P_2}
-
\sum _{(l_1,l_2)\in P_2\setminus P_1}.  \label{eqGP.1}
\end{equation}
Clearly,
$P_1\setminus P_2=\left\{ (l,k)\ \ |\ \ |l|\leq r/2\pi
<k\right\} $
and
$P_2\setminus P_1=\left\{ (l,k)\ \ |\ \ |k|\leq r/2\pi
<l\right\}$.
But the number of integer points $(l,k)$ lying on the line
$l-k=p$
inside the domain of summation in the subtrahend or the minuend in
(\ref{eqGP.1}) doesn't exceed $p$,
 whence
\begin{eqnarray*}
\left\| \sum _{P_1\setminus P_2}\right\| _A &
=&\sum_{P_1\setminus P_2} \left| \hat{F}(k)\right|
\cdot
\left| \hat{\gamma }(l-k)\right| \leq  \\
\sum_{p=-\infty }^\infty |p|\cdot \left| \hat{\gamma }%
(p)\right| \cdot \left\| F\right\| _{PF} &=&\left\| \gamma ^{\prime
}\right\| _A\cdot \left\| F\right\| _{PF}.
\end{eqnarray*}
It suffices now to verify the lemma's assertion for trigonometric
polynomials and then apply the
\INDEX{Ban}{S.Banach}\INDEX{Steinh}{G.Steinhaus}Banach-Steinhaus theorem.
\end{proof}
Theorem \ref{ThMinloc} follows immediately after taking a smooth function $%
\gamma $ which is identically $1$ in some neighborhood of the compact $K$.

The theorem's statement seems to be new. Usually the
\INDEX{Riemann}{B.Riemann}Riemann summability to zero of a general
trigonometric series is required. However, convergence in the sense of
generalized functions is more flexible and better suits for applications in
the theory of differential operators. Note also that the $A(K)-$ convergence
is stronger than the $C(K)-$ one.

\subsection{F.Sch\"afke's approach}

In the beginning of sixties
\INDEX{Schafke}{F.Sch\protect\"afke}F.Sch\"afke developped a general equiconvergence
principle and applied it to convergence of \ef expansions in the complex
domain. However we didn't succeed in translating his general
constructions to the case of an interval of thev real axis
and hence were unable to provide a
comparison with results of other researchers. Therefore
to our regret
we can only refer
the reader to his three thorough articles \cite{schaf60}.

\setcounter{page}{41}
\chapter{Equiconvergence on the whole interval}

\section{Introduction}

\subsection{Notations}
Let $n=2q$ and consider two
Birkhoff-regular
\INDEX{Birkhoff}{G.D.Birkhoff}
$n$-th order
differential operators $L_1$ and $L_2$
defined by two \bvp like
(\ref{Chap.intro}.\ref{eqB.1})-(\ref{Chap.intro}.\ref{eqB.2}).
In what follows set
$
 r_k=2\pi+\alpha,k=1,2,\ldots
$
Then under an appropriate choice of $\alpha>0$
condition
(\ref{Chap.intro}.\ref{eqB.16})
is fulfilled for sets of \chv of both operators $L_1, L_2$.
In the sequel $\|\,\|$ stands for $C(0,1)$-norm and
$\|\,\|_{(a,b)}$ for the norm in $C(a,b)$.

\subsection{Order two case}

At first, note that propositions
\ref{Chap.intro}.\ref{ThA}, \ref{Chap.intro}.\ref{ThB}
don't solve the problem of
equiconvergence on the whole interval $[0,1]$. And it is evident
that in this case some additional restrictions should be imposed
on the expanded function. Of course, it is possible to require more and
more
delicate smoothness conditions but this
process seems to be infinite without any hope for a final solution.
The situation completely changed in 1975
when
\INDEX{Khro}{A.P.Khromov}A.P.Khromov
made a decisive breakthrough \cite{hro75}.
His main idea was to
\textit{reduce}
the equiconvergence problem for general \ef expansions
to that of some model function system, namely, to the question of
\textit{uniform convergence of a trigonometric series}.
Let us state his result.
Its formulation is slightly modified but a proof we found is even
simpler than the original one and is given below in the section
\ref{sec:Two}
\begin{theorem}[Khromov]\label{Ohrowhole}
\INDEX{Khro}{A.P.Khromov}
Given two Birkhoff--regular
\INDEX{Birkhoff}{G.D.Birkhoff}
second order differential
operators
$L_1, L_2$
and a function $f\in C(0,1)$, assume that
\begin{equation}
\label{eqw:0.9}
f\in clos(D_{L_1})\cap clos(D_{L_2}).
\end{equation}
The closure is taken in $C(0,1)$. Then
\begin{equation}
\label{eqw:0.10}
\lim_{k\to\infty}\|S_{r_k}(f,L_1)-S_{r_k}(f,L_2)\|=0
\end{equation}
if and only if
\begin{equation}
\label{eqw:0.11}
\lim_{r\to\infty}\|\sigma_r(\Phi_0)\|_{(-\delta,0)}=
\lim_{r\to\infty}\|\sigma_r(\Phi_1)\|_{(-\delta,0)}
\end{equation}
for some auxiliary functions $\Phi_0,\Phi_1$
 and any fixed
$\delta,\ \ 0<\delta<\frac12$.
\end{theorem}
\begin{remark}
More precisely $\Phi_0,\Phi_1$ are linear combinations of
$f,\,f^\#$, where
\[
f^\#(\xi):=f(1-\xi),\quad 0\le\xi\le1,
\]
$\Phi_0,\Phi_1\equiv 0$ off $(0,1)$ and they are explicitly defined in
(\ref{eqw:0.17}).
\end{remark}

Thus the problem which is set above can be reduced to the
classical question of trigonometric series convergence.
We recall that there are a lot of
strong convergence criteria for these serii which go back to
\INDEX{Young}{W.H.Young}Young,
\INDEX{Lebesgue}{H.Lebesgue}Lebesgue,
\INDEX{Vall}{Ch.J. de la Vall\'ee-Poussin}de la Vall\'ee-Poussin
and others (see \cite{Zyg:ts}).
The
latest  belongs to
\INDEX{Wer}{E.Wermuth}E.Wermuth \cite{wer89}
and generalizes
all the previous
ones.\hfill

\subsection{$X$-equivalence}
 In his thesis
\INDEX{Wer}{E.Wermuth}E.Wermuth
\cite[p.61-73]{wer:phd} raised an
equiconvergence problem on the whole interval for two \ef
expansions
associated with
\INDEX{Birkhoff}{G.D.Birkhoff}Birkhoff-regular
operators $L_1, L_2$
{\it simultaneously for all functions}
in some class $X$.
More precisely, he introduced the following
\begin{definition}
Given two $n$-th order differential operators
$L_1, L_2 \in(R)$ in $L^2(0,1)$
and a function class
$X\subset L^1(0,1)$
we say that these operators
are $X$-equivalent if \hfill
\[
  \lim\|S_{r_k}(f,L_1)-S_{r_k}(f,L_2)\|=0\ \  \forall f\in X.
\]
\end{definition}
We shall also say that operators $L_1,L_2$
\textbf{essentially coincide} if the orders $r_j$ and
the leading parts $V_J$ of their boundary forms are identical.
Then
\INDEX{Wer}{E.Wermuth}E.Wermuth's
theorem \cite[Satz 13]{wer:phd} states that two
\INDEX{Birkhoff}{G.D.Birkhoff}Birkhoff-regular operators $L_1, L_2$
essentially coincide if and only if they are $L^1(0,1)$-equivalent.\hfill

This statement is also valid for $X=L^p(0,1)$,
where $p$ is a fixed number, $1\le p <\infty$
(see \cite[pp.63,72]{wer:phd}).
The less can be chosen the set $X$, the weaker is
the theorem's claim
and the stronger is the result.\hfill

The main difficulty in
\INDEX{Wer}{E.Wermuth}E.Wermuth's
problem consists of finding necessary
and sufficient conditions for uniform boundedness of the
family of operators:\hfill
\[
f\to(S_{r_k}(f,L_1)-S_{r_k}(f,L_2)),\quad k=1,2,3\ldots,
\]
acting from
$X$ to $C[0,1]$.
He solved it for $X=L^p(0,1),\ \ 1\le p<\infty$.
The extreme case $p=\infty$ and all the more $X=C(0,1)$ remained
open \cite[p.73]{wer:phd}.
 Of course, this result gives no answer to equiconvergence on
the whole interval for any {\it given fixed function} $f$
as would be desired if we intend to generalize the
\INDEX{Tam}{J.Tamarkin}Tamarkin-\INDEX{Sto}{M.Stone}Stone's
theorem
\ref{Chap.intro}.\ref{ThA}.
Obviously, this is a much more
subtle question than the analogous one for a class of functions.

\subsection{Higher order case}
Hence, in sharp contrast with the case of equiconvergence
in the internal points there is \textit{absolutely no results}
 concerning generalization
of the
\INDEX{Tam}{J.Tamarkin}Tamarkin-\INDEX{Sto}{M.Stone}Stone's
theorem to the whole interval when $n>2$.
Such state of affairs stems from the dificulty of the problem in question.
The standard resolvent's approach is good enough to obtain
{\em sufficient} conditions
provided $f$ is enough smooth but
fails to give necessary and sufficient ones, see, for instance
\cite{kau89,kau:phd,gora91} and others (the list may be
considerably increased).

At first, we give such criterion
provided\hfill
\begin{equation}
\label{eqw:0.12}
f=f_0\in C_0(0,1):=\{g\in C(0,1):\ g(0)=g(1)=0\}
\end{equation}
and $n$ is even, see theorems
\ref{Te:whole}, \ref{Ttrige1:whole}.
A more general situation with $f\in C(0,1)$
is reduced to this one in the section \ref{sec:zero}.
Odd order case is solved in the section \ref{sec:odd}.

To begin with let us briefly sketch the idea of the proof.
We have a quantity
(namely, a difference of two partial sums)
which tends to zero in $C(0,1)$. Expand then the numerator of
the
\INDEX{Green}{G.Green}Green function along the uppermost row
and each of the occurring
minors with respect to the leftmost column, containing the column-vector
$W$ (if any). In the appearing finite sum
some summands tend to zero uniformly on
$[0,1]$.
After eliminating them we obtain an equation with a sum of some
leading terms from the left and $o(1)$ from the right.
 In order to extract these terms we need, say, a system
of linear equations with a right-hand side $o(1)$.
But we have only one equation with many unknowns!

Here a new idea is invoked:
\textit{we differentiate the left-hand side $j$ times}
and divide the result by $r^j,\ \ j=0,\dots,n$.
Then each of the leading terms is factored by $\varepsilon_k^j$
(see below)
in the $j$th equation up to additional summands of the type $o(1)$.
When $j=n$ we return to the original equation.
Of course, the right-hand side is changed but it remains $o(1)$
and its concrete value is of no importance for us.
Next we employ another important idea, we invoke
\textit{inequality for derivatives}
 and thus arrive at the desired system
(after joining some terms pairwise but we omit details).

Further we shall need some notations.
Let
$n=2q>2$.
For any $n$-th order differential operator $L\in(R)$ let
$A_j,B_j$ be $n$-columns,
\begin{equation}
\label{eqw:0.13}
A_j:=[a_\nu \varepsilon_j^{\sigma_\nu}]_{\nu=0}^{n-1},\ \
B_j:=[b_\nu \varepsilon_j^{\sigma_\nu}]_{\nu=0}^{n-1},\ \
j=0,\ldots,n-1,
\end{equation}
the regularity determinant $\Theta=\Theta(b^0,b^1)\ne0$ reads as follows
\begin{equation}
\label{eqw:0.14}
\Theta:=\det[A_0...A_{q-1}B_{q}...B_{n-1}],
\end{equation}
$\Theta^{\nu k}$ be a cofactor of the $(\nu,k)$
entry of $\Theta$, indices $\nu,k$ vary from $0$ to $n-1$.
The same notation will be used for other matrices.
Let \hfill
\[
d_{m\nu}=
\left\{
 \begin{array}{rl}
  -a_\nu,           \quad           q\le m\le n-1\\
  \phantom{-}b_\nu, \quad           0\le m\le q-1
 \end{array}
\right.
\]
and set\hfill
\begin{equation}
\label{eqw:0.16}
\alpha_{mk}=\alpha_{mk}(L)=\frac{1}{2\pi}\sum_{\nu=0}^{n-1}
\varepsilon_m^{-(n-1-\sigma_\nu)}d_{m\nu}
\cdot \Theta^{\nu_k}/\Theta.
\end{equation}
 For $f\in C(0,1)$ introduce a collection of functions associated
with the given function $f$:\hfill
\begin{align}
\label{eqw:0.17}
\Phi_k(\xi,f,L)&:=\alpha_{qk}f(\xi)+\alpha_{0k}f^\#(\xi),      &
k&=0,\ldots,n-1;\\
\label{eqw:0.18}
\Psi_k(\xi,f,L)&:=\alpha_{n-1,k}f(\xi)+\alpha_{q-1,k}f^\#(\xi), &
k&=0,\ldots,n-1.
\end{align}
For the sake of brevity we shall write further
$r,\ \Gamma_r$ instead of $r_k$
and $\Gamma_{r_k}$, respectively.
\begin{theorem}
 \label{Te:whole}
Given a function $f_0$ satisfying (\ref{eqw:0.12}) and assume that\hfill
\begin{equation}
\label{eqw:0.19}
span(\varphi_0,\varphi_q)=span(\psi_q,\psi_{n-1})=span(f_0,f_0^\#)
\end{equation}
where\hfill
\begin{equation}
  \label{eqw:0.20}
 \left\{
\begin{split}
\varphi_k&:=\Phi_k(\cdot,f_0,L_1)-\Phi_k(\cdot,f_0,L_2),\\
\psi_k&:=\Psi_k(\cdot,f_0,L_1)-\Psi_k(\cdot,f_0,L_2),
\end{split}
 \right.
\ \ \ \ \
k=0,\ldots,n-1.
\end{equation}
Let\hfill
\begin{equation}
\label{eqw:0.21}
I_r^\pm(f_0)(x):=\int\limits_{1/r}^1
(x+\xi)^{-1}\exp(\pm ir\xi)f_0(\xi)\,d\xi.
\end{equation}
Then\hfill
\begin{equation}
\label{eqw:0.22}
\lim\limits_{r\to\infty}\|S_r(f_0,L_1)-S_r(f_0,L_2)\|=0
\end{equation}
if and only if\hfill
\begin{equation}
\label{eqw:0.23}
\lim\limits_{r\to\infty}\|I_r^\pm(g)\|=0,\ \ g=f_0, f_0^\#.
\end{equation}
\end{theorem}
Let us indicate that
in theorem \ref{Te:whole}
 only {\it a unique function\/} $f_0$ is considered,
i.e. we merely put a set
$X=\{ f_0 \}$ consisting of {\it one\/} element. Hence in our case
the operators
$L_1, L_2$ may be absolutely different and really we have established a
\textit{true generalization} of the Tamarkin-Stone's
\INDEX{Tam}{J.Tamarkin}Tamarkin-\INDEX{Sto}{M.Stone}
theorem.\hfill

\begin{remark}
It is possible to change
$\|\,\|$ in (\ref{eqw:0.23}) by
$\| \, \|_ {(0,\delta)} $
for any fixed $\delta,0<\delta<1$.
There is also a simple
sufficient condition for (\ref{eqw:0.19}) to be valid:\hfill
\begin{equation}
\label{eqw:0.24}
\begin{array}{l}
\det\left[
  \begin{array}{cc}
    \beta_{00}& \beta_{q0}\\
    \beta_{0q}& \beta_{qq}
  \end{array}
    \right]\ne0,
\ \ \ \
\det\left[
  \begin{array}{cc}
    \beta_{q-1,q-1}& \beta_{n-1,q}\\
    \beta_{q-1,n-1}& \beta_{n-1,n-1}
  \end{array}
    \right]\ne0,
\\
\end{array}
\end{equation}
where \hfill
\[
\beta_{mk}:=\alpha_{mk}(L_1)-\alpha_{mk}(L_2).
\]
 Of course, it is
equivalent to (\ref{eqw:0.19}) if the functions
$f_0$ and $f_0^\#$ are linearly
independent.\hfill
\end{remark}

\section{Green's function}\label{sec:wGreen}
\INDEX{Green}{G.Green}

\subsection{New fundamental system of solutions}
Consider a
\INDEX{Birkhoff}{G.D.Birkhoff}Birkhoff-regular
\bvp
\begin{equation}
\label{eqw:1.1}
D^ny=\lambda y+f,
\end{equation}
\begin{equation}
\label{eqw:1.2}
V_\nu(y)=0,\quad \nu=0,\ldots,n-1.
\end{equation}
Equation $D^n y=\lambda y$ has an evident
\fss, namely\hfill
\begin{equation}
\label{eqw:1.3}
y_k(x,\varrho)\equiv\exp(i\varrho\varepsilon_k x),\ \ k=0,\ldots,n-1.
\end{equation}
Introducing kernels\hfill
\begin{equation}
\label{eqw:1.4}
g(x,\xi,\varrho)=i\cdot
\left\{
  \begin{array}{rl}
  \sum\limits_{k=0}^{q-1}\varepsilon_k^{-(n-1)}y_k(x-\xi),\ \ x>\xi\\
  -\sum\limits_{k=q}^{n-1}\varepsilon_k^{-(n-1)}y_k(x-\xi),\ \ x>\xi\\
  \end{array}
\right.
\end{equation}
and\hfill
\begin{equation}
\label{eqw:1.5}
g_0(x,\xi,\varrho):=g(x,\xi,\varrho)/(n\varrho^{n-1})
\end{equation}
we get a particular solution $g_0(f)$ of (\ref{eqw:1.1}),
\[
g_0(f):=\int\limits_0^1g_0(x,\xi,\varrho)f(\xi))d\xi.
\]
In the sequel it will be more convenient to use another \fss
$\{z_k\}_{k=0}^{n-1}$, where\hfill
\begin{equation}
\label{eqw:1.6}
z_k(x,\varrho):=
\left\{
  \begin{array}{ll}
   y_k(x,\varrho),  &   \quad k=0,\ldots,q-1,\\
   y_k(x-1,\varrho),&   \quad k=q,\ldots,n-1.
  \end{array}
\right.
\end{equation}
This choice of a \fss is natural due to the fact that\hfill
\begin{equation}
\label{eqw:1.7}
z_k=\mbox{O}(1)\;,k=0,\ldots,n-1;\;\;\; g(x,\xi,\varrho)=\mbox{O}(1),
\;\;0\le x,\xi\le1
\end{equation}
for $\varrho\in S_0:=\{0\le\arg\varrho\le2\pi/n\}$.\hfill

\subsection{Green's function representation}
\INDEX{Green}{G.Green}
By variation of constants we get a well-known expression
(\ref{Chap.intro}.\ref{eqS.1})-(\ref{Chap.intro}.\ref{eqS.2})
for
the \INDEX{Green}{G.Green}Green's function
as a ratio of two determinants with the
\fss being chosen as in (\ref{eqw:1.6}).
Further on, canceling powers of $\varrho$
in the nominator and denominator, we get problem's
(\ref{eqw:1.1})--(\ref{eqw:1.2})
solution of the form :\hfill
\[
y=G(f):=\int\limits_0^1G(x,\xi,\varrho)f(\xi)\,d\xi,
\]
where\hfill
\begin{equation}
\label{eqw:1.8}
G(x,\xi,\varrho)=i\cdot\det H/(n\varrho^{n-1}\det\eta).
\end{equation}
Here\hfill
\begin{equation}
\label{eqw:1.9}
\eta=[\eta_{\nu k}]_0^{n-1},\quad
\eta_{\nu k} := \varepsilon_k^{\sigma_\nu}
( b_\nu z_k(1) + a_\nu z_k(0) ),
\end{equation}
\begin{equation}
\label{eqw:1.10}
H(x,\xi,\varrho):=
\left[
 \begin{array}{cc}
  g & z^T\\
  W & \eta
 \end{array}
\right],
\end{equation}
square brackets here denote a matrix and
$z^T$ stands for the transposed of the column-vector
$z=(z_k(x,\varrho))_0^{n-1}$,\hfill
\begin{equation}
\label{eqw:1.11}
W=(W_\nu)_{\nu=0}^{n-1},\quad
W_\nu=\sum\limits_{m=0}^{n-1}
\varepsilon_m^{-(n-1-\sigma_\nu)}d_{m\nu}u_m(\xi,\varrho)
\end{equation}
and at last\hfill
\begin{equation}
\label{eqw:1.12}
u_m(\xi,\varrho):=
\left\{
 \begin{array}{ll}
  y_m(1-\xi,\varrho),& m=0,\ldots,q-1,\quad 0\le\xi\le1\\
  y_m(-\xi,\varrho), & m=q,\ldots,n-1.
 \end{array}
\right.
\end{equation}
Changing a little bit notation from \cite[p.1185]{efs90b} we use an
abbreviation:\hfill
\[
[[q]]:=
q+O(e^{-\delta|Im\varrho|})
 + O(e^{-\delta|\varepsilon_{q-1}Im\varrho|})
 +O(\frac1\varrho),\ \
\varrho\in\Gamma_r.
\]
The quantity $q$ may vary with $\varrho$.\hfill
\begin{lemma}\label{L1.1}
The following relation is valid\hfill
\begin{equation}
\label{eqw:1.13}
(\det\eta)^{-1}=[[\Theta^{-1}]].
\end{equation}
\end{lemma}
\begin{proof}
We have that\hfill
\begin{equation}
\label{eqw:1.14}
z_k(0)=[[0]],\ \ k=q,\ldots,n-1;\;\;
z_k(1)=[[0]],\ \ k=0,\ldots,q-1.
\end{equation}
Let us expand $\det\eta$ into a sum of determinants with only $z_k(0)$ or
$z_k(1)$ in each column. Then all of them become $[[0]]$ except one
which coincides with $\Theta$. Moreover, from \cite[p.77--78]{Nai:ldo} it
follows that\hfill
\begin{equation}
\label{eqw:1.15}
|\det\eta|\ge C, \quad
\varrho\in S_\delta:=
S\setminus\{|\varrho-\varrho_j^i|\le\delta\}_{j=1}^\infty,
\quad i=1,2
\end{equation}
for sufficiently small $\delta>0$.
Here $\varrho_j^i$ stands for the \chv of the operator $L_i,\ \ i=1,2$.
Therefore,\hfill
\begin{equation}
\label{eqw:1.16}
(\det\eta)^{-1}-\Theta^{-1}=(\Theta-\det\eta)/(\Theta\det\eta)=[[0]].
\end{equation}
\end{proof}
A similar calculation can be found in \cite[p.1186]{efs90b}.\hfill

\subsection{A partial sum's formula}
Let us consider a partial sum\hfill
\begin{equation}
\label{eqw:1.17}
S_r(f):=(-2\pi i)^{-1}\int_{\Gamma_r} G(f)(x,\varrho)n\varrho^{n-1}\,d\varrho
\end{equation}
of \ef expansion for the \bvp
(\ref{eqw:1.1})--(\ref{eqw:1.2}).
Taking into account (\ref{eqw:1.8}),(\ref{eqw:1.14}) and
expanding $\det H$ with respect to the first row and each of the
appearing minors except the first one
with respect to the column $W$ we arrive at
 identity\hfill
\begin{equation}
\label{eqw:1.18}
S_{r}(f)
\equiv
S_{r,0}(f)+
\sum\limits_{k,\nu,m=0}^{n-1}J_{r,m,\nu,k}(f),
\end{equation}
where\hfill
\begin{equation}
\label{eqw:1.19}
S_{r,0}(f):=(-2\pi i)^{-1}\int_{\Gamma_r} g(f)\,d\varrho,
\end{equation}
\[
g(f):=\int\limits_0^1g(x,\xi,\varrho)f(\xi)\,d\xi
\]
and\hfill
\begin{eqnarray}
\label{eqw:1.20}
J_{r,m,\nu,k}(f):=
&(-2\pi)^{-1}\int_{\Gamma_r}\{z_k(x,\varrho)
(-1)^{k+3+k}
([[\Theta^{\nu k}]]/\Theta)\hfill
&\\
&\cdot\left\{\int\limits_0^1f(\xi)u_m(\xi,\varrho)\,d\xi
\varepsilon_m^{-(n-1-\sigma_\nu)}
\,d_{m\nu}\right\}
d\varrho.& \nonumber
\end{eqnarray}
Here we also  took into account the identity
$\eta^{\nu k}\equiv[[\Theta^{\nu k}]]$.
It can be easily deduced expanding $\det H$ along the uppermost row
and each of the occurring minors along the leftmost column.

 It will be helpful to rewrite (\ref{eqw:1.20}) in the following way:\hfill
\begin{equation}
\label{eqw:1.21}
J_{r,m,\nu,k}(f)=\alpha_{m\nu k}\cdot\tau_{rmk}(f[[1]])
\end{equation}
where\hfill
\begin{equation}
\label{eqw:1.22}
\alpha_{m\nu k}:=(2\pi)^{-1}\varepsilon_m^{-(n-1-\sigma_\nu)}d_{m\nu}
(\Theta^{\nu k}/\Theta)
\end{equation}
and\hfill
\begin{equation}
\label{eqw:1.23}
\tau_{rmk}(f):=
\int_{\Gamma_r} z_k(x,\varrho)\,
\left(\int\limits_0^1 f(\xi)u_m(\xi,\varrho)\,\right) \,d\varrho.
\end{equation}
\begin{remark}
 By an elementary computation
(see also analogous result in
\textrm{\cite[p.1191]{efs90b}})
we get that\hfill
\begin{equation}
\label{eqw:1.24}
S_{r,0}(f)\equiv\sigma_r(f).
\end{equation}
\end{remark}

\section{Equi\-con\-ver\-gence with a trigono\-met\-ric
 Fou\-rier integral}
\INDEX{Fou}{J.Fourier}

\subsection{Simplifications} \label{subsec:simpl}
\begin{lemma}
\label{L2.1}
Let us represent the second factor in (\ref{eqw:1.21}) as a sum\hfill
\begin{equation}
\label{eqw:2.1}
\tau_{rmk}(f)+\tau_{rmk}(f[[0]]).
\end{equation}
Then the following estimate is valid\hfill
\begin{equation}
\label{eqw:2.2}
\|\tau_{rmk}(f[[0]])\|=o(1)\; \mbox{as }\; r\to\infty,\;f\in L(0,1).
\end{equation}
\end{lemma}
\begin{proof}
We shall use a well-known relation:\hfill
\begin{equation}
\label{eqw:2.3}
\sup\limits_{r>0}\int\limits_{\Gamma_r}\exp(-\varepsilon|Im\varrho|)|
\,d\varrho|<\infty.
\end{equation}
Thus we have a uniformly bounded family of linear operators acting
from $L(0,1)$ to $C(0,1)$:\hfill
\begin{equation}
\label{eqw:2.4}
f\to\tau_{rmk}(f[[0]]).
\end{equation}
Obviously it tends to zero on a dense lineal $C_0^\infty(0,1)$
since the factor $[[0]]$ depends only on $\varrho$, not on $\xi$.
Therefore it
remains to apply the
\INDEX{Ban}{S.Banach}\INDEX{Steinh}{G.Steinhaus}Banach-Steinhaus theorem.
\end{proof}
\begin{lemma}
\label{L2.2}
For any
$
m\notin\{0,q-1,q,n-1\}\;\;\;
\|\tau_{rmk}(f)\|=\mbox{o}(1)$ as $r\to\infty.
$
\end{lemma}
\begin{proof}
Assume for defficiency that $0<m<q-1$. Then\hfill
\[
\begin{split}
|u_m(\xi,\varrho)|&= |y_m(1-\xi,\varrho)|\\
                  &=\exp
\bigl(-Im(\varrho\varepsilon_m)(1-\xi) \bigr)\\
                  & \le
\exp\bigl(-\varepsilon|\varrho|(1-\xi)   \bigr)
\end{split}
\]
with some positive $\varepsilon$. Hence\hfill
\begin{eqnarray}
\label{eqw:2.5}
\int_{\Gamma_r}\int\limits_0^1|u_m(\xi,\varrho)|\,d\xi\,|d\varrho|
\le
\int_{\Gamma_r}|\,d\varrho
|\int\limits_0^1exp\{-\varepsilon|\varrho|(1-\xi)\}\,d\xi
\\
=(2\pi/n)\int\limits_0^1r\exp\{-\varepsilon r(1-\xi)\}\,d\xi
\le
\varepsilon^{-1}2\pi/n. \nonumber
\end{eqnarray}
Since $f\in C_0(0,1)$, it suffices to recall (\ref{eqw:1.7}) and ap\-ply
the
\INDEX{Ban}{S.Banach}\INDEX{Steinh}{G.Steinhaus}Ba\-nach--Stein\-haus
theorem to the family
$r\to\tau_{rmk}(f)$ of
operators
acting from $C_0(0,1)$ to $C(0,1)$.
\end{proof}

\subsection{Remainder formula} \label{subsec:remain}
Now we consider the sum\hfill
\begin{equation}
\label{eqw:2.6}
\sum\limits_{k,\nu,m=0}^{n-1}J_{r,m,\nu,k}(f)
\end{equation}
and join pairwise all its summands with $m=0$ and $m=q$ or $m=q-1$
and $m=n-1$, respectively.
Lemmas \ref{L2.1}, \ref{L2.2} together with
(\ref{eqw:1.19}), (\ref{eqw:1.22}) yield an important identity\hfill
\begin{align}
\label{eqw:2.7}
S_r(f)-\sigma_r(f) &\equiv error\\
+\sum\limits_{k=0}^{n-1}
\bigl[&\ \{\alpha_{0k}\tau_{r0k(f)}+
\alpha_{qk}\tau_{rqk}(f)\}\nonumber\\
&+
\{\alpha_{q-1,k}\tau_{r,q-1,k}(f)+
\alpha_{n-1,k}\tau_{r,n-1,k}(f)\}
\bigr], \nonumber
\end{align}
where {\it error} stands for the summands tending to zero in
$C(0,1)$ as
$r\to\infty$. \hfill
\begin{lemma}
\label{L2.4}
The following identities are valid\hfill
\begin{align}
\label{eqw:2.8}
\alpha_{0k}\int\limits_0^1f(\xi)u_0(\xi,\varrho)\,d\xi
&+
\alpha_{qk}\int\limits_0^1 f(\xi)u_q(\xi,\varrho)\,d\xi\nonumber\\
&\equiv
\int\limits_0^1\Phi_k(\xi)y_0(\xi,\varrho)d\xi,\\
\label{eqw:2.9}
\alpha_{q-1,k}\int\limits_0^1f(\xi)u_{q-1}(\xi,\varrho)\,d\xi
&+
\alpha_{n-1,k}\int\limits_0^1f(\xi)u_{n-1}(\xi,\varrho)\,d\xi\nonumber\\
&\equiv
\int\limits_0^1\Psi_k(\xi)y_{q-1}(\xi,\varrho)d\xi,
\end{align}
\end{lemma}
\begin{proof}
It suffices to notice that
$y_{j+q}(\xi,\varrho)\equiv y_j(-\xi,\varrho)$
and to make a sub\-sti\-tu\-tion:
$\xi\to 1-\xi$ if needed. \end{proof}
\begin{corollary}\label{wcor1}
Formulas (\ref{eqw:2.7})--(\ref{eqw:2.9}) yield a final relation:\hfill
\begin{align}
\label{eqw:2.10}
S_r(f)-\sigma_r(f)
&\equiv error\\
&+\sum\limits_{k=0}^{n-1}
\bigl\{
\int_{\Gamma_r}z_k(x,\varrho)
\int\limits_0^1
\Phi_k(\xi)y_0(\xi,\varrho)d\xi\,d\varrho\nonumber\\
&+\int_{\Gamma_r}z_k(x,\varrho)
\int\limits_0^1\Psi_k(\xi)y_{q-1}(\xi,\varrho)
d\xi\,d\varrho
\bigr\} \nonumber\\
&=error+\sum\limits_{k=0}^{n-1}\{\eta_k(x,r,f)+\zeta_k(x,r,f)\},\nonumber
\end{align}
where $\|error\|\to 0$ as $r\to\infty$.\hfill
\end{corollary}

\subsection{Preliminary transformations}
Set\hfill
\begin{equation}
\label{eqw:2.12}
\gamma_0=\gamma_0(x,r,f):=\sum\limits_{k=0}^{q-1}
[\eta_k(x,r,f)+\zeta_k(x,r,f)],
\end{equation}
\begin{equation}
\label{eqw:2.13}
\gamma_1=\gamma_1(x,r,f):=\sum\limits_{k=q}^{n-1}
[\eta_k(x,r,f)+\zeta_k(x,r,f)].
\end{equation}
Then\hfill
\[
 S_r(f)-\sigma_r(f)=\gamma_0+\gamma_1+error.
\]
\begin{lemma}
An estimate\hfill
\[
\|S_r(f)-\sigma_r(f)\|=o(1),\ \ r\to\infty
\]
is valid if and only if for any fixed
$\delta,\quad 0<\delta<1/2$ \hfill
\begin{equation}
\label{eqw:2.14}
\|\gamma_0\|_{(0,\delta)}=o(1),\ \
\|\gamma_1\|_{(1-\delta,1)}=o(1)\;\; \mbox{as } r\to\infty.
\end{equation}
\end{lemma}
\begin{proof}
 It suffices to show that\hfill
\begin{align}
\label{eqw:2.15}
 \|\eta_k\|_{(\delta,1)}&=o(1),       &
 \|\zeta_k\|_{(\delta,1)}&=o(1),      &
                                      & k=0,\ldots,q-1,\\
\label{eqw:2.16}
 \|\eta_j\|_{(0,1-\delta)}&=o(1),     &
 \|\zeta_j\|_{(0,1-\delta)}&=o(1),    &
                                      & j=q,\ldots,n-1
\end{align}
as $r\to\infty$. To be definite we shall consider only (\ref{eqw:2.15}).
But\hfill
\[
\|z_k(\cdot,\varrho)\|_{(\delta,1)}
=O(\exp(-Im(\varrho\varepsilon_k)\delta))
=O(\exp(-\delta_1Im\varrho)
\]
for
$\varrho\in\Gamma_r,\ \ k=0,\ldots,q-2$ and some positive $\delta_1$.
Quite analogously,\hfill
\[
\|z_k(\cdot,\varrho)\|_{(\delta,1)}=O(exp(-\delta_1Im\varrho))
\]
for\hfill
\[
\varrho\in\varepsilon_{q-1}\Gamma_r
:=\{\varrho\varepsilon_{q-1}:\varrho\in\Gamma_r\}.
\]
The circular arc $\varepsilon_{q-1}\Gamma_r$
lies in the upper half-plane.
Therefore, we can use (\ref{eqw:2.3}) in any case.
At last the same hint with the
\INDEX{Ban}{S.Banach}\INDEX{Steinh}{G.Steinhaus}Banach-Steinhaus
theorem completes
the proof.
\end{proof}

\subsection{Behaviour of the main terms under differentiation}

Let us consider a difference\hfill
\begin{equation}
\label{eqw:2.17}
D^n\eta_k-r^n\eta_k
=\int_{\Gamma_r}(\varrho^n-r^n)z_k(x,\varrho)
\int\limits_0^1\Phi_k(\xi)y_0(\xi)d\xi\,d\varrho.
\end{equation}
\begin{lemma}\label{L2.6}
The following relation holds uniformly with respect to
$x$,\hfill\\ $0\le x\le1$:\hfill
\begin{equation}
\label{eqw:2.18}
D^j\eta_k-(r\varepsilon_k)^j\eta_k=o(r^j),\ \
j=1,\ldots,n;
\quad
k=0,\ldots,n-1.
\end{equation}
\end{lemma}
\begin{proof}
Since \hfill
\[
D^jz_k(x,\varrho)\equiv(\varrho\varepsilon_k)^jz_k(x,\varrho)
\]
we see
that the left-hand side in (\ref{eqw:2.18}) differs from
$\eta_k$
by the additional factor
$\varepsilon_k^j(\varrho^j-r^j)$ in (\ref{eqw:2.10}).
Then, taking  into account
(\ref{eqw:1.7}) we get that\hfill
\begin{equation}
\label{eqw:2.19}
\|D^j\eta_k-(r\varepsilon_k)^j\eta_k\|
\le
C\int_{\Gamma_r} |\varrho^j - r^j|
\int\limits_0^1 |y_0(\xi)|\,d\xi
\, |d\varrho| \cdot\|f\|.
\end{equation}
In the meantime,\hfill
\[
|y_0(\xi,\varrho)|\equiv\exp(-Im\varrho\,\xi)\le\exp(-2r\varphi\xi/\pi),
\]
$\varrho=r\exp(i\varphi),\quad 0\le\varphi<2\pi/n$.
Then the right-hand side in (\ref{eqw:2.19})
is less or equal to\hfill
\begin{equation}
\label{eqw:2.20}
 \begin{array}{l}
  C\|f\|r^{j+1}
\int\limits_0^{2\pi/n}\,d\varphi
\int\limits_0^1\,d\xi\exp(-2r\varphi\xi/\pi)|\exp(ij\varphi)-1| \\
  \le C\|f\|r^{j+1}\int\limits_0^{2\pi/n}\,d\varphi
\int\limits_0^1\,d\xi\,\exp(-2r\varphi\xi/\pi)O(j\varphi).
\end{array}
\end{equation}
Replacing the internal integral in (\ref{eqw:2.20}) by\hfill
\[
\int\limits_0^\infty j\varphi\exp(-2r\varphi\xi/\pi)=\pi j/2r
\]
we come to (\ref{eqw:2.18}) with $O(r^j)$ instead of $o(r^j)$.
It remains now to
apply the
\INDEX{Ban}{S.Banach}\INDEX{Steinh}{G.Steinhaus}Banach-Steinhaus theorem.
\end{proof}
\begin{remark}\label{Rorj}
Proceeding in the same way we find that\hfill
\[
D^j\zeta_k-(r\varepsilon_1\varepsilon_k)^j\zeta_k
=o(r^j),\ \
j=1,\ldots,n;\ \
k=0,\ldots,n-1
\]
as $r\to\infty$. In this case we have only to take into account
that the modulus\hfill
\[
|y_{q-1}(\xi,\varrho)|,\quad \varrho\in\Gamma_r
\]
attains its maximum at the
point $r\varepsilon_1$, not at $r$ as $|y_0(\xi,\varrho)|$.\hfill
\end{remark}
\begin{lemma} \label{L2.7}
Let $\delta$ be any fixed number, $0<\delta<1/2$. Then the
following inequalities hold for $r\ge1$ and $j=1,\ldots,n$\hfill
\begin{equation}
\label{eqw:2.22}
\|D^j\gamma_0\|_{(0,\delta)}\le Cr^j\|\gamma_0\|_{(0,\delta)}+o(r^j),
\end{equation}
\begin{equation}
\label{eqw:2.23}
\|D^j\gamma_1\|_{(1-\delta,1)}\le Cr^j\|\gamma_1\|_{(1-\delta,1)}+o(r^j).
\end{equation}
\end{lemma}
\begin{proof}
For the sake of being definite consider only (\ref{eqw:2.22}).
According to lemma \ref{L2.6}
 and remark \ref{Rorj}
 we have that\hfill
\begin{equation}
\label{eqw:2.24}
D^j\gamma_0=\sum\limits_{k=0}^{q-1}\{
(r\varepsilon_k)^j\eta_k+
(r\varepsilon_1\varepsilon_k)^j\zeta_k\}
+o(r^j),\ \ j\ge1
\end{equation}
uniformly with respect to $x,\ \ 0\le x\le\delta$.
Thus for $j=n$ we have that\hfill
\begin{equation}
\label{eqw:2.25}
D^n\gamma_0=r^n\gamma_0+o(r^n),\quad 0\le x\le\delta.
\end{equation}
Now it is time to recall an inequality for derivatives
\cite[p.131]{BeIlNi:int}.
Let $g$ be any continuously differentiable function in $(0,\delta)$.
Then\hfill
\[
\|g^{(j)}\|_{(0,\delta)}
\le C_1(\|g\|_{(0,\delta)})^{(n-j)/n}
(\|g^{(n)}\|_{(0,\delta)})^{j/n}+C_2\|g\|_{(0,\delta)}.
\]
Applying it to $g=\gamma_0|_{(0,\delta)}$
and taking into account (\ref{eqw:2.25})
we come to an inequality\hfill
\[
\|D^j\gamma_0\|_{(0,\delta)}
\le
C_1(\|\gamma_0\|_{(0,\delta)})^{(n-j)/n}
(\|\gamma_0^n\|_{(0,\delta)})^j/n
+
C_2\|\gamma_0\|_{(0,\delta)}.
\]
This yields (\ref{eqw:2.22}) after considering two cases:\hfill
\[
i)\|\gamma_0\|_{(0,\delta)}
\le\varepsilon,\;\;\;
ii)\|\gamma_0\|_{(0,\delta)}
\ge
\varepsilon,
\]
$\varepsilon$ being sufficiently small.
\end{proof}

\subsection{Main criterion}
According to lemmas \ref{L2.6} and \ref{L2.7} we have that\hfill
\[
\|S_r(f)-\sigma_r(f)\|=o(1)
\]
if and only if\hfill
\begin{align}
\label{eqw:2.26}
i)\ \   &r^{-j}\|D^j\gamma_0\|_{(0,\delta)}=o(1),  \ \ &
  j=0,\ldots,q;\\[5pt]
\label{eqw:2.27}
ii)\   &r^{-j}\|D^j\gamma_1\|_{(1-\delta,1)}=o(1), \ \ &
  j=0,\ldots,q.
\end{align}
Of course, we can replace the number $q$ in
(\ref{eqw:2.26})-(\ref{eqw:2.27})
by any other nonnegative one but our choice is suitable for
further purposes.\hfill

Substituting (\ref{eqw:2.24}) into (\ref{eqw:2.26}) we obtain a system of
equations with respect to the variables $\eta_k,\zeta_k$\hfill
\[
\left\{
\sum\limits_{k=0}^{q-1}
(\varepsilon_k^j\eta_{k}+\varepsilon_{k}^j\zeta_{k})
=o(1),\ \
0\le x\le\delta,\ \
j=0,\ldots,q
\right.
\]
or in a modified form\hfill
\begin{equation}
\label{eqw:2.28}
\left\{
\begin{array}{rl}
\varepsilon_0^j\eta_0 + \varepsilon_1^j(\eta_1+\zeta_0) + \ldots\\[5pt]
+ \varepsilon_{q-1}^j(\eta_{q-1} + \zeta_{q-2})+
+ \varepsilon_{q}^j\zeta_{q-1}
=o(1), \\[5pt]
j=0,\ldots,q,\quad 0\le x\le\delta.&
\end{array}
\right.
\end{equation}
Since the system's determinant coincides with a
\INDEX{Van}{A.T.Vandermonde}Vandermonde one\hfill
\[
|\varepsilon_k^j|_{j,k=0}^q,
\]
we solve it immediately:\hfill
\begin{equation}
\label{eqw:2.29}
\eta_0=o(1),\ \
\eta_1+\zeta_0=o(1),\ \
\ldots,\eta_{q-1}+\zeta_{q-2}
=o(1)
\end{equation}
uniformly with respect to $x,\ \ 0\le x\le\delta$.
Proceeding in the same way, we derive from (\ref{eqw:2.27}) that\hfill
\begin{equation}
\label{eqw:2.30}
  \eta_q=o(1),\
  \eta_{q+1} + \zeta_{q}=o(1),\ldots,\eta_{n-1} + \zeta_{n-2}=o(1),\
  \zeta_{n-1}=o(1)
\end{equation}
uniformly with respect to $x,\ \ 1-\delta\le x\le1$.\hfill
\begin{theorem}
\label{Ttrige:whole}
Let $L$ be a
\INDEX{Birkhoff}{G.D.Birkhoff}Birkhoff-reqular differential operator
of the form
(\ref{Chap.intro}.\ref{eqB.1})-(\ref{Chap.intro}.\ref{eqB.2}),
$\delta$ be any positive number $\le1$ and $f_0\in C_0(0,1)$.
Then\hfill
\[
\lim\limits_{r\to\infty}\|S_r(f_0)-\sigma_r(f_0)\|=0
\]
if and only if relations (\ref{eqw:2.29})-(\ref{eqw:2.30}) hold.
\end{theorem}
\begin{proof}
Use (\ref{eqw:2.29})--(\ref{eqw:2.30}) and take into account
(\ref{eqw:2.15})--(\ref{eqw:2.16}). \end{proof}
\begin{corollary}\label{wcor2}
\label{C2}
All variables in (\ref{eqw:2.29}) are $o(1)$ for $\delta\le x\le1$.
So we can replace
$\delta$ in (\ref{eqw:2.29}) by $1$ and similarly by $0$ in (\ref{eqw:2.30}).
\end{corollary}

\section{Modification of the criterion}\label{sec:wmodif}

\subsection{Preliminary transformations}
It is difficult to use theorem \ref{Ttrige:whole} directly.
Therefore in this section we
simplify its hypotheses.\hfill
\begin{lemma}
\label{L3.1}
Expressions
$\|\eta_0\|, \ \|\eta_{q}\|,\ \|\zeta_{q-1}\|$
and
$\|\zeta_{n-1}\|$ tend to zero as $r\to\infty$ if and only if
the same is true for
\[
\|I_r^+(\varphi_0)\|,\
\|I_r^+(\varphi_q)\|,\
\|I_r^-(\Psi_{q-1})\|\
\mbox{and}\
\|I_r^-(\Psi_{n-1})\|,
\]
respectively.\hfill
\end{lemma}
\begin{proof}
Recall that\hfill
\begin{equation}
\label{eqw:3.1}
\eta_k
=
 \int\limits_0^1 P_k(x,\sigma,r)\Phi_k(\xi)\,d\xi,\;\;
\zeta_k
=
 \int\limits_0^1 Q_k(x,\sigma,r)\Psi_k(\xi)\,d\xi,
\end{equation}
where\hfill
\begin{equation}
\label{eqw:3.2}
P_k(x,\sigma,r)
=\int_{\Gamma_r}z_k(x,\varrho)y_0(\xi,\varrho)\,d\varrho,
\end{equation}
\begin{equation}
\label{eqw:3.3}
Q_k(x,\sigma,r)
=\int_{\Gamma_r}z_k(x,\varrho)y_{q-1}(\xi,\varrho)\,d\varrho.
\end{equation}
All the four cases in the lemma can be proved in one and the same way.
Therefore we will
consider only the quantity $\|\eta_0\|$.
Formulas (\ref{eqw:3.2}) and (\ref{eqw:1.7})
yield that\hfill
\[
P_0(x,\xi,r)=O(r),\quad 0\le x,\le1.
\]
Replace representation
(\ref{eqw:3.1}) for $\eta_0$ by the integral\hfill
\begin{equation}
\label{eqw:3.4}
\int\limits_{1/r}^1P_0(x,\xi,r)\varphi_0(\xi)\,d\xi
\end{equation}
with an error $o(1)$ as $r\to\infty$.
 Then a direct calculation shows that\hfill
\begin{equation}
\label{eqw:3.5}
P_0(x,\xi,r)=
[exp(ir\varepsilon_1(x+\xi))-exp(ir(x+\xi))]/[i(x+\xi)].
\end{equation}
Clearly\hfill
\[
|exp(ir\varepsilon_1(x+\xi))|=exp(-rh(x+\xi))
\]
with
$h=Im\varepsilon_1>0$
and
$(x+\xi)^{-1}\le r$ because $r^{-1}\le\xi\le1,\ \ x\ge0$.
Since $\Phi_0\in C_0(0,1)$ we have that\hfill
\[
\begin{array}{l}
\|\int\limits_{1/r}^1(x+\xi)^{-1}
exp(ir\varepsilon_1(x+\xi))\Phi_0(\xi)d\xi\|\\
=O(r\int\limits_{1/r}^1\exp(-rh\xi))|\Phi_0(\xi)|d\xi)
=o(1).
\end{array}
\]
Therefore, $\quad\|\eta_0\|=o(1)$ if and only if $\|I_r^+(\Phi_0)\|=o(1)$.
\end{proof}

\subsection{Kernels' calculation}
Let\hfill
\begin{equation}
\label{eqw:3.6}
u_{kp}(x,\xi)
=\varepsilon_kx+\varepsilon_p\xi,\quad
v_{kp}(x,\xi)
=\varepsilon_k(x-1)+
\varepsilon_p\xi.
\end{equation}
Then\hfill
\begin{equation}
\label{eqw:3.7}
P_k(x,\xi,r)
=(iu_{k0})^{-1}[\exp(ir\varepsilon_1u_{k0})-\exp(iru_{k0})]
\end{equation}
for $k=1,\ldots,q-1$;
\begin{equation}
\label{eqw:3.8}
  Q_k(x,\xi,r)=
               (iv_{k,q-1})^{-1}
               [\exp(ir\varepsilon_1u_{k,q-1})-\exp(iru_{k,q-1})]
\end{equation}
for $k=0,\ldots,q-2$.
Analogous formulas hold for \hfill
\[
P_k,\ \ k=q+1,\ldots,n-1;\ \
Q_k,\ \ k=q,\ldots,n-2
\]
with $v_{kp}$ instead of $u_{kp}$ in the exponentials in square
brackets.
Factors $(\ldots)^{-1}$ remain unchanged.\hfill
\begin{lemma}
The following relations are valid\hfill
\begin{align}
\label{eqw:3.9}
\eta_k      &=
i\exp(ir\varepsilon_kx)I_r^+(\Phi_k)+o(1),&
1\le &k \le q-1;\\[10pt]
\label{eqw:3.10}
\eta_k      &=
i\exp(ir\varepsilon_k(x-1))I_r^+(\Phi_k)+o(1),&
q+1\le &k \le n-1;\\[10pt]
\label{eqw:3.11}
\zeta_j &=
(i\varepsilon_1)\exp(ir\varepsilon_{j+1}x)I_r^-(\Psi_j)+o(1),&
0\le &j \le q-2;\\[10pt]
\label{eqw:3.12}
\zeta_j &=
(i\varepsilon_1)\exp(ir\varepsilon_{j+1}(x-1))I_r^-(\Psi_j)+o(1),&
q\le &j\le n-2.
\end{align}
\end{lemma}
\begin{proof}
For the sake of brevity consider
$\eta_k$ for some $k,\ \ 1\le k \leq-1$.
Exponential factor in (\ref{eqw:3.9}) is bounded.
Then, repeating lemma's \ref{L3.1} proof, we get that\hfill
\[
\eta_k
=i\exp(ir\varepsilon_kx)
\int\limits_{1/r}^1(\varepsilon_kx+\xi)^{-1}
\exp(ir\xi)\Phi_k(\xi)d\xi+o(1).
\]
Taking into account an evident relation\hfill
\[
(\varepsilon_kx+\xi)^{-1}-(x+\xi)^{-1}=
\left\{
 \begin{array}{rlc}
  O(x^2/\xi),\quad \xi\ge x,\\
  O(1/x),    \quad \xi<x,
 \end{array}
\right.
\]
we come to (\ref{eqw:3.9}).
Formulas (\ref{eqw:3.10})--(\ref{eqw:3.12}) can be
attained in quite the same way. \end{proof}

\subsection{Equiconvergence with a trigonometric series}
We introduce now a nondegeneracy condition:\hfill
\begin{equation}
\label{eqw:3.13}
span(\Phi_0,\Phi_q)
=span(\Psi_{q-1},\Psi_{n-1})
=span(f_0,f_0^\#).
\end{equation}
Then all the expressions below are $o(1)$,\hfill
\begin{equation}
\label{eqw:3.14}
\|I_r^+(\Phi_0)\|,        \ \
\|I_r^+(\Phi_q)\|,     \ \
\|I_r^-(\Psi_{q-1})\|, \ \
\|I_r^-(\Psi_{n-1})\|=o(1)
\end{equation}
if and only if (\ref{eqw:0.23}) holds. Moreover, (\ref{eqw:3.14})
yields that\hfill
\[
\|\zeta_j\|,\|\eta_k\|=o(1)\ \ \forall j,k.
\]
It remains now to compare lemma
\ref{L3.1} with theorem \ref{Ttrige:whole} and corollary \ref{C2}
and we come to\hfill
\begin{theorem}
\label{Ttrige1:whole}
Consiger a \INDEX{Birkhoff}{G.D.Birkhoff}Birkhoff-reqular operator $L$
defined by the \bvp
(\ref{Chap.intro}.\ref{eqB.1})-(\ref{Chap.intro}.\ref{eqB.2}).
Let $f_0\in C_0(0,1)$ and (\ref{eqw:3.13}) be satisfied.
Then\hfill
\begin{equation}
 \label{eqw:3.14a}
\lim_{r\to\infty}\|S_r(f_0)-\sigma_r(f_0)\|=0
\end{equation}
if and only if (\ref{eqw:0.23}) is true.
\end{theorem}

\subsection{End of theorem's \protect\ref{Te:whole} proof.}
Applying (\ref{eqw:2.7}) to the difference\hfill
$
S_r(f_0,L_1)-S_r(f_0,L_2)
$
we get an analogous
formula with new coefficients\hfill
\[
\beta_{kj}:=\alpha_{kj}(L_1)-\alpha_{kj}(L_2).
\]
It suffices now to repeat  theorem's
\ref{Ttrige1:whole} proof word by word,
replacing
(\ref{eqw:3.13}) by (\ref{eqw:0.19}).
\qed

\section{Functions, satisfying zero-order conditions}
\label{sec:zero}

Let now $n=2q\ge2$ and consider an $n$-th order operator $L\in(R)$.
Let
\begin{equation}
 \label{eqw:zero1}
f\in C[0,1],\ \ f\in clos_{C[0,1]}D_L,
\end{equation}
that is $f$ satisfies
zero-order normalized boundary conditions (if any).
Set
\begin{gather} \label{eqw:f0}
  P(x,f)=f(0)\cdot x + f(1)\cdot (1-x),\quad
  f_0(x):=f(x)-P(x,f),\\
\tilde{f}(x)=\left\{
\begin{array}{ll}
		    f(0), & x<0;\\
                    f(x), & 0\le x \le1;\\
		    f(1), & x>1.
\end{array}
\right.
\end{gather}
A direct calculation (we omit details) shows that
\[
  S_r(P)-\sigma_r(P)=f(0)\cdot\sigma_r(\chi_1)(x) +
f(1)\cdot\sigma_r(\chi_2)(x) +o(1),
\]
where $\chi_1$ and $\chi_2$ are characteristic functions of the intervals
$(-\infty,0)$ and $(1,\infty)$, respectively.

Hence, one obtains a refinement of the theorem \ref{Ttrige1:whole}
assuming validity of (\ref{eqw:zero1}) and replacing (\ref{eqw:3.14a}) with
\begin{equation}
 \label{eqw:3.15}
\lim_{r\to\infty}\|S_r(f)-\sigma_r(\tilde{f})\|=0.
\end{equation}
Quite analogously, theorem \ref{Te:whole} is also true for continuous
functions $f$, subject to (\ref{eqw:0.9}) instead of (\ref{eqw:0.12}).
It is only needed to
replace $f_0$ by $f$ in (\ref{eqw:0.22}).

\section{Order two case}\label{sec:Two}

In this section we present a short proof of
\INDEX{Khro}{A.P.Khromov}A.P.Khromov's
theorem
\ref{Ohrowhole} when \linebreak
$f=f_0\in C_0[0,1]$. The general case is covered
by corollaries \ref{wcor3}--\ref{wcor4}.
\begin{proof}
Repeating considerations from subsections
\ref{subsec:simpl}--\ref{subsec:remain}
word by word
we come to the formula (cf. (\ref{eqw:2.7}) )
\begin{align}
\label{eqw:Two.1}
S_r(f)-\sigma_r(f) &\equiv error\\
+\sum\limits_{k=0}^{1}
\bigl[&\ \{\alpha_{0k}\tau_{r0k(f)}+
\alpha_{1k}\tau_{r1k}(f)\}\nonumber\\
& = error + \sum\limits_{k=0}^{1} J_k, \nonumber
\end{align}
because for $n=2\ \ q-1=0$ and $q=n-1=1$
whereas the sum in (\ref{eqw:Two.1}) contains only two summands instead of
four for $n=2q>2$. Then an analogue of lemma \ref{L2.4} gives
\[
  J_k=\tau_{r0k}(\Phi_k),\ \ k=0,1.
\]
Further, we can not argue as before differentiating the
right-hand side of (\ref{eqw:Two.1}) since now the arc $\Gamma_r$ is a
semicircle with both endpoints on the real axis.
Informally speaking these endpoints both affect an integral over $d\varrho$.
Therefore it is not easy to evaluate the main part of such integral
after differentiation.

Instead, we merely shrink the path of integration to the interval
$[-r,r]$
and deduce that
\begin{align}
J_0 & =  \int_{-r}^r e^{i\varrho x}\,d\varrho \cdot
\int_0^1 \Phi_0(\xi) e^{i\varrho \xi}\,d\xi\\
J_1 & =  \int_{-r}^r e^{i\varrho (1-x)}\,d\varrho \cdot
\int_0^1 \Phi_1(\xi) e^{i\varrho \xi}\,d\xi.
\end{align}
Extending both functions $\Phi_0, \Phi_1$ by zero off $[0,1]$
we readily  obtain
\[
  J_0=2\pi\cdot \sigma_r(\Phi_0)(-x),\ \
  J_1=2\pi\cdot \sigma_r(\Phi_1)(x-1),\ \ 0\leq x \leq 1.
\]
Then
\begin{equation}
\label{eqw:Two.2}
\|  S_r(f)-\sigma_r(f) \| \to 0 \Longleftrightarrow
\|  \sigma_r(\Phi_0)(-x) + \sigma_r(\Phi_1)(x-1) \| \to 0.
\end{equation}
Fix any $\delta,\ 0<\delta<1$ and observe that
\begin{equation}
\label{eqw:Two.3}
\max_{\delta\leq x\leq1}
\|  \sigma_r(\Phi_0)(-x)\| \to 0,\ \
\max_{0\leq x\leq 1-\delta}
\| \sigma_r(\Phi_1)(x-1) \| \to 0
\end{equation}
according to the localization principle for trigonometric series.
At last the proof completes by combining
(\ref{eqw:Two.2}) and (\ref{eqw:Two.3}).

\begin{corollary}\label{wcor3}
Incidentally we also established the following useful relation
for a second order operator $L\in(R)$:
\begin{equation}
\label{eqw:Two.4}
  S_r(f)(x)-\sigma_r(f)(x)=
2\pi\left(
 \sigma_r(\Phi_0)(-x) + \sigma_r(\Phi_1)(x-1)
\right)
+ o(1), \ \ 0\leq x\leq1.
\end{equation}
\end{corollary}
\begin{corollary}\label{wcor4}
Given two second order operators $L_1, L_2\in (R)$,
the following relations
are equivalent:
\begin{align}
 \label{eqw:Two5}
 1. & \lim_{r\to\infty} \|  S_r(f,L_1) -S_r(f,L_2) \| \to 0, \\
 2. & \lim_{r\to\infty} \max_{-\delta\leq x\leq0}
|\sigma_r(\varphi_k)(x)|=0,\ \ k=0,1,
\end{align}
provided that $f$ obeys (\ref{eqw:0.9}).
Here
$\varphi_k$
are defined in (\ref{eqw:0.20}),
functions $f$ and $f_0$ are related by (\ref{eqw:f0}),
 and the proof stems
immediately from  (\ref{eqw:Two.4}).
\end{corollary}

\end{proof}

\section{Odd order operators}
\label{sec:odd}

In this section let $m$ be an odd number, $m=2q+1$,
$L\in(R)$ be an $m$th order differential
operator in $L^2(0,1)$ defined by the
\bvp
(\ref{Chap.intro}.\ref{eqB.1})-(\ref{Chap.intro}.\ref{eqB.2}).
Then theorem \ref{Chap.intro}.\ref{Thsquare} reduces the equiconvergence
problem for operator $L$ to the analogous one for its square.
Moreover, the $\alpha$-numbers for  $L^2$
possess a remarkable property:
\(
  \alpha_{tk}(L^2)=0
\)
\textit{if the indices $t$ and $k$ are of different parity}.

Set  $n=2m$ and (see,
(\ref{Chap.intro}.\ref{eqB.5}) and \ref{Chap.intro}.\ref{eq:total}) )
\begin{gather}
  \delta(L):=\exp(\chi/m)
\Omega(L):=\frac{\theta (b^1,b^0,L)}{\theta (b^0,b^1,L)}
\cdot\frac1{\delta(L)^{q+1}}.
\end{gather}
Then the four most important  $\alpha$-numbers may be written as follows:
\begin{align*}
  \alpha_{00}(L^2)       & = \frac1{2\pi} \delta(L)\cdot\Omega(L),
& \alpha_{m-1,m-1}(L^2)  & = \frac1{2\pi}  \varepsilon_q \cdot\Omega(L), \\
  \alpha_{mm}(L^2)       & = \frac1{2\pi}  \Omega(L)^{-1},
& \alpha_{n-1,n-1}(L^2)  & = \frac1{2\pi}  \cdot\varepsilon_{m-\frac12}
\cdot\delta(L)^{-1}\cdot\Omega(L)^{-1}.
\end{align*}
Further, results of the section \ref{sec:zero} combined with
theorems from our article \cite{min95a} yield
\begin{theorem}
1.Let $L\in (R)$ be an $m$th order differential operator and
$f$ satisfy (\ref{eqw:zero1}).
 Then (\ref{eqw:3.15}) is valid if and only if (\ref{eqw:0.23})
is fulfilled.

2.Let $L_1, L_2\in (R)$ be two $m$th order differential operators and
$f$ satisfy (\ref{eqw:0.9}).
 Then (\ref{eqw:Two5}) is valid if and only if
(\ref{eqw:0.23})
is fulfilled.
\end{theorem}

In this section we improved formulation of the corresponding statements
in \cite{min95a} and checked several misprints there.

\newpage
\setcounter{page}{134}

\newpage
\input{equi.ind}

\begin{thebibliography}{99}
\addcontentsline{toc}{abcd}{\protect\bibname}

\bibitem{AkhGl:II}
{\em Akhiezer N.I.,Glazman I.M.} Theory of linear operators
in hilbert space--- Kharkov: Vitscha schkola, 1978.--- 288 p.

\bibitem{alilnik76}
{\em Alimov Sh.A., {I}l\'{}in V.A., {N}ikishin E.M.} Questions
of convergence of multiple trigonometric serii and spectral
decompositions.I.
 \raz Uspekhi Mat. Nauk--- 1976.--- 31, vyp.6(192).--- P. 28--83


\bibitem{amv83}
{\em Amvrosova O.A.} Eigenvalues asymptotics and
equiconvergence theorems for operators with power
singularities in boundary conditions, \raz In the book:
Functsional'nii analiz. Ulyanovsk--- Ulyanosk gos. ped.
inst. press--- 1983.---  vyp. 21--- P. 3--11

\bibitem{amv84}
{\em \bysame} On one boundary value problem with a
power singularity in boundary condition, \raz Issled. po
sovremennim problemam matematiki. Saratov Univ.--- Saratov
univ. press--- 1984.---  P. 31--37

\bibitem{At:discont}
{\em Atkinson F.V.}
{D}iscrete and {C}ontinuous {B}oundary {P}roblems
--- M.: Mir, 1968.--- 749 p.

\bibitem{Bar:ts}
{\em Bari N.K.} Trigonometric serii
--- M.: Fizmatgiz, 1961.--- 936 p.

\bibitem{baka88}
{\em Baskakov A.G.,Katzaran T.K.} Spectral analysis of
integro-diffe\-ren\-tial ope\-ra\-tors with non\-local boun\-dary
conditions
 \raz Differentsial'nye Uravneniya--- 1988.--- 24, \no 8.---
P. 1424--1433

\bibitem{ben70}
{\em Benzinger H.E.} Green's Function for Ordinary
Differential Operators
 \raz J. Differential Equations--- 1970.--- 7, \no 3.--- P. 478--496

\bibitem{Ber:eigen}
{\em Berezansii Yu.M}
{E}igenfunction expansions for self-adjoint operators
--- Kiev: Naukova Dumka, 1965.--- 798 p.

\bibitem{BeIlNi:int}
{\em Besov O.V., Il'yn V.P., Nikol'skii S.M.}
Integral representation of functions and imbedding theorems
---M.: Nauka, 1975.--- 480 p.

\bibitem{bir08a}
{\em Birkhoff G.D.} On the asymptotic character of the
solutions of certain linear differential equations
containing a parameter
 \raz Trans. Amer. Math. Soc.--- 1908.--- 9.--- P. 219--231

\bibitem{bir08b}
{\em \bysame} Boundary value and expansion problems of
ordinary linear differential equations
 \raz Trans. Amer. Math. Soc.--- 1908.--- 9.---
P. 373--395

\bibitem{DunSchw:II}
{\em Dunford N.,Schwartz J.T.} Linear operators. Part II
(Spectral theory. Self-adjoint operators in hilbert space)
--- M.: Mir, 1966.--- 1063 p.

\bibitem{DunSchw:III}
{\em \bysame,\bysame} Linear operators. Spectral
operators
--- M.: Mir, 1974.--- 661 p.

\bibitem{dya96}
{\em Dyadechko A.V.} To the question of equiconvergence for
matrix differential operators with matrix-dia\-go\-nal
eigenvalue
 \raz Differentsial'nye Uravneniya--- 1996.--- 32, \no 2.---
P. 161--170

\bibitem{ebe64}
{\em Eberhard W.} \language=2
Das asymptotische Verhalten der Greenschen Funktion
irregul\"arer
Eigenvertprobleme mit zerfallenden Randbrdingungen
\language=0
 \raz Math. Z.--- 1964.--- 86.--- P. 45--53

\bibitem{ebfr78}
{\em \bysame,Freiling G.} \language=2
Stone-regul\"are
Eigen\-vert\-prob\-leme
\language=0
 \raz Math. Z.--- 1978.--- 160.--- P.139-161

\bibitem{ef74}   \language=2
{\em \bysame,\bysame} Das Verhalten der Greenschen
Matrix und der Entwicklungen nach Eigenfunktionen
N-regul\"arer Eigenwertprobleme \raz Math. Z.--- 1974
--- 136.--- P. 13--30
\language=0

\bibitem{efs90b}
{\em \bysame,\bysame,Schneider A.} Expansion
theorems for a class of regular indefinite eigenvalue
problems \raz J. Differential Integr. Equat.--- 1990.--- 3,
\no 6.--- P. 1181--1200

\bibitem{efs90a}
{\em \bysame,\bysame,\bysame} On the
distribution of the eigenvalues of a class of regular
indefinite eigenvalue problems \raz J. Differential Integr.
Equat.--- 1990.--- 3, \no 6.--- P. 1167--1179

\bibitem{fie72}
{\em Fiedler H.} \language=2
Zur Regularit\"at selbst\-adjun\-gier\-ter
Rand\-wert\-auf\-gaben
\language=0 \raz Manuscripta Math.--- 1972.--- 7.---
P. 185--196

\bibitem{fre81}
{\em Freiling G.} Irregular Multipoint Eigenvalue Problems
 \raz Math. Methods Appl. Sci.--- 1981.--- 3.---
P. 88--103


\bibitem{frry95}
{\em \bysame, V.Rykhlov}
On a general class of Birkhoff-regular
eigenvalue problems \raz Differential and Integral
equations--- 1995.--- 8,--- \no 8.--- P. 2157--2176

\bibitem{gora91}
{\em Gomilko A.M.,Radzievskii G.V.} Equiconvergence of
series in eigenfunctions of ordinary functional-differential
operators
 \raz Dokl. AN SSSR--- 1991.--- 316,--- \no 2.--- P. 265--270;
Engl. transl. in:
Soviet Math. Dokl.--- 1991.--- 43, \no 1.--- P. 47--52

\bibitem{haa10}  \language=2
{\em Haar A.} Zur {T}heorie der orthogonalen
{F}unktionensysteme.I.
 \raz Math. Ann.--- 1910.--- 69.--- P. 331--371;
II. 1911.--- 71.--- P. 38--53


\language=0

\bibitem{hob08}
{\em Hobson E.W.} On a general convergence theorem, and the
theory of the representation of a function by a series of
normal functions
 \raz Proc. London Math. Soc. (3)--- 1908.--- 6, ser.2.
--- P. 349--395

\bibitem{HNP81}
{\em Hru\v{s}\v{c}\"ev S.V.,{N}ikol\'{}skii N.K.,{P}avlov
B.S.} Unconditional bases of exponentials and of reproducing
kernels, \raz Complex Analysis and Spectral Theory. \\ Lecture
Notes Math.--- Springer-Verlag. --- 1981.--- 864. --- P. 214--335

\bibitem{Il:spectral}
{\em Il\'{}in V.A.} Spectral theory of differential
operators.--- M.: Nauka, 1991.--- 367 p.

\bibitem{il68}
{\em \bysame} Problems of localization and convergence
for Fourier series in fundamental functions of the Laplace
operator
 \raz Uspekhi Mat. Nauk--- 1968.--- 23, \no 2.--- P. 61--120

\bibitem{il76a}
{\em \bysame} Necessary and sufficient conditions of
basicity of a subsystem of eigen- and associated functions
of the M.V.Keldysh' pencil of ordinary differential
operators
 \raz Dokl. Akad. Nauk SSSR--- 1976.--- 227, \no 4.---
P. 796--799

\bibitem{il77}
{\em \bysame} On convergence of eigenfunction
expansions in the points of the coefficients' discontinuity
of a differential operator
 \raz Mat. Zametki--- 1977.--- 22, \no 5.--- P. 679--698

\bibitem{il80a}
{\em \bysame} Necessary and sufficient conditions of
basicity
and equiconvergence with a trigonometric series of spectral
decompositions.I
 \raz Differentsial'nye Uravneniya--- 1980.--- 16, \no 5.---
P. 771--794;II. 1980.--- 16, \no 6.--- P. 980--1009


\bibitem{il82a}
{\em \bysame} On sharp in order relations between norm
estimates for eigen- and associated functions of a second
order elliptic operator
 \raz Differentsial'nye Uravneniya--- 1982.--- 18, \no 1.---
P. 30--37

\bibitem{il83a}
{\em \bysame} Necessary and sufficient conditions of
basicity
in $L_p$ and of equiconvergence with a trigonometric series
of spectral expansions and decompositions in exponential
series
 \raz Dokl. Akad. Nauk SSSR--- 1983.--- 273, \no 4.---
P. 789--793

\bibitem{il95a}
{\em \bysame} Uniform equiconvergence on the whole
line
$\mathbb{R}$
with Fourier integral of the spectral decomposition,
corresponding to a self-adjoint extension of Schr\"odinger
operator with uniformly locally summable potential
 \raz Differentsial'nye Uravneniya
--- 1995.--- 31, \no 12.--- P. 1947--1956


\bibitem{ilan95}
{\em \bysame,Antoniu I.} On uniform equiconvergence
with Fourier integral
on the whole line
$\mathbb{R}$
 of the spectral decomposition of arbitrary function from
the class
$L_p(\mathbb{R})$,
 corresponding to a self-adjoint extension of the Hill
operator
 \raz Differentsial'nye Uravneniya
--- 1995.--- 31, \no 8.--- P. 1310--1322


\bibitem{ilan96}
{\em \bysame,\bysame} On spectral decompositions
corresponding to a liouvillian, generated by the
 Schr\"odinger operator with uniformly locally summable
potential
 \raz Differentsial'nye Uravneniya
--- 1996.--- 32, \no 4.--- P. 435--440


\bibitem{iljoo79}
{\em \bysame,Jo\'{o} I.} Uniform eigenfunctions'
estimate and estimate from above of the number of
eigenvalues of the Sturm-Liouville operator with a class
$L^p$ potential
 \raz Differentsial'nye Uravneniya
--- 1979.--- 15, \no 7.--- P. 1164--1174


\bibitem{ilkri95}
{\em \bysame,Kritskov L.V.} Uniform estimate on the
whole line  of generalized eigenfunctions of one-dimensional
SChr\"odinger operator with a uniformly locally summable
potential
 \raz Differentsial'nye Uravneniya
--- 1995.--- 31, \no 8.--- P. 1323--1329


\bibitem{ilmo83}
{\em \bysame,{M}oiseev E.I.} Sharp in order maximum
moduli estimates of eigen- and associated functions of
elliptic operators
 \raz Mat. Zametki
--- 1983.--- 34, \no 5.--- P. 683--692


\bibitem{ilmo92}
{\em \bysame,\bysame} On the systems consisting of
subsets of root functions of two distinct boundary value
problems
 \raz Trudy Mat. Inst. RAN.--- 1992.--- 201.--- P. 219--230

\bibitem{ima93}
{\em Imamberdiev V.I.} Spectral function's asymptotic of
ordinary odd order differential operator in a general case
 \raz Uspekhi Mat. Nauk--- 1993.--- 48, \no 1.--- P. 165--166


\bibitem{joo84}
{\em Jo\'{o} I.} Remarks to a paper of V. Komornik
 \raz Acta Sci. Math. (Szeged)
--- 1984.--- 47, --- P. 201--204


\bibitem{kab88}
{\em Kabanov S.N.} Equiconvergence theorem for a $n$-th
order differentiation operator with boundary
conditions generated by linear functionals, \raz Mathematika
i ee prilojeniya.Saratov univ.--- Saratov univ. press.
 --- 1988.--- P. 4--6

\bibitem{kab90a}
{\em \bysame} Equiconvergence theorem for differential
operators with a general form boundary condition, \raz
Teoriya functsii i pribl. (Trudy  4 Sarat. zimn. schkoly)---
Saratov univ. press,
--- 1990.--- Part II--- P. 108--110

\bibitem{kab90b}
{\em \bysame} Equiconvergence theorem for one
integro-differential operator
 \raz Sarat. univ.--- Saratov, 1990.--- 21 p. --- Bibliogr.  items. 5
---Rus.--- Dep. VINITI 30.07.90, ü 4312-B90

\bibitem{Kah:abs}
{\em Kahane J.P.} S\'eries de {F}ourier ab\-so\-lu\-ment
con\-ver\-gentes.--- M: Mir, 1970.--- 206 p.

\bibitem{kats56}
{\em Kats I.S.} On integral representations of analytic
functions mapping an upper half-plane to its part
 \raz Uspekhi Mat. Nauk
--- 1956.--- 11, vyp.3(69)--- P. 139--144


\bibitem{kalu94}
{\em Kaufmann F.J., W.J.Luter
} Degree of convergence of Birkhoff serii,
direct and inverse theorems \raz J. Math. Anal. Appl.
---1994.--- 1.--- P. 156--168


\bibitem{kau89}
{\em \bysame}
\newblock On the degree of convergence of {B}irkhoff's series
for functions of
  bounded variation \raz Analysis
--- 1989.--- 9.--- P. 303--315


\bibitem{kau:phd}\language=2
{\em \bysame}
\newblock {Abgeleitete {B}irkhoff-{R}eihen bei {R}andeigenwertproblemen zu
  ${N}(y)=\lambda {P}(y)$ mit $\lambda$-abh\"angigen {R}andbedingungen}.
\newblock PhD thesis, Aachen, 1989


\bibitem{ker86a} 
{\em Kerimov N.B.} Asymptotical formulas for eigen- and
associated functions of ordinary differential operators
 \raz Moscow univ.--- Moscow, 1986.--- 19 p. --- Bibliogr. 3 items.---
Rus.--- Dep. VINITI 25.12.1986, ü 663-B86

\bibitem{ker86b}  
{\em \bysame} Some properties of eigen- and associated
functions of ordinary differential operators
 \raz Dokl. Akad. Nauk SSSR
--- 1986.--- 291, \no 5.--- P. 1054--1056


\bibitem{hro62}
{\em Khromov A.P.} Eigenfunction expansion of ordinary
linear differential operators in a finite interval
 \raz Dokl. Akad. Nauk SSSR
--- 1962.--- 146, \no 6.--- P. 1294--1297


\bibitem{hro66}
{\em \bysame} Eigenfunction expansion of ordinary
linear differential
operators with decomposing boundary conditions
 \raz Math. USSR-Sb.
--- 1966.--- 70, \no 3.--- P. 310--329


\bibitem{hro75}
{\em \bysame} On equiconvergence for eigenfunction
expansion associated with second order differential operators,
\raz Differentsial'nye Uravneniya i Vitchyslistel'naya
Mathematika. Saratov univ.--- Saratov univ. press,
--- 1975.--- 5, Part II.--- P. 3--20

\bibitem{hro76}
{\em \bysame} Differential operator with irregular
decomposing boundary conditions
 \raz Mat. Zametki
--- 1976.--- 19, \no 5.--- P. 763--772


\bibitem{hro81}
{\em \bysame} Equiconvergence theorems for
integro--differential
and integral operators
 \raz Math. USSR-Sb.
--- 1981.--- 114(156), \no 3.--- P. 378--405


\bibitem{hro95a}
{\em \bysame} Spectral analysis of differential
operators in a finite interval
 \raz Differentsial'nye Uravneniya--- 1995.--- 31, \no 10.---
P. 1691--1696


\bibitem{hro95b}
{\em \bysame} Equiconvergence of spectral
decompositions, \raz Teoriya functsii i pribl.(Trudy  Sarat.
zimn. schkoly)--- Saratov univ. press.
 --- 1995. Part 1.--- P. 86--96


\bibitem{hro:phd}
{\em \bysame} Expansion in eigenfunctions of ordinary
linear differential operators in a finite interval, Phd
Thesis, Saratov, --- 1963

\bibitem{kogro75}
{\em Kogan V.I., Rofe-Beketov F.S.} On Square-integrable
Solutions of Symmetric Systems of Differential Equations of
Arbitrary Order \raz Proc. Roy. Soc. Edinburgh Sect. A
---1974/75.--- 74A, \no 1.--- P. 5--40


\bibitem{kom84a}
{\em Komornik V.} Upper estimates for eigenfunctions \raz
Ann. Univ. Sci. Budapest. E{\"o}tv{\"o}s Sect. Math.
---1984.--- 27. --- P. 125--135


\bibitem{kom84b}
{\em \bysame} Generalisation of a theorem of \Io
 \raz Ann. Univ. Sci. Budapest. E{\"o}tv{\"o}s Sect.
Math.
--- 1984(1985).--- 27. --- P. 59--64


\bibitem{kom85a}
{\em \bysame} Some new estimates for the eigenfunctions
of
higher order \raz Acta Math. Hungar.
--- 1985.--- 45, \no 3-4.--- P. 451--457


\bibitem{kom85b}
{\em \bysame} Lower estimates for the eigenfunctions
\raz Acta Math. Hungar.
--- 1985.--- 45, \no 1-2.--- P. 189--193


\bibitem{kom85c}
{\em \bysame} Local upper estimates for the
eigenfunctions of a linear differential operator
 \raz Acta Sci. Math. (Szeged)
--- 1985.--- 48, \no 1-4.---P. 243--256


\bibitem{kom86a}
{\em \bysame} The asymptotic behavior of the
eigenfunctions
of higher order of a linear differential operator \raz
Studia Sci. Math.
--- 1986.--- 5. --- P. 1075--1077


\bibitem{kos68}
{\em Kostuchenko A.G.} Asymptotic behaviour of the spectral
function of self-adjoint elliptic operators, \raz Fourth
summer mathematical school
--- 1968.  --- P. 42--117

\bibitem{kos:dr}
{\em \bysame} On some spectral properties of
differential operators. Dis. ... dokt. phys.-mat. nauk.,
M.: MGU.--- 1966.

\bibitem{kra68}
{\em Krall A.M.} Differential operators and their adjoints
under integral and multiple point boundary conditions \raz
J. Differential Equations
--- 1968.--- 4.--- P. 327--336


\bibitem{kra75}
{\em \bysame} The development of general differential and
general differential boundary systems
 \raz Rocky Mountain J. Math.
--- 1975.--- 5, \no 4.--- P. 493--542


\bibitem{kri92}
{\em Kritskov L.V.} Representation and estimates of root
functions of a singular differential operator in an
interval.I
 \raz Differentsial'nye Uravneniya
--- 1992.--- 28, \no 8.--- P. 2294--2305;
II.--- 1993.--- 29, \no 1.--- P. 64--73




\bibitem{nkup67a}
{\em Kuptsov N.P.} Equiconvergence theorem for Fourier
expansions in
Banach spaces
 \raz Mat. Zametki
--- 1967.--- 1, \no 4.--- P. 469--474


\bibitem{nkup67b}
{\em \bysame} Localization of equiconvergence theorems
 \raz Math. USSR-Sb.
--- 1967.--- 74(116), \no 4.--- P. 554--564


\bibitem{kur96}
{\em Kurkina A.B.} Uniform component equiconvergence on the
whole line
with Fourier integral of the spectral decomposition,
corresponding to the Schr\"odinger operator with a matrix
potential satisfying Kato condition
 \raz Differentsial'nye Uravneniya
--- 1996.--- 32, \no 6.--- P. 759--768


\bibitem{lan26}
{\em Langer R.} On the theory of integral operators with
discontinuous kernels
 \raz Trans. Amer. Math. Soc.
--- 1926.--- 28, \no 4.--- P. 585--639


\bibitem{lev53}
{\em Levitan B.M.} On asymptotic behaviour of spectral
function and eigenfunction expansion of second order
self-adjoint differential equation.I
 \raz Izv. Akad. Nauk SSSR Ser. Mat.
--- 1953.--- 17, \no 4--- P. 331--364;
II --- 1955.--- 19, \no 1--- P. 33--58




\bibitem{lom85}
{\em Lomov I.S.} Estimates of eigenfunctions and generalized
eigenfunctions of ordinary differential operators
 \raz Differentsial'nye Uravneniya
--- 1985.--- 21, \no 5.--- P. 903--906


\bibitem{lom95a}
{\em \bysame} On function approximation on an interval by
spectral expansions of Schr\"odinger operator
 \raz Vestnik Moskov. Univ. Ser. I Mat. Mekh.
--- 1995.--- 1, \no 4.--- P. 43--54


\bibitem{mar55}
{\em Marchenko V.A.} Tauberian type theorems in spectral
analysis of differential operators
 \raz Izv. Akad. Nauk SSSR Ser. Mat.
--- 1955.--- 19, \no 6.--- P. 381--422


\bibitem{min77a}
{\em Minkin A.M.} Regularity of self-adjoint boundary
conditions, \raz Mat. Zametki
---  1977.---22, \no 6.--- P. 835--846

\bibitem{min77b}
{\em \bysame} Equiconvergence theorem for a normal
integral operator with a Green function-type kernel, \raz
Vychisl. Metody i Programmirovanie--- Saratov univ. press.
--- 1977.--- vyp. 1--- P. 181--190

\bibitem{min80a}
{\em \bysame} Equiconvergence theorem for
expansions associated with a generalized spectral function
of a symmetric differential operator and in Fourier integral,
\raz In the book: Functsional'nii analiz. Ulyanovsk---
Ulyanosk gos. ped. inst. press. Spectral theory.---
1980.---  14.--- P. 109--112

\bibitem{min80b}
{\em \bysame} Equiconvergence theorem for singular
self-adjoint  differential operators, \raz All-Union symposium
on function approximation in the complex domain.---
Ufa.--- 1980.--- P. 95--96

\bibitem{min80c}
{\em \bysame} Localization principle for series
in eigenfunctions of ordinary differential operators, \raz
Differentsial'nye Uravneniya i teoriya functsii--- Saratov
univ. press.--- 1980.--- 3.--- P. 68--80


\bibitem{min90}  
{\em \bysame}
Eigenfunction expansion for one class of nonsmooth
differential operators \raz Differentsial'nye Uravneniya
--- 1990.--- 26, \no 2.--- P. 356--358; Manuscript completely
deposited at VINITI 11.08.89, ü 5407-'89 by the editorial board
of the journal. Minsk. 1989.--- 54 p. ---
Bibliogr.  29 items---Rus.


\bibitem{min91}
{\em \bysame} Reflection of exponents and unconditional
bases from exponentials
 \raz Algebra i analyz
--- 1991.--- 3, \no 5.--- P. 110--135;
Engl. transl. in St. Petersburg Math. J.
--- 3, \no 5.--- P. 1043--1068



\bibitem{min93a}
{\em \bysame} On the class of regular
boundary conditions
 \raz Re\-sults in Mathematics
--- 1993.--- 24, $n^0$ 3/4--- P. 274--279


\bibitem{min93b}
{\em \bysame} {A}lmost orthogonality of {B}irkhoff's
{S}olutions
 \raz Re\-sults in Math\-e\-matics
--- 1993.--- 24, $n^0$ 3/4.--- P. 280--287


\bibitem{min95a}
{\em \bysame} Odd and Even cases of Birkhoff-regularity
 \raz Math. Nachr.
--- 1995.--- 174.--- P. 219--230



\bibitem{min:phd}
{\em \bysame} Equiconvergence theorems for
differential  operators, Phd Thesis, Saratov, SGU.--- 1982

\bibitem{mishu91}
{\em \bysame,{S}huster L.} Spectrum distribution and
convergence of spectral expansions for a Schr\"odinger
operator
 \raz Differentsial'nye Uravneniya
--- 1991.--- 27, \no 10.--- P. 1778--1789


\bibitem{moi80}
{\em Moiseev E.I.} Asymptotical mean value formulas for
regular solution of differential equation
 \raz Differentsial'nye Uravneniya
--- 1980.--- 16, \no 5.--- P. 827--844


\bibitem{moi84}
{\em \bysame} On basicity of sine and cosine systems
 \raz Dokl. Akad. Nauk SSSR
--- 1984.--- 275, \no 4.--- P. 794--798


\bibitem{moi87}
{\em \bysame} On basicity of some sine system
 \raz Differentsial'nye Uravneniya
--- 1987.--- 23, \no 1.--- P.177--179


\bibitem{Nai:ldo}
{\em Naimark M.A.}
{L}inear differential operators
.--- M.: Nauka, 1969.--- 528 p.

\bibitem{osk73}
{\em Os'kina G.P.} Asymptotic formulas for partial sums of
Fourier eigenfunction series of ordinary
differential operators, \raz Issled. po different. uravn. i
theorii functsii--- Saratov univ. press.---
1973.--- P. 40--54

\bibitem{pal72}
{\em Pal'tsev B.V.} Eigenfunction expansion of integral
convolution operators in a finite interval with rational
kernel Fourier transform
 \raz Izv. Akad. Nauk SSSR Ser. Mat.
--- 1972.--- 36, \no 3.--- P. 591--634


\bibitem{pav73}
{\em Pavlov B.S.} Spectral analysis of a differentiation
operator with a "smeared" boundary condition
 \raz Problemi matem. fiziki
--- 1973.--- 6.--- P. 101--119


\bibitem{rad94a}
{\em Radzievskii G.V.} The rate of convergence of
decompositions of ordinary functional-differential operators
by eigenfunctions.
In a book:
Some problems of the modern theory of differential
equations.
Preprint 94.29. Ukraianian national academy of sciences.
Institute of mathematics. Kiev. 1994. \raz
--- 1994.--- P. 14--27


\bibitem{rad95a}
{\em \bysame} Boundary value problems and
associated moduli of continuity
 \raz Funktsional. Anal. i Prilozhen.
--- 1995.--- 29, \no 3.--- P. 87--90


\bibitem{rad96a}
{\em \bysame} Eigenvalues' asymptotics of a regular
boundary value problem
 \raz Ukrain. Mat. Zh.
--- 1996.--- 48, \no 4.--- P. 483--519


\bibitem{ryh77}
{\em Rykhlov V.S.} Eigenfunction expansion for one class
of quasidifferential operators, \raz Differentsial'nye
Uravneniya i teoriya functsii--- Saratov univ. press.
--- 1977.--- vyp. 1--- P. 151--169

\bibitem{ryh83}
{\em \bysame} Asymptotics of a system of solutions of
a quasidifferential operator, \raz Differentsial'nye
Uravneniya i teoriya functsii--- Saratov univ. press.
--- 1983.--- vyp. 5--- P. 51--59

\bibitem{ryh84}
{\em \bysame} On the rate of equiconvergence for
differential operators with a nonzero coefficient by the
$n-1$th derivative
 \raz Dokl. Akad. Nauk SSSR
--- 1984.--- 279, \no 5.--- P. 1053--1056


\bibitem{ryh96a}
{\em \bysame} Equiconvergence rate in terms of general
moduli of continuity for differential operators
 \raz Results in Mathematics
--- 1996.--- 29.--- P. 153--168


\bibitem{ryh:phd}
{\em \bysame} Eigenfunction expansions for
quasidifferential and integral operators, Phd Thesis,
Saratov, SGU. --- 1981

\bibitem{sal68}
{\em Salaff S.} {R}egular {B}oundary {C}onditions for
{O}rdinary {D}ifferential {O}perators
 \raz Trans. Amer. Math. Soc.
--- 1968.--- 134, \no 2.--- P. 355--373


\bibitem{schaf60}
{\em Sch\"afke F.} \language=2
Rei\-hen\-ent\-wick\-lun\-gen ana\-ly\-ti\-schen
Funk\-tio\-nen nach
Bi\-or\-tho\-go\-nal\-sys\-te\-men spe\-zi\-el\-ler
Funk\-tio\-nen.I
\language=0 \raz Math. Z.
--- 1960.--- 74, \no 4.--- P. 436--470;
II---1961.--- 75, \no 1.--- P. 154--191;
III---1963.--- 80, \no 4.--- P. 400--442






\bibitem{bschu79a}
{\em Schultze B.} On the definition of Stone-Regularity \raz
J. Dif\-feren\-tial Equ\-a\-tions
--- 1979.--- 31.--- P. 224--229


\bibitem{bschu79b}
{\em \bysame} Strongly irregular boundary value problems
 \raz Proc. Roy. Soc. Edinburgh Sect. A
--- 1979.--- 82A.--- P. 291--303


\bibitem{sed75}
{\em Sedletskii A.M.} On equiconvergence and equisummability
of nonharmonic Fourier expansions with ordinary
trigonometric series
 \raz Mat. Zametki
--- 1975.--- 18, \no 1.--- P. 9--17


\bibitem{sed82}
{\em \bysame} Biorthogonal expansions of functions
in exponential serii on the intervals of the real axis
 \raz Uspekhi Mat. Nauk
--- 1982.--- 37, \no 5(227)--- P. 51--95


\bibitem{sed91}
{\em \bysame} On uniform convergence of nonharmonic
Fourier series
 \raz Proc. Steklov Inst. Math.
--- 1991.--- 200, \no --- P. 299--309


\bibitem{sed94}
{\em \bysame} Expansion in eigenfunctions  of a
differentiation operator with a smeared boundary condition
 \raz Differentsial'nye Uravneniya
--- 1994.--- 30, \no 1.--- P. 70--76


\bibitem{sed95b}
{\em \bysame} Approximative properties of
exponential systems in $L^p(a,b)$
 \raz Differentsial'nye Uravneniya
--- 1995.--- 31, \no 10.--- P. 1675--1681


\bibitem{Shi:sec}
{\em Shilov G.E.} Mathematical analysis. Second special
course
.--- M.: Nauka, 1965.--- 328 p.

\bibitem{shi72}
{\em Shilovskaya O.K.} Expansion in eigenfunctions of a
second order differential operator in the case of
irregular  boundary conditions, \raz Differential and Integral
equations--- Saratov univ. press.--- 1972.--- P. 53--79

\bibitem{shk83}
{\em Shkalikov A.A.} Boundary value problems for ordinary
differential equations with a parameter in the boundary
conditions
 \raz Trudy Sem. Petrovsk.
--- 1983.---9.--- P. 190--229. Engl. transl. in:
J. Soviet Math.--- 1986.--- 33.--- P. 1311--1342




\bibitem{shktr94}
{\em \bysame,Tretter C.} Kamke Problems. Properties
of Eigenfunctions
 \raz Math. Nachr.
--- 1994.--- 170.--- P. 251--275


\bibitem{sht57}
{\em Shtraus A.V.} On generalized resolvents and spectral
functions of even order differential operators
 \raz Izv. Akad. Nauk SSSR Ser. Mat.
--- 1957.--- 21, \no 6.--- P. 785--808


\bibitem{sht65}
{\em \bysame} On extensions of a symmetric operator
depending on a parameter
 \raz Izv. Akad. Nauk SSSR Ser. Mat.
--- 1965.--- 29, \no 6.--- P. 1389--1416


\bibitem{shu84}
{\em Shuster L.A.} Uniform eigenfunctions estimates for one
differential operator
 \raz Vestnik AN Kazakh.SSR.
---Alma-ata, 1984.--- 29 p. --- Bibliogr. 10 items.--- Rus.---
Dep. VINITI 27.04.1984, \no 3303-84

\bibitem{ste07}
{\em Stekloff V.A.} Sur les expressions asymptotiques de
certains fonctions d\'efinies par des \'equations
diff\'er\'entielles lin\'eaieres de deuxi\'eme ordre, et
leurs applications au probl$\grave{e}$me
du developpement d'une fonction arbitraires en s\'eries
proc\'edant suivantes dites fonctions
 \raz Kharkov. Soobtscheniya matem. obtschesstva
---1907--1909.--- 10, \no (2-6).--- P. 97--199


\bibitem{ste10}
{\em \bysame} Solution g\'en\'erale du
probl$\grave{e}$me
de developpement d'une fonction arbitraire en s\'eries
suivant les fonctions fondamentales de
Sturm-Liouville
 \raz Rend. Acad. Lincei
--- 1910.--- 19.--- P. 490--496


\bibitem{sto26}
{\em Stone M.H.} A comparison of the series of Fourier and
Birkhoff
 \raz Trans. Amer. Math. Soc.
--- 1926.--- 28.--- P. 695--761


\bibitem{sto27}
{\em \bysame} Irregular differential systems of order two
and the related expansion problems
 \raz Trans. Amer. Math. Soc.
--- 1927.--- 29.--- P. 23--53


\bibitem{such79}
{\em Suchkov M.V.} On uniform equiconvergence on any compact
of
an integral Fourier expansion and of spectral decomposition,
corresponding to a nonself-adjoint differential operator on
the half-axis
 \raz Differentsial'nye Uravneniya
--- 1979.--- 15, \no 12.--- P. 2161--2167


\bibitem{Tam:expan}
{\em Tamarkin J.D.} On some general problems of the theory
of ordinary linear differential operators and on expansion
of arbitrary function into serii
--- Petrograd. 1917.--- 308 p.

\bibitem{tam12}
{\em \bysame} Sur quelques points de la th\'e\-o\-rie des
\'equa\-tions
diff\'e\-ren\-tielles lin\'eaires ordinaires et sur la
g\'en\'eralisation de la s\'erie de Fourier
 \raz Rend. Circ. Mat. Palermo (2)
--- 1912.--- 34.--- P. 345--382


\bibitem{tam28}
{\em \bysame} Some general problems of the theory of
linear differential equations and expansions of an arbitrary
functions in series of fundamental functions
 \raz Math. Z.
--- 1928.--- 27, \no 1.--- P. 1--54


\bibitem{wer:phd}
{\em Wermuth E.} \language=2
Konvergenz\-unter\-suchungen bei
{E}igen\-funktions\-ent\-wick\-lungen zu {R}and\-eigen\-wert\-problemen
$n$-ter {O}rdnung mit
parameter\-abh\"an\-gigen {R}and\-bedin\-gungen. Thesis, Aachen, RWTH,
Mathematisch-Natur\-wissen\-schaft. fak. --- 1984
\language=0

\bibitem{wer89}
{\em \bysame} A generalization of Lebesgue's convergence
criterion for Fourier series
 \raz Results in Mathematics
--- 1989.--- 15.--- P. 186--195


\bibitem{wol83}
{\em Wolter M.} \language=2
Das asymptotische Verhalten der Greenschen Funktion
$N$-irregul\"arer
Eigenwertprobleme mit zerfallenden Randbedingungen
\language=0
 \raz Math. Methods Appl. Sci.
--- 1983.--- 5.--- P. 331--345


\bibitem{Zyg:ts}
{\em Zygmund A.} Trigonometric series. I,II.
--- M.: Mir, 1965.--- 615 p.;537 p.

\end{thebibliography}
\end{document}